\documentclass[11pt,fleqn,a4paper]{article}
\usepackage{amsfonts, amsmath, amssymb, graphicx, latexsym, curves, ebezier}

    \advance\textheight by \topskip

    \numberwithin{equation}{section}

    \newtheorem{theorem}{Theorem}[section]
    \newtheorem{lemma}[theorem]{Lemma}
    
    \newtheorem{proposition}[theorem]{Proposition}
    
    \newtheorem{Definition}[theorem]{Definition}
    
    \newtheorem{Remark}[theorem]{Remark}
    \newenvironment{remark}{\begin{Remark}\rm}{\end{Remark}}
    \newtheorem{Example}[theorem]{Example}
    
    \newenvironment{proof}
        { \rm \trivlist \item[\hskip \labelsep{\bf Proof. }] }
        { \hspace*{\fill}$\Box$\endtrivlist }
    \newenvironment{varproof}
        { \rm \trivlist \item[\hskip \labelsep{\bf Proof}] }
        { \hspace*{\fill}$\Box$\endtrivlist }

    \newcommand{\Ai}{\operatorname{Ai}}
    \newcommand{\I}{\operatorname{I}}
    \newcommand{\II}{\operatorname{II}}
    \newcommand{\III}{\operatorname{III}}
    \newcommand{\IV}{\operatorname{IV}}
    \newcommand{\supp}{\operatorname{supp}}
    \newcommand{\PVint}
        {
            \mathop{
                \setbox0\hbox{$\displaystyle\intop$}
                \hskip0.2\wd0
                \vcenter{\hrule width0.6\wd0height0.5pt depth0.5pt}
                \hskip-0.8\wd0
            }
            \mskip-\thinmuskip\intop\nolimits
        }
    \renewcommand{\Re}{\operatorname{Re}}
    \renewcommand{\Im}{\operatorname{Im}}
    \def\bigO{{\cal O}}

    \renewcommand{\qbeziermax}{1000}

\begin{document}

\begin{center} \Large\bf
    Strong asymptotics of Laguerre-type orthogonal polynomials and
    applications in random matrix theory
\end{center} \

\begin{center}
    \large M. Vanlessen\footnote{
        Postdoctoral Fellow of the Fund for Scientific Research
        -- Flanders (Belgium). This work was done while visiting
        the Department of Mathematics of the Ruhr Universit\"at
        Bochum and is supported in part by the SFB/TR12 of the Deutsche Forschungsgemeinschaft.} \\[1ex]
    \normalsize \em
        Department of Mathematics, Katholieke Universiteit Leuven, \\
        Celestijnenlaan 200 B, 3001 Leuven, Belgium \\
            and \\
        Fakult\"at f\"ur Mathematik, Ruhr Universit\"at Bochum \\
        Universit\"atsstrasse 150, 44801 Bochum, Germany \\[1ex]
    \rm maarten.vanlessen@wis.kuleuven.be
\end{center}\ \\[1ex]

\begin{abstract}
    We consider polynomials orthogonal on $[0,\infty)$ with
    respect to Laguerre-type weights $w(x)=x^\alpha e^{-Q(x)}$,
    where $\alpha>-1$ and where $Q$ denotes a polynomial with
    positive leading coefficient. The main purpose of this paper
    is to determine Plancherel-Rotach type asymptotics in the
    entire complex plane for the orthonormal polynomials with
    respect to $w$, as well as asymptotics of the corresponding
    recurrence coefficients and of the leading coefficients of
    the orthonormal polynomials. As an application we will use
    these asymptotics to prove universality results in random
    matrix theory.

    We will prove our results by using the characterization of
    orthogonal polynomials via a $2\times 2$ matrix valued
    Riemann-Hilbert problem, due to Fokas, Its and Kitaev,
    together with an application of the Deift-Zhou steepest
    descent method to analyze the Riemann-Hilbert problem
    asymptotically.
\end{abstract}

\section{Introduction}

We consider Laguerre-type weights,
\begin{equation}\label{weight}
    w(x)=x^\alpha e^{-Q(x)}, \qquad\mbox{for $x\in[0,\infty)$,}
\end{equation}
where $\alpha>-1$, and where
\begin{equation}\label{definition: Q}
    Q(x)=\sum_{k=0}^m q_k x^k,\qquad q_m>0,
\end{equation}
denotes a polynomial of degree $m$ with positive leading
coefficient. In case $Q(x)=x$, the weight $w$ is the classical
Laguerre weight. Since all the moments of $w$ exist we have a
sequence of orthogonal polynomials. We use
$p_n(x)=p_n(x;w)=\gamma_n x^n+\cdots$ with leading coefficient
$\gamma_n>0$, to denote the $n$-th degree orthonormal polynomial
with respect to $w$, i.e.
\[
    \int_0^\infty p_k(x) p_n(x) w(x)dx=\delta_{k,n},
        \qquad\mbox{for $k,n\in\mathbb{N}$.}
\]
The orthonormal polynomials $p_n$ satisfy a three-term recurrence
relation of the form,
\begin{equation}\label{three-term recurrence relation}
    xp_n(x)= b_n p_{n+1}(x) + a_n p_n(x) + b_{n-1} p_{n-1}(x),
        \qquad \mbox{$n\in\mathbb{N}$,}
\end{equation}
and the coefficients $a_n$ and $b_{n-1}$ are called the recurrence
coefficients.

The goal of this paper is to determine the asymptotic behavior (as
$n\to\infty$) of the recurrence coefficients $a_n, b_{n-1}$ and of
the leading coefficient $\gamma_n$, and to determine
Plancherel-Rotach type asymptotics for the orthonormal polynomials
$p_n$ in the entire complex plane. In addition, we use these
asymptotics to prove universality results for Laguerre-type
unitary ensembles.

\medskip

The weights considered in this paper fall into the class of Freud weights on
$[0,\infty)$. These weights are also closely related to symmetric Freud weights
on $\mathbb{R}$ by a simple transformation so that our results can also be
applied to Freud weights of the form $|x|^\gamma e^{-Q(x^2)}$ on $\mathbb{R}$
with $\gamma>-1$. There is a vast literature on asymptotic questions with
respect to Freud weights and we will now briefly summarize the results which
are already known with respect to weights (\ref{weight}).

The leading order asymptotics of the recurrence coefficients associated with
symmetric weights $|x|^\gamma e^{-Q(x^2)}$ on $\mathbb{R}$, where
$Q$ is a polynomial with non-negative coefficients, have been
derived by Damelin \cite{Damelin}. This result can then be applied
to get the leading order asymptotics of the recurrence coefficients associated
to the corresponding weights (\ref{weight}).

Full asymptotic description of the orthogonal polynomials and
related quantities are known in the following special cases.
First, the classical Laguerre case, i.e.\ $Q(x)=x$ and
$\alpha>-1$, is described in \cite{Szego}. Next, for $\alpha=\pm
1/2$ and $Q$ a polynomial, the asymptotics are known from results
of Deift, Kriecherbauer, McLaughlin, Venakides and Zhou in
\cite{DKMVZ1}. Indeed, the weights
$|x|^{\pm 1/2}e^{-Q(x)}$ on $[0,\infty)$ are related to the
symmetric weights $e^{-Q(x^2)}$ on $\mathbb{R}$, which belong to
the class of weights considered in \cite{DKMVZ1}. See also the
references therein for older results. Finally, in the special case
that $Q(x)=q_m x^m$ is a monomial, one can verify that the weights
$x^\alpha e^{-q_m x^m}$ on $[0,\infty)$ relate to the varying
weights $x^\alpha e^{-n q_m x^{2m}}$ on $\mathbb{R}$, which fall
into the class of weights considered by the author together with
Arno Kuijlaars in \cite{KV2}. In that paper, the main goal was to
prove universality of spectral correlations of the associated random
matrix ensembles at the hard edge $\lambda=0$. Explicit formulae for
orthogonal polynomials were not derived in \cite{KV2} but can be
deduced from their proof.

\medskip

One strong motivation for deriving the asymptotics of the
orthogonal polynomials comes from random matrix theory. Indeed,
the Laguerre-type weights (\ref{weight}) are related to random
matrix ensembles used in the statistical description of mesoscopic
quantum systems with disorder. Depending on the symmetries of the
physical system one distinguishes ten different classes of
ensembles, see for example \cite{AltlandZirnbauer,Ivanov,TBFM} and
references therein. Besides the three standard classes
(orthogonal, unitary and symplectic) there are seven novel
classes, which have been introduced by physicists during the past
ten years. In the Gaussian case the joint probability distribution
of the positive eigenvalues $(E_1,\ldots ,E_n)\in\mathbb{R}^n_+$
can be written in the following form for all the seven novel
symmetry classes,
\begin{equation}\label{probability distribution Ej}
    \frac{1}{Z_n} \prod_{1\leq j< k \leq
    n}|E_j^2-E_k^2|^\beta \left(\prod_{\ell=1}^n E_\ell ^\gamma
    e^{-E_\ell^2}\right)dE_1\ldots dE_n,
\end{equation}
where $\beta\in\{1,2,4\}$, $\gamma\in\mathbb{N}=\{0,1,2,\ldots\}$
and $Z_n$ is a normalization constant. We note that for these ensembles
the eigenvalues come in pairs $\pm E_j$, i.e.\ for every positive eigenvalue
$E_j$ there corresponds a negavitive eigenvalue $-E_j$.
In the spirit of the universality conjecture in random matrix
theory we replace the term $e^{-E_\ell^2}$ in
(\ref{probability distribution Ej}) by $e^{-Q(E_\ell^2)}$ where
$Q$ denotes a polynomial of the form (\ref{definition: Q}).
A simple transformation $\lambda_j=E_j^2$ then leads to the
following probability density function,
\[
    F_n(\lambda_1,\ldots ,\lambda_n)=\frac{1}{\hat Z_n} \prod_{1\leq j<k\leq
    n}|\lambda_j-\lambda_k|^\beta \prod_{\ell=1}^n
    \lambda_\ell^{\frac{1}{2}(\gamma-1)}e^{-Q(\lambda_\ell)},
    \qquad \mbox{for $(\lambda_1,\ldots ,\lambda_n)\in \mathbb R_+^n$.}
\]
It is a beautiful and important observation (see for example
\cite{Deift,Mehta}) in random matrix theory that for $\beta=2$,
\[
    F_n(\lambda_1,\ldots ,\lambda_n)=\frac{1}{n!}\det(K_n(\lambda_j,\lambda_k))_{1\leq j,k\leq
    n},
\]
where $K_n$ is constructed out of orthogonal polynomials with respect to the
weight function $w(x)=x^{\frac{1}{2}(\gamma-1)}e^{-Q(x)}$. Moreover, all
related statistical quantities of interest can be expressed in terms of $K_n$.
Using our results on the asymptotics of the orthogonal polynomials we will
prove in this paper that for $\beta=2$ the local eigenvalue statistics have
universal behavior (as $n\to\infty$) in different regimes of the spectrum. By
{\em universal} we mean that the behavior is independent of $Q$ as long as $Q$
is chosen to be a polynomial. In work in progress, together with Deift, Gioev
and Kriecherbauer \cite{DGKV}, we will use the results of this paper, together
with the approach developed in \cite{DeiftGioev,DeiftGioev2,Widom}, to prove
universality for $\beta=1,4$.

\medskip

In order to get our results, we are inspired by the papers
\cite{DKMVZ2,DKMVZ1} of Deift et al. They considered orthogonal
polynomials with respect to varying weights $e^{-nV(x)}$ on
$\mathbb{R}$ where $V$ is real analytic and has enough increase at
infinity \cite{DKMVZ2}, as well as with respect to non-varying
weights $e^{-Q(x)}$ where $Q$ is an even polynomial with positive
leading coefficient \cite{DKMVZ1}. See also \cite{Deift} for an
excellent exposition. As in these papers, we will characterize the
orthogonal polynomials via the well-known $2\times 2$ matrix
valued Fokas-Its-Kitaev Riemann-Hilbert (RH) problem
\cite{FokasItsKitaev}, and apply the Deift-Zhou steepest descent
method, introduced in \cite{DeiftZhou} and further developed in
\cite{DVZ1,DVZ2,DZ2}, to analyze this RH problem (and thus also
the orthogonal polynomials) asymptotically.

The main difference between the orthogonal polynomials with
respect to the Laguerre-type weights (\ref{weight}) and the
weights $e^{-nV(x)}$ and $e^{-Q(x)}$ on $\mathbb{R}$, lies in the
behavior near the fixed endpoint 0. Near this endpoint we have to
do a local analysis using Bessel functions. This is analogous to
the local analysis near the fixed endpoints $\pm 1$ of the
modified Jacobi weight considered in \cite{KMVV}, see also
\cite{Kuijlaars}.

\medskip

The present paper is organized as follows. In the next section we
state our results. In Section 3 we state the RH problem for
orthogonal polynomials and apply the Deift-Zhou steepest descent
method to analyze the RH problem asymptotically. Afterwards, we
use the asymptotics of the solution of this RH problem to
determine the asymptotics of $a_n$, $b_{n-1}$ and $\gamma_n$ (in
Section 4), to determine the Plancherel-Rotach type asymptotics
for $p_n$ (in Section 5), and to prove universality results in
random matrix theory (in Section 6).

\section{Statement of results}

It is well-known in the theory of orthogonal polynomials that there
are two ingredients needed for the asymptotic (as $n\to\infty$) description.
First, we need the sequence of Mashkar-Rahmanov-Saff (MRS) numbers, which we will
denote by $\beta_n$. For the weigths we consider in this paper, $\beta_n$ is
uniquely determined for $n$ sufficiently large by the equation
\begin{equation}\label{introduction: MRS number}
    \frac{1}{2\pi} \int_0^{\beta_n} Q'(x)
        \sqrt{\frac{x}{\beta_n-x}}dx=n,
\end{equation}
see Proposition \ref{proposition: MRS}. From now on, we assume $n$
to be sufficiently large such that the MRS numbers exist. The
other ingredient we need, is the equilibrium measure
\cite{SaffTotik} of $[0,\infty)$ in the presence of the rescaled
field $V_n(x)=\frac{1}{n}Q(\beta_n x)$. This measure is defined as
the unique minimizer of the following minimization problem,
\begin{equation}\label{minimization problem}
    \inf
    \left(
        \iint\log\frac{1}{|s-t|}d\mu(s)d\mu(t)+\int V_n(t)d\mu(t)
    \right),
\end{equation}
where the infimum is taken over all probability measures on
$[0,\infty)$. We will show in Section \ref{section: equilibrium
measure} below, that the minimizer is of the form,
\begin{equation}
    d\mu_n(x) = \hat\psi_n(x)dx
        = \frac{1}{2\pi}\sqrt{\frac{1-x}{x}}h_n(x)\chi_{(0,1]}dx,
\end{equation}
where $\chi_{(0,1]}$ denotes the indicator function of the set
$(0,1]$, and where $h_n$ is a polynomial of degree $m-1$,
\begin{equation}
    h_n(x)=\sum_{k=0}^{m-1}h_{n,k}x^k,
\end{equation}
with real coefficients $h_{n,k}$, given by (\ref{definition: hn}),
(\ref{expression: Vn}) and (\ref{definition: Q}), and which have
an explicitly computable power series in $n^{-1/m}$, see
(\ref{definition: hk0}).

\subsection{Asymptotics of the recurrence coefficients
    $a_n,b_{n-1}$ and the leading coefficient $\gamma_n$}

To state our first result we need the quantity
$\ell_n=2\int_0^\infty \hat\psi_n(y)\log ydy - V_n(0)$.

\begin{theorem}\label{theorem: asymptotics coefficients}
    The recurrence coefficients $a_n$ and $b_{n-1}$ in the
    three-term recurrence relation {\rm (\ref{three-term recurrence relation})}
    of orthogonal polynomials with respect to the Laguerre-type weight
    {\rm (\ref{weight})}, have the following asymptotic behavior,
    \begin{align}
        \label{theorem: asymptotics coefficients: bn-1}
        & \frac{b_{n-1}}{\beta_n} = \frac{1}{4}+\frac{\alpha}{2h_n(1)n}+\bigO(1/n^2),
            \qquad\mbox{as $n\to\infty$},
        \\[1ex]
        \label{theorem: asymptotics coefficients: an}
        & \frac{a_n}{\beta_n} = \frac{1}{2}+\frac{\alpha+1}{h_n(1)n}+\bigO(1/n^2),
            \qquad\mbox{as $n\to\infty$}.
    \end{align}
    The leading coefficients $\gamma_n$ of the orthonormal
    polynomials with respect to the Laguerre-type weight {\rm (\ref{weight})}
    have the following asymptotic behavior,
    \begin{multline}\label{theorem: asymptotics coefficients: gamman}
        \gamma_n = \beta_n^{-(n+\frac{\alpha}{2}+\frac{1}{2})}
            e^{-\frac{1}{2}n\ell_n} \sqrt{\frac{2}{\pi}} \, 2^\alpha
        \\[1ex]
            \times\,
            \left[
                1-\left(
                    \frac{4\alpha^2-1}{8h_n(0)}
                    + \frac{12\alpha^2+24\alpha+11}{24h_n(1)} - \frac{h_n'(1)}{8h_n(1)^2}
                \right) \frac{1}{n} + \bigO(1/n^2)
            \right],
            \quad\mbox{as $n\to\infty$.}
    \end{multline}
    Each of the number sequences $\beta_n$, $\ell_n$, $h_n(0)$,
    $h_n(1)$ and $h_n'(1)$ as well as each of the error terms
    have an asymptotic expansion in powers of
    $n^{-1/m}$ which can be calculated explicitly.
\end{theorem}

\begin{remark}
    For the convenience of the reader we summarize here the
    leading order behavior of the quantities appearing in
    (\ref{theorem: asymptotics coefficients: bn-1})--(\ref{theorem: asymptotics coefficients:
    gamman}):
    \begin{align*}
        & \beta_n=n^{1/m}\sum_{j=0}^\infty \beta^{(j)}n^{-j/m},
            && \beta^{(0)}=\bigl(\frac{1}{2} m q_m A_m\bigr)^{-1/m},
            \qquad A_m=\prod_{j=1}^m \frac{2j-1}{2j}, \\
        & \ell_n=\sum_{j=0}^\infty
        \ell^{(j)}n^{-j/m},&& \ell^{(0)}=-\frac{2}{m}-4\log 2,
        \\
        & h_n(0)=\sum_{j=0}^\infty h^{(j)}n^{-j/m},&& h^{(0)}=
        \frac{4m}{2m-1},
        \\
        & h_n(1)=\sum_{j=0}^\infty \tilde h^{(j)}n^{-j/m},&& \tilde
        h^{(0)}=4m,
        \\
        & h_n'(1)=\sum_{j=0}^\infty \hat h^{(j)}n^{-j/m},&& \hat
        h^{(0)}=\frac{8}{3}m(m-1),
    \end{align*}
    see Proposition \ref{proposition: MRS}, Remark \ref{remark: xin and
    ln}, and equations (\ref{proof hn polynomial: eq1})--(\ref{h(0)h(1)h'(1)}).
\end{remark}

\begin{remark}\label{remark: special case an bn-1 gamman: introduction}
    The special case $Q(x)=q_m x^m$: In this case, it will be clear from
    Remarks \ref{remark: special case betan}, \ref{remark: special case hn}
    and \ref{remark: special case Hn and ln} that the values of
    $n^{-1/m}\beta_n$, $\ell_n$, $h_n(0)$, $h_n(1)$ and $h_n'(1)$ do not
    depend on $n$ and we have $\beta_n=n^{1/m}\beta^{(0)}$, $\ell_n=\ell^{(0)}$,
    $h_n(0)=h^{(0)}$, $h_n(1)=\tilde h^{(0)}$ and $h_n'(1)=\hat h^{(0)}$.
    Furthermore, the asymptotic expansions of the error terms in the
    theorem are given in powers of $1/n$ rather than in powers of $n^{-1/m}$.
\end{remark}

\subsection{Plancherel-Rotach type asymptotics for the orthonormal
    polynomials $p_n$}

In order to state the asymptotic behavior (as $n\to\infty$) of the
rescaled orthonormal polynomials $p_n(\beta_n z)$, with $\beta_n$
the MRS number, we need to introduce some more notation. Let
\begin{equation}\label{definition: psin: introduction}
    \psi_n(z) = \frac{1}{2\pi i} \frac{(z-1)^{1/2}}{z^{1/2}} h_n(z),
        \qquad\mbox{for $z\in\mathbb{C}\setminus[0,1]$.}
\end{equation}
Throughout this paper we always take the principal branch of the
power. It is clear that
$\psi_{n,+}(x)=-\psi_{n,-}(x)=\hat\psi_n(x)$ for $x\in (0,1)$,
where $\hat\psi_n$ is the density of $\mu_n$.

Further, we will show in Sections \ref{section: parametrix near 1}
and \ref{section: parametrix near 0} below, that there exist
biholomorphic maps $f_n$ and $\tilde f_n$ in a neighborhood of 1
and 0, respectively, satisfying
\begin{align}
    \label{statement of results: fn}
    & \frac{2}{3}f_n(z)^{3/2} = - \pi in\int_z^1 \psi_n(s)ds,
        \qquad \mbox{for $|z-1|$ small and $z\notin(-\infty,1]$,}
    \\[1ex]
    \label{statement of results: fntilde}
    & 2\tilde f_n(z)^{1/2} = -\pi in\int_0^z\psi_n(s)ds,
        \qquad \mbox{for $|z|$ small and $z\notin[0,\infty)$.}
\end{align}
See (\ref{definition: fn}), (\ref{definition: phin}) and
(\ref{definition: phinhat}) for an explicit expression of $f_n$,
and see (\ref{definition: fntilde}), (\ref{definition: phintilde})
and (\ref{definition: phintildehat}) for the explicit expression
of $\tilde f_n$.

A last function which we will need is the conformal map $\varphi$
from $\mathbb{C}\setminus[0,1]$ onto the exterior of the unit
circle,
\begin{equation}\label{definition: varphi: eq1}
    \varphi(z) = 2(z-1/2)+2z^{1/2}(z-1)^{1/2},
        \qquad\mbox{for $z\in\mathbb{C}\setminus[0,1]$.}
\end{equation}

\medskip

\begin{figure}[t]
\begin{center}
    \setlength{\unitlength}{1mm}
    \begin{picture}(100,30)(0,13)
        \put(25,20){\thicklines\circle*{.8}}
        \put(75,20){\thicklines\circle*{.8}}
        \put(5,20){\line(1,0){90}}

        \cCircle(25,20){12}[t]
        \cCircle(75,20){12}[t]
        \put(33.06,29){\line(1,0){34}}

        \put(47.5,35){\small $A_\delta$}
        \put(47.5,23){\small $B_\delta$}
        \put(73,23){\small $C_\delta$}
        \put(23,23){\small $D_\delta$}

        \put(10.5,15){\small $-\delta$}
        \put(24,15){\small 0}
        \put(36,15){\small $\delta$}
        \put(58.7,15){\small $1-\delta$}
        \put(74.5,15){\small 1}
        \put(84.3,15){\small $1+\delta$}
    \end{picture}
    \caption{The different asymptotic regions for
        $p_n(\beta_n z)$ in the upper half-plane.}
    \label{figure: asymptotic regions}
\end{center}
\end{figure}
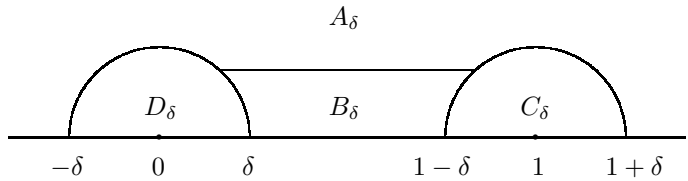

Because of the symmetry $p_n(z)=\overline{p_n(\overline z)}$ we
only need to present the asymptotics of $p_n(\beta_n z)$ in the
closed upper half-plane $\overline{\mathbb{C}_+}$. We state the
asymptotic formulae in the four closed regions $A_\delta$,
$B_\delta$, $C_\delta$ and $D_\delta$, depicted in Figure
\ref{figure: asymptotic regions}.

\begin{theorem}\label{theorem: asymptotics pn}
    Let $p_n$ be the $n$-th degree orthonormal polynomial with
    respect to the Laguerre-type weight {\rm (\ref{weight})}.
    There exists $\delta_0>0$ such that for all $\delta\in(0,\delta_0)$
    the $p_n$ have the following asymptotic behavior (as $n\to\infty$).
    \begin{itemize}
    \item[{\rm (a)}] For $z$ in the outside region $A_\delta$,
        \begin{multline} \label{theorem: asymptotics pn: eq1}
            p_n(\beta_n z) = (\beta_n z)^{-\frac{\alpha}{2}}
                e^{\frac{1}{2}Q(\beta_n z)}
            \\[1ex] \times\,
                \sqrt{\frac{2}{\pi\beta_n}}
                \frac{\varphi(z)^{\frac{1}{2}(\alpha+1)}}{2z^{1/4}(z-1)^{1/4}}
                \exp\left(-\pi in\int_1^z\psi_n(s)ds\right) (1+\bigO(1/n)).
        \end{multline}
    \item[{\rm (b)}] For $z$ in the bulk region $B_\delta$,
        \begin{multline} \label{theorem: asymptotics pn: eq2}
            p_n(\beta_n z) = (\beta_n z)^{-\frac{\alpha}{2}}
                e^{\frac{1}{2}Q(\beta_n z)} \sqrt{\frac{2}{\pi\beta_n}}
                \frac{1}{z^{1/4}(1-z)^{1/4}}
            \\[1ex] \times\,
                \left[
                \cos\left( \frac{1}{2}(\alpha+1)\arccos(2z-1)
                    - \pi n \int_1^z\psi_n(s)ds
                    - \frac{\pi}{4} \right) (1+\bigO(1/n))
                \right. \\[1ex]
                \left. +\,
                \cos\left( \frac{1}{2}(\alpha-1)\arccos(2z-1)
                    - \pi n \int_1^z\psi_n(s)ds
                    - \frac{\pi}{4} \right) \bigO(1/n)
                \right],
        \end{multline}
    \item[{\rm (c)}] For $z$ in the Airy region $C_\delta$,
        \begin{multline}\label{theorem: asymptotics pn: eq3}
            p_n(\beta_n z) = (\beta_n z)^{-\frac{\alpha}{2}}
                e^{\frac{1}{2}Q(\beta_n z)}\sqrt{\frac{2}{\beta_n}}
                \frac{1}{z^{1/4}(z-1)^{1/4}}
            \\[1ex] \times\,
                \left[
                \cos\left(\frac{1}{2}(\alpha+1)\arccos(2z-1)\right)
                f_n(z)^{1/4} \Ai(f_n(z)) (1+\bigO(1/n))
                \right. \\[1ex]
                \left. -\,
                i \sin\left(\frac{1}{2}(\alpha+1)\arccos(2z-1)\right)
                f_n(z)^{-1/4} \Ai'(f_n(z)) (1+\bigO(1/n))
                \right],
        \end{multline}
        with $\Ai$ the Airy function.
    \item[{\rm (d)}] For $z$ in the Bessel region $D_\delta$,
        \begin{multline}
            p_n(\beta_n z) = (\beta_n z)^{-\frac{\alpha}{2}}
                e^{\frac{1}{2}Q(\beta_n z)}
                (-1)^n \sqrt{\frac{1}{\beta_n}}
                \frac{\sqrt 2(-\tilde f_n(z))^{1/4}}{z^{1/4}(1-z)^{1/4}}
            \\[1ex] \times\,
                \left[
                \sin\zeta_1(z) J_\alpha\left(2(-\tilde f_n(z))^{1/2}\right) (1+\bigO(1/n))
                \right. \\
                \left. +\,
                \cos\zeta_1(z) J_\alpha'\left(2(-\tilde f_n(z))^{1/2}\right) (1+\bigO(1/n))
                \right],
        \end{multline}
        with $J_\alpha$ the $J$-Bessel function of order $\alpha$,
        and with
        \begin{equation} \label{theorem: asymptotics pn: eq5}
            \zeta_1(z) =
                \frac{1}{2}(\alpha+1)\arccos(2z-1) - \frac{\pi\alpha}{2}.
        \end{equation}
    \end{itemize}
    All the error terms are uniform for $\delta$ in compact
    subsets of $(0,\delta_0)$ and for $z\in X_\delta$ with
    $X\in\{A,B,C,D\}$. The error terms have an explicit asymptotic
    expansion in $n^{-1/m}$.
\end{theorem}

\begin{remark}\label{remark: arccos}\
    \begin{enumerate}
    \item The function $\arccos z$ which appears in (\ref{theorem: asymptotics pn: eq2}),
        (\ref{theorem: asymptotics pn: eq3}) and (\ref{theorem: asymptotics pn: eq5}) is
        defined as the inverse function of
        $\cos z: \{0<\Re z<\pi\}\to \mathbb{C}\setminus\left((-\infty,-1]\cup[1,\infty)\right)$,
        and is given by,
        \[
            \arccos z=\int_z^1\frac{ds}{(1-s)^{1/2}(1+s)^{1/2}},\qquad
            \mbox{for $z\in\mathbb{C}\setminus\left((-\infty,-1]\cup[1,\infty)\right)$.}
        \]
    \item
        Some of the expressions in Theorem \ref{theorem: asymptotics pn} are not well defined
        for all $z\in\mathbb{R}$. For example:
        $z^{-\frac{\alpha}{2}}$, $\arccos(2z-1)$,
        $\varphi(z)^{\frac{1}{2}(\alpha+1)}$, $(z-1)^{1/4}$,
        $\int_1^z\psi_n(s)ds$, etc. For these expressions we then take
        the limiting values as $z$ is approached from the upper
        half-plane.
    \item
        The function $\int_1^z\psi_n(s)ds$ which appears in
        (\ref{theorem: asymptotics pn: eq1}) and (\ref{theorem: asymptotics pn:
        eq2}) is explicitly computable. The result is, see Remark
        \ref{remark: xin and ln} below,
        \begin{equation}\label{introduction: intpsin}
            \int_1^z \psi_n(s)ds=\frac{1}{2\pi} H_n(z) z^{1/2}
            (1-z)^{1/2} -\frac{2}{\pi}\arccos z^{1/2},\qquad\mbox{for
            $z\in\mathbb{C}_+$,}
        \end{equation}
        where $H_n$ is a polynomial with real coefficients of degree $m-1$, given by
        (\ref{definition: Hn}).
    \end{enumerate}
\end{remark}

\begin{remark}
    The above theorem is in agreement with \cite[Theorem 8.22.8]{Szego} where the
    Plancherel-Rotach type asymptotics for the classical Laguerre polynomials
    were stated.
\end{remark}

\begin{remark}
    The special case $Q(x)=q_m x^m$: In this case $\beta_n=n^{1/m}(\frac{1}{2}m q_m A_m)^{-1/m}$,
    cf.\ Remark \ref{remark: special case betan},
    and since $h_n$ is $n$-independent it follows from (\ref{definition: psin: introduction})
    that $\psi_n$, and thus also $\int_1^z \psi_n(s)ds$, is $n$-independent.
    The latter integral is explicitly given by (\ref{introduction:
    intpsin}) where $H_n(z)=\frac{1}{m}h(z)$ with $h$ defined
    in Remark \ref{remark: special case Hn and ln}.
    Furthermore, the asymptotic expansions of the error terms in the theorem are given in
    powers of $1/n$ rather than in powers of $n^{-1/m}$.
\end{remark}

\subsection{Applications in random matrix theory}

Here, we consider matrices taken from the random matrix ensembles which induce
the following probability distribution on the positive eigenvalues $(E_1,\ldots E_n)\in\mathbb{R}^n_+$,
\begin{equation}\label{statement results: probability
distribution}
    P^{(n)}(E_1,\ldots ,E_n)d^n E=\frac{1}{Z_n}\prod_{1\leq j<k\leq n}|E_j^2-E_k^2|^2
    \left(\prod_{\ell=1}^n E_\ell^\gamma e^{-Q(E_\ell^2)})\right)d^n E,
\end{equation}
where we use $d^n E$ to denote $dE_1\ldots dE_n$. As noted in the introduction, after a
transformation $\lambda_j=E_j^2$ the eigenvalue statistics
can be expressed in terms of a scalar 2-point kernel $K_n$ constructed out of orthogonal polynomials $p_k$
with respect to the
weight function $w(x)=x^{\frac{1}{2}(\gamma-1)e^{-Q(x)}}$,
\begin{align}\label{definition: Kn}
    \nonumber K_n(x,y)
    & = \sqrt{w(x)}\sqrt{w(y)}
        \sum_{k=0}^{n-1} p_k(x) p_k(y) \\[1ex]
    & = \sqrt{w(x)}\sqrt{w(y)}\frac{\gamma_{n-1}}{\gamma_n}\frac{p_n(x)p_{n-1}(y)-p_{n-1}(x)p_n(y)}{x-y}.
\end{align}
For example, consider the probability $P_n(L_E\leq x)$ that the
smallest eigenvalue $L_E=\min_j(E_j)$ is less than or equal to
$x\in\mathbb{R}_+$. Obviously we have $P_n(L_E\leq
x)=P_n(L_\lambda\leq x^2)$, with $L_\lambda=\min_j(\lambda_j)$ and
$\lambda_j=E_j^2$. Since it is well known in random matrix theory, see e.g.\ \cite{Deift},
that $P_n(L_\lambda\leq x^2)=1-\det(I-\mathcal K_n)$, where $\mathcal K_n$ is the trace class operator
with integral kernel $K_n(x,y)$ acting on $L^2(0,x^2)$
and where $\det(I-\mathcal K_n)$ is a Fredholm determinant, we then obtain
\begin{equation}\label{smallest eigenvalue}
    P_n(L_E\leq x)=1-\det(I-\mathcal K_n).
\end{equation}


For the case $Q(x)=x$ it is known that the behavior of $K_n$ is
described by the sine kernel in the bulk of the spectrum
\cite{FoxKahn, NagaoWadati}, by the Airy kernel at the soft edge
of the spectrum \cite{Forrester,TW1}, and by the Bessel kernel at
the hard edge of the spectrum \cite{Forrester,TracyWidom}. In the
next theorem, we state that this behavior persists for any
polynomial $Q$ of the form (\ref{definition: Q}).

\begin{theorem}\label{theorem: universality}
    Let $w$ be the Laguerre-type weight {\rm (\ref{weight})}
    and let $K_n$ be the kernel {\rm (\ref{definition: Kn})}
    associated to $w$. Then the following holds.
    \begin{itemize}
        \item[{\rm (a)}] The bulk of the spectrum: With
        $\hat\psi_n(x)=\frac{1}{2\pi}\sqrt{\frac{1-x}{x}}h_n(x)$ the density of the
        equilibrium measure $\mu_n$,
        \begin{multline}
            \frac{\beta_n}{n\hat\psi_n(x)}K_n\left[\beta_n\left(x+\frac{u}{n\hat\psi_n(x)}\right),
                \beta_n\left(x+\frac{v}{n\hat\psi_n(x)}\right)\right]\\[1ex]
                =\frac{\sin(\pi(u-v))}{\pi(u-v)}+\bigO\left(\frac{1}{n}\right),
        \end{multline}
        as $n\to\infty$, uniformly for $u,v$ in compact subsets of
        $\mathbb{R}$ and $x$ in compact subsets of $(0,1)$.
        \item[{\rm (b)}] The soft edge of the spectrum: With $c_n=\left(\frac{1}{2}h_n(1)\right)^{2/3}$,
        \begin{multline}
            \frac{\beta_n}{c_n n^{2/3}}
            K_n\left[\beta_n\left(1+\frac{u}{c_n n^{2/3}}\right),
            \beta_n\left(1+\frac{v}{c_n n^{2/3}}\right)\right]\\[1ex]
            =\frac{\Ai(u)\Ai'(v)-\Ai(v)\Ai'(u)}{u-v}+\bigO\left(\frac{1}{n^{1/3}}\right),
        \end{multline}
        as $n\to\infty$, uniformly for $u,v$ in compact subsets of
        $\mathbb{R}$.
        \item[{\rm (c)}] The hard edge of the spectrum: With $\tilde
        c_n=\left(\frac{1}{2}h_n(0)\right)^2$,
        \begin{multline}\label{theorem: universality: eq 3}
            \frac{\beta_n}{4\tilde c_n n^2}
            K_n\left[\beta_n \frac{u}{4\tilde c_n n^2},
            \beta_n\frac{v}{4\tilde c_n n^2}\right] \\[1ex]
            =\frac{J_\alpha(\sqrt u)\sqrt v J_\alpha'(\sqrt v)-J_\alpha(\sqrt v)\sqrt u J_\alpha'(\sqrt
            u)}{2(u-v)}+\bigO\left(\frac{u^{\alpha/2}v^{\alpha/2}}{n}\right),
        \end{multline}
        as $n\to\infty$, uniformly for $u,v$ in bounded subsets of
        $(0,\infty)$.
    \end{itemize}
\end{theorem}

\begin{remark}
    As in \cite[Section 3.4]{KV1} one can show, using (\ref{smallest eigenvalue}) and
    (\ref{theorem: universality: eq 3}), that
    the smallest eigenvalue $L_E=\min(E_1,\ldots ,E_n)$, where the
    $E_1,\ldots ,E_n$ have probability distribution (\ref{statement results: probability distribution}),
    satisfies
    \[
        \lim_{n\to\infty}P_n\left(\frac{2\tilde c_n^{1/2} n}{\beta_n^{1/2}}L_E \leq s\right)
            =1-\det\left(I-\mathbb{J}_{\frac{1}{2}(\gamma-1),s^2}\right).
    \]
    Here, we use $\mathbb{J}_{\alpha,s}$ to denote the integral operator with kernel
    $\frac{J_\alpha(\sqrt u)\sqrt v J_\alpha'(\sqrt v)-J_\alpha(\sqrt
    v)\sqrt u J_\alpha'(\sqrt u)}{2(u-v)}$
    acting on $L^2(0,s)$. Tracy and Widom \cite{TracyWidom} have determined an explicit expression
    for the Fredholm determinant in the
    latter equation in terms of a Painlev\'e V transcendent. Let
    $q_\alpha(s)$ be the unique solution of the differential equation,
    \[
        s(q_\alpha^2-1)(s q_\alpha')'=q_\alpha(s
        q_\alpha')^2+\frac{1}{4}(s-\alpha^2)q_\alpha+\frac{1}{4}sq_\alpha^3(q_\alpha^2-2),\qquad
        ('=\frac{d}{ds}),
    \]
    with the boundary condition,
    \[
        q_\alpha(s)\sim\frac{1}{2^\alpha
        \Gamma(1+\alpha)}s^{\alpha/2},\qquad \mbox{as $s\to 0$.}
    \]
    Note that this equation is reducible to a special case of the Painlev\'e
    V differential equation, see for example \cite{TracyWidom}.
    Further, let
    \[
        F_\alpha(s)=\exp\left(-\frac{1}{4}\int_0^s
        \log\frac{s}{x}q_\alpha(x)dx\right).
    \]
    Then, it has been shown in \cite{TracyWidom} that
    $\det(I-\mathbb{J}_{\alpha,s})=F_\alpha(s)$, so that
    \begin{equation}
        \lim_{n\to\infty}P_n\left(\frac{2\tilde c_n^{1/2}n}{\beta_n^{1/2}}L_E \leq s\right)
            =1-F_{\frac{1}{2}(\gamma-1)}(s^2).
    \end{equation}
\end{remark}

\medskip

Another application which we will present here concerns averages of products
and ratios of characteristic polynomials $\det(xI-M)$ of random matrices $M$
taken from the random matrix ensemble which induces the probability
distribution (\ref{statement results: probability distribution}) on the
positive eigenvalues $E_1,\ldots E_n$. Recall that the eigenvalues of $M$ come
in pairs $\pm E_j$, so that $\det(xI-M)=\prod_{j=1}^n(x^2-E_j^2)$. The averages
are then intimately related with kernels $W_{\I,n}$, $W_{\II,n}$ and
$W_{\III,n}$ constructed out of orthogonal polynomials (with respect to the
weight function $w(x)=x^{\frac{1}{2}(\gamma-1)}e^{-Q(x)}$) and their Cauchy
transforms, see \cite{BaikDeiftStrahov,StrahovFyodorov}. Which kernel we have
to use depends on whether the characteristic polynomials are only in the
numerator ($W_{\I,n}$), only in the denominator ($W_{\III,n}$) or in both the
numerator and denominator ($W_{\II,n}$). We will not state for brevity reasons
the explicit formulae of this connection here, but refer the reader to
\cite{BaikDeiftStrahov,StrahovFyodorov} (see also \cite{v5}). The three kernels
are given by,
\begin{align}
    \label{definition: WI}
        & W_{\I,n}(u,v) =
            \frac{\pi_n(u)\pi_{n-1}(v)-\pi_{n-1}(u)\pi_n(v)}{u-v},
        \\[1ex]
    \label{definition: WII}
        & W_{\II,n}(u,v) =
            \frac{C(\pi_n w)(u)\pi_{n-1}(v)-C(\pi_{n-1}w)(u)\pi_n(v)}{u-v},
\end{align}
and
\begin{equation}\label{definition: WIII}
    W_{\III,n}(u,v) =
        \frac{C(\pi_n w)(u)C(\pi_{n-1}w)(v)-C(\pi_{n-1}w)(u)C(\pi_n w)(v)}{u-v},
\end{equation}
with $\pi_n=\gamma_n^{-1}p_n=z^n+\cdots$ the $n$-th degree monic
orthogonal polynomial with respect to $w$, and $Cf$ the Cauchy
transform of $f$ on the positive real line, i.e.
\[
    Cf(z)=\frac{1}{2\pi i}\int_0^\infty \frac{f(x)}{x-z}dx,
        \qquad\mbox{for $z\in\mathbb{C}\setminus[0,\infty)$.}
\]

In the next theorem, we state that these kernels (and thus also
the averages of characteristic polynomials) have universal
behavior, as $n\to\infty$, at the hard edge of the spectrum in
terms of Bessel functions.

\begin{table}
\begin{center}
    \begin{tabular}{|l|c|}
        \hline
        \multicolumn{2}{|c|}{} \\[-0,5ex]
        \multicolumn{2}{|c|}{Limiting Bessel kernels}\\
        \multicolumn{2}{|c|}{} \\[-0,5ex]
        \hline
        & \\[-0,5ex]
            $\mathbb{J}_{\alpha,\I}(u,v)$ &
            $u^{-\frac{\alpha}{2}}v^{-\frac{\alpha}{2}}
                \frac{J_\alpha(u^{1/2}) v^{1/2} J_\alpha'(v^{1/2})-
                J_\alpha(v^{1/2}) u^{1/2}J_\alpha'(u^{1/2})}{2(u-v)}$
                \\
        & \\[-0,5ex]
        \hline
        & \\[-0,5ex]
            $\mathbb{J}_{\alpha,\II}^{+}(u,v)$ &
            $u^{\frac{\alpha}{2}} v^{-\frac{\alpha}{2}}
                \frac{H_\alpha^{(1)}(u^{1/2}) v^{1/2} J_\alpha'(v^{1/2})-
                J_\alpha(v^{1/2}) u^{1/2}(H_\alpha^{(1)})'(u^{1/2})}{4(u-v)}$ \\[4ex]
            $\mathbb{J}_{\alpha,\II}^{-}(u,v)$ &
            $-u^{\frac{\alpha}{2}} v^{-\frac{\alpha}{2}}
                \frac{H_\alpha^{(2)}(u^{1/2}) v^{1/2} J_\alpha'(v^{1/2})-
                J_\alpha(v^{1/2}) u^{1/2}(H_\alpha^{(2)})'(u^{1/2})}{4(u-v)}$ \\
        & \\[-0,5ex]
        \hline
        & \\[-0,5ex]
            $\mathbb{J}_{\alpha,\III}^{+}(u,v)$ &
            $u^{\frac{\alpha}{2}} v^{\frac{\alpha}{2}}
                \frac{H_\alpha^{(1)}(u^{1/2}) v^{1/2} (H_\alpha^{(1)})'(v^{1/2})-
                H_\alpha^{(1)}(v^{1/2}) u^{1/2}(H_\alpha^{(1)})'(u^{1/2})}{8(u-v)}$
            \\[4ex]
            $\mathbb{J}_{\alpha,\III}^{\pm}(u,v)$ &
            $-u^{\frac{\alpha}{2}} v^{\frac{\alpha}{2}}
                \frac{H_\alpha^{(1)}(u^{1/2}) v^{1/2} (H_\alpha^{(2)})'(v^{1/2})-
                H_\alpha^{(2)}(v^{1/2}) u^{1/2}(H_\alpha^{(1)})'(u^{1/2})}{8(u-v)}$
            \\[4ex]
            $\mathbb{J}_{\alpha,\III}^{-}(u,v)$ &
            $u^{\frac{\alpha}{2}} v^{\frac{\alpha}{2}}
                \frac{H_\alpha^{(2)}(u^{1/2}) v^{1/2} (H_\alpha^{(2)})'(v^{1/2})-
                H_\alpha^{(2)}(v^{1/2})
                u^{1/2}(H_\alpha^{(2)})'(u^{1/2})}{8(u-v)}$
            \\[-0,5ex]
        & \\
        \hline
    \end{tabular}
    \caption{Expressions for the limiting Bessel kernels. Here, $J_\alpha$
        is the usual $J$-Bessel function of order
        $\alpha$, and
        $H_{\alpha}^{(1)}$ and $H_{\alpha}^{(2)}$ are
        the Hankel functions of order $\alpha$ of the first and second kind, respectively.}
    \label{table: limiting-Bessel-kernels}
\end{center}
\end{table}

\begin{theorem}\label{theorem: averages}
    Let $w$ be the Laguerre-type weight {\rm (\ref{weight})} and
    let $W_{\I,n}, W_{\II,n}$ and $W_{\III,n}$ be the kernels
    {\rm (\ref{definition: WI})}--{\rm (\ref{definition: WIII})}
    associated to $w$. Further, let $\mathbb{J}_{\alpha,\I}$,
    $\mathbb{J}_{\alpha,\II}^+, \ldots$ be
    the Bessel kernels given by Table
    {\rm \ref{table: limiting-Bessel-kernels}}, and let
    $\tilde c_n=\bigl(\frac{1}{2}h_n(0)\bigr)^2$. Then the following holds.
    \begin{itemize}
    \item[{\rm (a)}]
        The kernel $W_{\I,n}$ satisfies,
        \begin{multline}
            \gamma_{n-1}^2\frac{\beta_n}{4\tilde c_n n^2}
            W_{\I,n}\left(\beta_n \frac{u}{4\tilde c_n n^2}
                ,\beta_n \frac{v}{4\tilde c_n n^2}\right) \\[1ex]
            = \left(\frac{\beta_n}{4\tilde c_n n^2}\right)^{-\alpha}e^{Q(0)}
                \left(\mathbb{J}_{\alpha,\I}(u,v) + \bigO\left(\frac{1}{n}\right)\right),
        \end{multline}
        as $n\to\infty$, uniformly for $u,v$ in compact
        subsets of $\mathbb{C}$.
    \item[{\rm (b)}]
        The kernel $W_{\II,n}$ satisfies,
        \begin{equation}
            \gamma_{n-1}^2\frac{\beta_n}{4\tilde c_n n^2}
            W_{\II,n}\left(\beta_n \frac{u}{4\tilde c_n n^2}
                ,\beta_n \frac{v}{4\tilde c_n n^2}\right)
            = \mathbb{J}_{\alpha,\II}^+(u,v) + \bigO\left(\frac{1}{(u-v)n}\right),
        \end{equation}
        as $n\to\infty$, uniformly for $u$ and $v$ in compact
        subsets of $\mathbb{C}_+$ and $\mathbb{C}$, respectively.
        The corresponding behavior of $W_{\II,n}$ that holds uniformly for $u$ and
        $v$ in compact subsets of $\mathbb{C}_-$ and $\mathbb{C}$,
        respectively, is stated by replacing $\mathbb{J}_{\II,\alpha}^+$ with
        $\mathbb{J}_{\II,\alpha}^-$.
    \item[{\rm (c)}]
        The kernel $W_{\III,n}$ satisfies,
        \begin{multline}
            \gamma_{n-1}^2\frac{\beta_n}{4\tilde c_n n^2}
                W_{\III,n}\left(\beta_n \frac{u}{4\tilde c_n n^2}
                    ,\beta_n \frac{v}{4\tilde c_n n^2}\right)\\[1ex]
            =\left(\frac{\beta_n}{4\tilde c_n n^{2}}\right)^\alpha e^{-Q(0)}
                \left(\mathbb{J}_{\alpha,\III}^\pm(u,v)+\bigO\left(\frac{1}{n}\right)\right),
        \end{multline}
        as $n\to\infty$, uniformly for $u$ and $v$ in compact
        subsets of $\mathbb{C}_+$ and $\mathbb{C}_-$, respectively.
        The corresponding behavior of $W_{\III,n}$ that holds uniformly for $u,v$
        in compact subsets of $\mathbb{C}_+$ is stated by replacing $\mathbb{J}_{\III,\alpha}^\pm$
        with $\mathbb{J}_{\III,\alpha}^+$. Further,
        the behavior of $W_{\III,n}$ that holds uniformly for $u,v$
        in compact subsets of $\mathbb{C}_-$ is stated by replacing $\mathbb{J}_{\III,\alpha}^\pm$
        with $\mathbb{J}_{\III,\alpha}^-$.
    \end{itemize}
\end{theorem}

\begin{remark}
    From our analysis we can also prove that the three kernels
    $W_{\I,n}$, $W_{\II,n}$ and $W_{\III,n}$ (and thus also the averages of characteristic
    polynomials) have universal behavior at the soft edge 1 of the spectrum, and in the bulk
    of the spectrum. The proof of these facts is similar to the proof
    of Theorem \ref{theorem: averages}. In the bulk of the
    spectrum the results will be analogous to the results of Strahov
    and Fyodorov in \cite{StrahovFyodorov}.
\end{remark}

\subsection{The varying weights $x^\alpha e^{-nV(x)}$}

Consider the varying weights $x^\alpha e^{-nV(x)}$ on
$[0,\infty)$, where $V$ is real analytic and has sufficient growth
at infinity. We now briefly explain how one can obtain the
analogues of the results of this paper for these weights. There are two ways to proceed.
\medskip

One approach is to link the orthogonal polynomials via a simple
transformation to orthogonal polynomials with respect to symmetric
weights of the form $|x|^\gamma e^{-n V(x^2)}$ on $\mathbb{R}$.
The latter weights fall into the class of weights considered in
\cite{KV2}. Recall that in \cite{KV2} the RH approach was used to
prove universality for the associated random matrix ensembles, and
not to determine asymptotics of the orthogonal polynomials.
However, one can use the techniques of \cite{KV2} to determine the
asymptotics of the orthogonal polynomials.

\medskip

If one is interested in applications in random matrix theory, a second approach
is more convenient since the kernel $K_n$, see (\ref{definition: Kn}),
associated to the weight $x^\alpha e^{-nV(x)}$ on $[0,\infty)$ doesn't relate
directly to the kernel associated to the weight $|x|^\gamma e^{-n V(x^2)}$ on
$\mathbb{R}$ (they only relate via the orthogonal polynomials). One can apply
the RH approach directly to the varying weights on the half real line. The
analysis is then very similar as in the present paper, although the situation
is more complicated. The complications come from the fact that the support
$S_V$ of the equilibrium measure $\mu_V$ of $[0,\infty)$ in external field $V$
can consist of more than one interval,
\[
    S_V=[0,a_1]\cup \bigcup_{j=1}^L [b_j,a_{j+1}]\qquad \mbox{or}
    \qquad S_V=\bigcup_{j=1}^L [b_j,a_{j+1}], \qquad
    a_i,b_j>0,
\]
and that $\mu_V$ can have degenerate behavior.

One distinguishes four types of degenerate behavior depending at which point in
$\mathbb R_+$ the degenerate behavior occurs, cf.\
\cite{DKMVZ2,KuijlaarsMcLaughlin}. Type I: points in
$[0,\infty)\setminus S_V$ where equality holds in the variational
condition (\ref{variational condition 2}) below. Type II: interior
points of $S_V$ where the density $\psi_V$ of $\mu_V$ vanishes.
Type III: edge points $a_i,b_j$ where $\psi_V$ vanishes to higher
order than a square root. And type IV: the edge point 0 if $0\in S_V$
and $\psi_V(0)=0$. Near each of these singular points we have
to construct special local parametrices. For type I-III
singularities, the existence of these parametrices is established
in \cite[Section 5]{DKMVZ2}.

In the nondegenerate case, which is generic, the construction of the parametrices at the
endpoints is as in the present paper.  In particular, for
regular soft edge points $a_i, b_j$ this will be done using Airy
functions, and for a regular edge point 0 using Bessel functions.

If $S_V$ consists of more than one interval, the construction of
the so-called parametrix for the outside region (see Section
\ref{section: parametrix outside region} below for the
one-interval case) uses $\Theta$-functions as in \cite[Section
4.2]{DKMVZ2}.

\section{Asymptotic analysis of the RH problem for orthogonal polynomials}
    \label{section: asymptotic analysis}

Here, we will perform the asymptotic analysis of the
Fokas-Its-Kitaev RH problem for orthogonal polynomials
\cite{FokasItsKitaev}. This RH problem is the following. Seek a
$2\times 2$ matrix valued function $Y(z)=Y(z;n,w)$ which satisfies
the following conditions.

\subsubsection*{RH problem for $Y$:}

\begin{itemize}
    \item[(a)] $Y:\mathbb{C}\setminus [0,\infty)\to\mathbb{C}^{2\times
    2}$ is analytic.
    \item[(b)] $Y$ possesses continuous boundary values for
    $x\in(0,\infty)$ denoted by $Y_+(x)$  and $Y_-(x)$, where
    $Y_+(x)$ and $Y_-(x)$ denote the limiting values of $Y(z)$ as
    $z$ approaches $x$ from above and below, respectively, and
    \begin{equation}
        Y_+(x)=Y_-(x)
        \begin{pmatrix}
            1 & x^\alpha e^{-Q(x)}\\
            0 & 1
        \end{pmatrix},\qquad\mbox{for $x\in(0,\infty)$.}
    \end{equation}
    \item[(c)] $Y$ has the following asymptotic behavior at
    infinity,
    \begin{equation}
        Y(z)
        \begin{pmatrix}
            z^{-n} & 0 \\
            0 & z^n
        \end{pmatrix}=I+\bigO(1/z),\qquad\mbox{as $z\to\infty$.}
    \end{equation}
    \item[(d)] $Y$ has the following behavior near $z=0$,
    \begin{equation}\label{RHP Y: d}
        Y(z)=
        \begin{cases}
            \bigO\begin{pmatrix}
                1 & z^\alpha \\
                1 & z^\alpha
            \end{pmatrix}, & \mbox{if $\alpha<0$,} \\[3ex]
            \bigO\begin{pmatrix}
                1 & \log z \\
                1 & \log z
            \end{pmatrix}, & \mbox{if $\alpha=0$,} \\[3ex]
            \bigO\begin{pmatrix}
                1 & 1 \\
                1 & 1
            \end{pmatrix}, & \mbox{if $\alpha>0$,}
        \end{cases}
    \end{equation}
    as $z\to 0, z\in\mathbb{C}\setminus[0,\infty)$.
\end{itemize}

\begin{remark}
    The $\bigO$-terms in (\ref{RHP Y: d}) are to be taken entrywise. So
    for example $Y(z)=
    \bigO\left(\begin{smallmatrix}
        1 & z^\alpha \\
        1 & z^\alpha
    \end{smallmatrix}\right)$ means that $Y_{11}=\bigO(1)$,
    $Y_{12}(z)=\bigO(z^\alpha)$, etc.
\end{remark}

\begin{remark}
    In contrast to the case of weights on the whole real line,
    considered by Deift et al.~\cite{DKMVZ2,DKMVZ1} we now have an
    extra condition (d) near the origin, cf.~\cite{Kuijlaars,KMVV}.
    This condition is used to control the behavior near the origin.
\end{remark}

The unique solution of the RH problem for $Y$, see
\cite{FokasItsKitaev} (for condition (d) see
\cite{Kuijlaars,KMVV}), is then given by,
\begin{equation}\label{definition: Y}
    Y(z)=
    \begin{pmatrix}
        \frac{1}{\gamma_n}p_n(z) & \frac{1}{\gamma_n} C(p_n w)(z) \\[1ex]
        -2\pi i\gamma_{n-1} p_{n-1}(z) & -2\pi i\gamma_{n-1}
        C(p_{n-1}w)(z)
    \end{pmatrix}, \qquad\mbox{for $z\in\mathbb{C}\setminus[0,\infty)$,}
\end{equation}
where $p_n(x)=\gamma_n x^n+\cdots$ is the $n$-th degree
orthonormal polynomial with respect to the weight $w(x)=x^\alpha
e^{-Q(x)}$, where $\gamma_n>0$ is the leading coefficient of
$p_n$, and where $Cf(z)$ is the Cauchy transform of $f$ on the
positive real line,
\[
    Cf(z)=\frac{1}{2\pi i}\int_0^\infty
    \frac{f(s)}{s-z}ds,\qquad\mbox{for $z\in\mathbb{C}\setminus[0,\infty)$.}
\]

\begin{remark}\label{remark: kernels in Y}
    From (\ref{definition: Y}), we immediately see that the kernels
    $K_n$, $W_{\I,n}$, $W_{\II,n}$ and $W_{\III,n}$,
    given by (\ref{definition: Kn}) and (\ref{definition: WI})--(\ref{definition: WIII}),
    can be written in terms of $Y$. The kernels $K_n$ and $W_{I,n}$ depend only
    on the first column of $Y$, the kernel $W_{\II,n}$ on both the
    first and the second column, and the kernel $W_{\III,n}$ only
    on the second column. Using the fact that $\det Y\equiv 1$, it is easy to check that
    \begin{equation}\label{Kn in Y}
        K_n(x,y)=x^{\frac{\alpha}{2}}e^{-\frac{1}{2}Q(x)} y^{\frac{\alpha}{2}}e^{-\frac{1}{2}Q(y)}
        \frac{1}{2\pi i(x-y)}\begin{pmatrix} 0 & 1 \end{pmatrix}
        Y^{-1}(y)Y(x) \begin{pmatrix} 1 \\ 0 \end{pmatrix},\qquad
    \mbox{$x,y\in\mathbb{R}$,}
    \end{equation}
    and
    \begin{equation}\label{Wn in Y}
        \begin{pmatrix}
                W_{\II,n}(v,u) & W_{\III,n}(u,v) \\
                -W_{\I,n}(u,v) & -W_{\II,n}(u,v)
            \end{pmatrix}
            =\frac{1}{-2\pi i\gamma_{n-1}^2 (u-v)}Y^{-1}(v) Y(u),
    \qquad\mbox{$u,v\in\mathbb{C}\setminus[0,\infty)$.}
    \end{equation}
    Note in (\ref{Kn in Y}) that  $\begin{pmatrix}0 & 1\end{pmatrix}Y^{-1}$ and
    $Y\begin{pmatrix} 1 \\ 0 \end{pmatrix}$ have an analytic continuation to $\mathbb{R}$.
\end{remark}

\medskip

Now, we will do the asymptotic analysis of the RH problem for $Y$.
As in \cite{DKMVZ2,DKMVZ1}, see also \cite{Deift}, we will use the
Deift-Zhou steepest descent method \cite{DeiftZhou}, and apply a
series of transformations $Y\mapsto U\mapsto T\mapsto S\mapsto R$
to arrive at a RH problem for $R$ with jump matrix uniformly close
to the identity matrix. Then, one can show \cite{DKMVZ1} that $R$
is also uniformly close to the identity matrix. By going back in
the series of transformations we then have the asymptotic behavior
of $Y$ in all regions of the complex plane.

The transformations are analogous to the corresponding ones in
\cite{Deift,DKMVZ2,DKMVZ1}. Yet, there are some technical
differences which come from the factor $x^\alpha$ in the weight.
As noted in the introduction, the main difference lies in the fact
that we have to construct a parametrix near the origin out of
Bessel functions. The construction of this parametrix is analogous
to the construction of the parametrix near the endpoints of the
modified Jacobi weight, as done in \cite{KMVV}, see also
\cite{Kuijlaars}.

\subsection{MRS number $\beta_n$ and rescaling: $Y\mapsto U$}

The first step in the asymptotic analysis of the RH problem for
$Y$ will be a rescaling $Y\mapsto U$. To do this rescaling, we
will use the MRS number $\beta_n$ satisfying equation
(\ref{introduction: MRS number}) above, and which will be
constructed in the following proposition for sufficiently large
$n$, cf.\ \cite[Proposition 5.2]{DKMVZ1}.

\begin{proposition}\label{proposition: MRS}
    There is $n_1\in\mathbb{N}$ such that for all $n\geq n_1$
    there exists a constant $\beta_n\in\mathbb{R}$ satisfying
    the condition,
    \begin{equation}\label{MRS condition}
        \frac{1}{2\pi}
            \int_0^{\beta_n}Q'(x)\sqrt{\frac{x}{\beta_n-x}}dx=n.
    \end{equation}
    The number $\beta_n$ has a convergent power series of the
    form
    \begin{equation}\label{power series betan}
        \beta_n=n^{1/m}\sum_{k=0}^\infty\beta^{(k)}n^{-k/m},
    \end{equation}
    with coefficients $\beta^{(k)}$ that can be expressed
    explicitly in terms of the coefficients $q_0,\ldots,q_m$
    of the polynomial $Q$. In particular, the first two
    coefficients $\beta^{(0)}$ and $\beta^{(1)}$ are,
    \begin{equation}\label{definition: beta0beta1}
        \beta^{(0)}=\bigl(\frac{1}{2}mq_mA_m\bigr)^{-1/m},
            \qquad \beta^{(1)}=-\frac{2(m-1)q_{m-1}}{m(2m-1)q_m},
    \end{equation}
    where
    \begin{equation}\label{definition: Am: bis}
        A_m = \frac{1}{\pi}\int_0^1 x^{m-1}\sqrt{\frac{x}{1-x}}dx
            = \prod_{j=1}^m \frac{2j-1}{2j}.
    \end{equation}
\end{proposition}

\begin{proof}
    The construction of $\beta_n$ such that it satisfies
    (\ref{MRS condition}) is analogous to the construction of the
    MRS numbers $\alpha_n,\beta_n$ in \cite[Proposition
    5.2]{DKMVZ1}.
    Introduce the auxiliary function,
    \begin{equation}\label{proof proposition MRS number: eq1}
        G(\beta,\varepsilon)=\frac{1}{2\pi}\int_0^1\sum_{k=0}^m k q_k
            \beta^k\varepsilon^{m-k}x^{k-1}\sqrt{\frac{x}{1-x}}dx
            = \sum_{k=0}^m \frac{1}{2} k q_k A_k \beta^k\varepsilon^{m-k}.
    \end{equation}
    An easy calculation shows that,
    \[
        G(\beta n^{-1/m},n^{-1/m}) = \frac{1}{2\pi n}\int_0^\beta
        Q'(x)\sqrt{\frac{x}{\beta-x}}dx.
    \]
    So, we need to construct $\beta_n$ such that $G(\beta_n n^{-1/m},n^{-1/m})=1$.

    Now, with $\beta^{(0)}$ defined in (\ref{definition: beta0beta1}),
    we have by (\ref{proof proposition MRS number: eq1})
    \[
        G(\beta^{(0)},0) = 1, \qquad \mbox{and}\qquad
        \left.\frac{d}{d\beta}G(\beta,0)\right|_{\beta=\beta^{(0)}}
            = m\bigl(\frac{1}{2}m q_m A_m\bigr)^{1/m}
            \neq 0.
    \]
    Therefore, using the implicit function theorem, there exists
    $\varepsilon_0>0$ and a real analytic function
    $\beta:(-\varepsilon_0,\varepsilon_0)\to\mathbb{R}$, such that
    $G(\beta(\varepsilon),\varepsilon)=1$.
    We then define $\beta_n=n^{1/m}\beta(n^{-1/m})$
    for $n$ sufficiently large such that $n^{-1/m}<\varepsilon_0$,
    and we obtain $G(\beta_n n^{-1/m},n^{-1/m})=1$. We now have constructed
    $\beta_n$ (for large enough $n$) to satisfy (\ref{MRS condition}).

    Since the function $\beta$ is analytic near 0, it follows that
    $\beta_n=n^{1/m}\beta(n^{-1/m})$ has a convergent power series of the form
    (\ref{power series betan}).
    The coefficients $\beta^{(1)},\beta^{(2)},\ldots$ can
    be expressed explicitly in terms of the coefficients $q_0,\ldots ,q_m$
    by an inductive argument using the fact that
    \[
        0=\left.\frac{d^j}{d\varepsilon^j}G(\beta(\varepsilon),\varepsilon)\right|_{\varepsilon=0}
        =\sum_{k=0}^m \frac{1}{2} k q_k A_k
        \left.\frac{d^j}{d\varepsilon^j}(\beta(\varepsilon)^k\varepsilon^{m-k})\right|_{\varepsilon=0},
        \qquad\mbox{for $j=1,2,\ldots.$}
    \]
    For example, for $j=1$ this gives $m^2q_mA_m\beta^{(1)}+(m-1) q_{m-1}A_{m-1}=0$, so
    that $\beta^{(1)}$ is given by (\ref{definition: beta0beta1}).
\end{proof}

\begin{remark}\label{remark: special case betan}
    The special case $Q(x)=q_m x^m$: In this case $\beta_n$
    takes a simple form and exists for all $n\in\mathbb{N}$.
    By a straightforward calculation one can verify that
    $\beta_n=n^{1/m} (\frac{1}{2}mq_m A_m)^{-1/m}$ solves (\ref{MRS condition}) in
    this special case.
\end{remark}

\medskip

In order to rescale the RH problem for $Y$, we also have to
introduce the rescaled field $V_n$. Define, for all $n\geq n_1$,
\begin{equation}\label{definition: Vn}
    V_n(x)=\frac{1}{n}Q(\beta_n x)=\sum_{k=0}^m \left(\frac{1}{n}q_k\beta_n^k\right) x^k.
\end{equation}
So, $V_n$ is again a polynomial of degree $m$ with real
coefficients. The coefficients of $V_n$ are $n$-dependent and have
an explicitly computable power series in $n^{-1/m}$. In particular
one has
\begin{equation}\label{expression: Vn}
    V_n(x)=\sum_{k=0}^m v_{n,k}x^k, \qquad
    v_{n,k}=\frac{1}{n}q_k\beta_n^k=\sum_{l=m-k}^\infty v_k^{(l)}n^{-l/m}.
\end{equation}

\begin{remark}\label{remark: Vn}
    We immediately see, by using Proposition \ref{proposition: MRS},
    that the leading coefficient $v_{n,m}$ has the asymptotic behavior
    $v_{n,m}=(\frac{1}{2} mA_m)^{-1}+\bigO(n^{-1/m})$ as $n\to\infty$,
    whereas the other coefficients of $V_n$ tend to zero as $n\to\infty$.
\end{remark}

\begin{remark}\label{remark: special case Vn}
    The special case $Q(x)=q_m x^m$: In this case $V_n(x)=\frac{1}{n}q_m (\beta_n x)^m$
    so that by Remark \ref{remark: special case betan},
    $V_n(x)=(\frac{1}{2} mA_m)^{-1}x^m$. Thus, $v_{n,m}=(\frac{1}{2}mA_m)^{-1}$ and $v_{n,j}=0$ for
    $j=0,\ldots ,m-1$.
\end{remark}

Now, we are ready to rescale the RH problem for $Y$. Let $\beta_n$
be the MRS number constructed in Proposition \ref{proposition:
MRS}, and define for all $n\geq n_1$,
\begin{equation}\label{definition: U}
    U(z)=
        \beta_n^{-(n+\frac{\alpha}{2})\sigma_3}
        Y(\beta_n z)\beta_n^{\ \frac{1}{2}\alpha\sigma_3},
        \qquad \mbox{for $z\in\mathbb{C}\setminus[0,\infty)$,}
\end{equation}
with $\sigma_3=\left(\begin{smallmatrix} 1 & 0 \\ 0 & -1
\end{smallmatrix}\right)$ the third Pauli matrix.
Then, it is straightforward to check, using (\ref{definition: Vn})
and the conditions of the RH problem for $Y$, that $U$ is the
unique solution of the following equivalent RH problem.

\subsubsection*{RH problem for $U$:}
\begin{itemize}
    \item[(a)] $U:\mathbb{C}\setminus[0,\infty)\to\mathbb{C}^{2\times 2}$
        is analytic.
    \item[(b)] $U_+(x)=U_-(x)
                \begin{pmatrix}
                    1 & x^\alpha e^{-nV_n(x)}\\
                    0 & 1
                    \end{pmatrix}$,
                \qquad for $x\in(0,\infty)$.
    \item[(c)] $U(z)=(I+\bigO(1/z))
                \begin{pmatrix}
                    z^n & 0 \\
                    0 & z^{-n}
                \end{pmatrix}$, \qquad as $z\to\infty$.
    \item[(d)] $U$ satisfies the same behavior near $z=0$ as $Y$ does,
        given by (\ref{RHP Y: d}).
\end{itemize}

\begin{remark}
    The RH problem for $U$ is the RH problem for orthogonal
    polynomials corresponding to the rescaled weight
    $x^\alpha e^{-nV_n(x)}$ where $V_n$ is the rescaled field (\ref{definition: Vn}).

    The rescaling with $\beta_n$ is chosen to
    ensure that the equilibrium measure $\mu_n$ of $[0,\infty)$ in
    the presence of the external field $V_n$ will be supported on
    the interval $(0,1)$.
\end{remark}

\subsection{The equilibrium measure $\mu_n$ of $[0,\infty)$
    in the external field $V_n$} \label{section: equilibrium measure}

Here, we will determine the equilibrium measure $\mu_n$ of
$[0,\infty)$ in the presence of the external field $V_n$. This
measure will be used (via its log-transform) in the next
subsection to normalize the RH problem for $U$ at infinity. The
equilibrium measure $\mu_n$ is the unique minimizer of
(\ref{minimization problem}), and is characterized by the
Euler-Lagrange variational conditions: there exists
$\ell\in\mathbb{R}$ such that
\begin{align}
    \label{variational condition 1}
    & 2\int\log|x-y|d\mu(y)-V_n(x)-\ell=0,
        \qquad \mbox{for $x\in\supp(\mu)$,} \\
    \label{variational condition 2}
    & 2\int\log|x-y|d\mu(y)-V_n(x)-\ell\leq 0,
        \qquad \mbox{for $x\in[0,\infty)$.}
\end{align}

\medskip

The construction of $\mu_n$ is as in \cite[Section 5.2]{DKMVZ1}.
It involves an analytic scalar function $h_n$ defined for all
$n\geq n_1$ as,
\begin{equation}\label{definition: hn}
    h_n(z)=\frac{1}{2\pi i}\oint_{\Gamma_z}
    \frac{y^{1/2}}{(y-1)^{1/2}}V_n'(y)\frac{dy}{y-z},\qquad
    \mbox{for $z\in\mathbb{C}\setminus[0,1]$,}
\end{equation}
where $V_n$ is the rescaled field (\ref{definition: Vn}), and
$\Gamma_z$ is a positively oriented contour containing $[0,1]$ and
$z$ in its interior. This $h_n$ has the following properties, cf.\
\cite[Proposition 5.3]{DKMVZ1}.

\begin{proposition}\label{proposition: hn}
    The function $h_n$ is a polynomial of degree $m-1$ with real
    coefficients that have an explicitly computable power series in
    $n^{-1/m}$. Furthermore, there exists $n_2\geq n_1$ and a
    constant $h_0>0$ such that $h_n(x)>h_0$ for all $n\geq n_2$
    and $x\in [0,\infty)$.
\end{proposition}

\begin{proof}
    From taking the residue at infinity in (\ref{definition: hn})
    we obtain after a straightforward calculation that $h_n$ is a polynomial of degree $m-1$
    given by,
\begin{equation}\label{proof hn polynomial: eq1}
    h_n(z)=\sum_{k=0}^{m-1}h_{n,k}z^k,\qquad
    h_{n,k}=\sum_{j=k+1}^{m} j v_{n,j} A_{j-k-1},
\end{equation}
where the $A_j$ are defined by (\ref{definition: Am: bis}) and the
$v_{n,j}$ are the coefficients of the rescaled field $V_n$. Recall
that the $v_{n,j}$ are real with an explicitly computable power
series in $n^{-1/m}$. This yields that the $h_{n,k}$ are also real
with an explicitly computable power series in $n^{-1/m}$,
\begin{equation}\label{definition: hk0}
    h_{n,k} = \sum_{j=0}^\infty h_k^{(j)}n^{-j/m},
        \qquad h_k^{(0)} = 2\frac{A_{m-1-k}}{A_m}
            = 2\prod_{j=m-k}^m \frac{2j}{2j-1}.
\end{equation}
Here, the leading order behavior $h_k^{(0)}$ of $h_{n,k}$ has been
determined by using (\ref{proof hn polynomial: eq1}) and Remark
\ref{remark: Vn}. By (\ref{definition: hk0}) there exists $n_2\geq
n_1$ and a constant $h_0>0$ such that for every $n\geq n_2$, all
coefficients $h_{n,k}$ are positive and $h_n(0)>h_0$. Therefore,
$h_n(x)\geq h_n(0)>h_0$, for all $x\in[0,\infty)$ and $n\geq n_2$.
\end{proof}

\begin{remark}\label{remark: h}
    Using (\ref{proof hn polynomial: eq1}) and (\ref{definition: hk0})
    we obtain that,
    \begin{equation}\label{asymptotics hn}
        h_n(z) = h(z) + \bigO(n^{-1/m}),
            \qquad h(z) = 2\sum_{k=0}^{m-1} \frac{A_{m-1-k}}{A_m}z^k,
    \end{equation}
    as $n\to\infty$, uniformly for $z$ in compact subsets of $\mathbb{C}$.
    Using \cite[formula 15.4.1]{AbramowitzStegun} the function $h$ can be expressed in
    terms of a hypergeometric series,
    \begin{equation}
        h(z)= \frac{4m}{2m-1}\ _2F_1(1,-m+1;-m+\frac{3}{2};z).
    \end{equation}
    Furthermore, one can check (by induction on $m$) that
    \begin{equation}\label{h(0)h(1)h'(1)}
        h(0)=\frac{4m}{2m-1},\qquad h(1)=4m,\qquad\mbox{and}\qquad h'(1)=\frac{8}{3}m(m-1).
    \end{equation}
\end{remark}

\begin{remark}\label{remark: special case hn}
    The special case $Q(x)=q_m x^m$: In this case, it follows from
    (\ref{proof hn polynomial: eq1}) and Remark \ref{remark: special case Vn}
    that $h_n=h$ with $h$ defined in (\ref{asymptotics hn}).
\end{remark}

\medskip

Now, we determine the equilibrium measure $\mu_n$ (for all $n\geq
n_2$) in terms of the polynomial $h_n$, cf.\ \cite[Proposition
5.3]{DKMVZ1}.

\begin{proposition}\label{proposition: equilibrium measure}
    Define for all $n\geq n_2$,
    \begin{equation}\label{definition: psinhat}
        \hat\psi_n(x)=\frac{1}{2\pi}\sqrt{\frac{1-x}{x}}h_n(x)\chi_{(0,1]}(x),
    \end{equation}
    where $h_n$ is given by {\rm (\ref{definition: hn})}, and
    where $\chi_{(0,1]}$ is the indicator function of the set
    $(0,1]$. By Proposition {\rm \ref{proposition: hn}} it is clear
    that $\hat\psi_n$ is non-negative. Furthermore,
    \begin{equation}\label{intpsinhat=1}
        \int_0^1 \hat\psi_n(y)dy=1,
    \end{equation}
    and there exists a constant $\ell_n\in\mathbb{R}$ such that
    \begin{align}
        \label{variational conditions: eq1}
        & 2\int \log|x-y|\hat\psi_n(y)dy-V_n(x)-\ell_n=0,\qquad \mbox{for
        $x\in[0,1]$,} \\
        \label{variational conditions: eq2}
        & 2\int \log|x-y|\hat\psi_n(y)dy-V_n(x)-\ell_n < 0,\qquad
        \mbox{for $x\in(1,\infty)$.}
    \end{align}
    So, $d\mu_n(x)=\hat\psi_n(x)dx$ is the
    equilibrium measure of $[0,\infty)$ in the external field $V_n$.
\end{proposition}

\begin{proof}
    The proof is similar to \cite[Proof of Proposition
    5.3]{DKMVZ1} and is based on the auxiliary scalar
    function
    \begin{equation}\label{definition: hulpFn}
        F_n(z)=\frac{1}{2\pi
        i}\frac{(z-1)^{1/2}}{z^{1/2}}h_n(z)-\frac{1}{2\pi i}
        V_n'(z),\qquad\mbox{for $z\in\mathbb{C}\setminus[0,1]$.}
    \end{equation}

    In order to prove that $\int_0^1 \hat\psi_n(y)dy=1$ we will determine
    two representations for the asymptotic behavior of $F_n$ at infinity and compare them to each other.
    First, from the definition (\ref{definition: hn}) of $h_n$ it
    follows that
    \[
        F_n(z)=\frac{1}{\pi i}\frac{(z-1)^{1/2}}{z^{1/2}}\frac{1}{2\pi}\int_0^1
            \sqrt{\frac{y}{1-y}}V_n'(y) \frac{dy}{y-z},
            \qquad \mbox{for $z\in\mathbb{C}\setminus[0,1]$,}
    \]
    so that, since $V_n'(y)=\frac{1}{n} \beta_n Q'(\beta_n y)$ and
    by the condition (\ref{MRS condition}) on $\beta_n$,
    \begin{equation}\label{asymptotics hulpFn}
        F_n(z)=\frac{-1}{\pi i z}\, \frac{1}{2\pi n}\int_0^{\beta_n}
        \sqrt{\frac{y}{\beta_n-y}}Q'(y)dy+\bigO(z^{-2})=\frac{-1}{\pi i
        z}+\bigO(z^{-2}),\qquad \mbox{as $z\to\infty$.}
    \end{equation}
    Next, observe that $F_{n,+}(y)-F_{n,-}(y)=2\hat\psi_n(y)$ for $y\in(0,1)$.
    Since $F_n$ is analytic in $\mathbb{C}\setminus[0,1]$ and since
    $F_n(z)=\bigO(z^{-1})$ as $z\to\infty$, see (\ref{asymptotics hulpFn}),
    a standard complex analysis argument then shows that,
    \begin{equation}\label{proof proposition equilibrium measure: eq3}
        \frac{1}{\pi i}\int_0^1\frac{\hat\psi_n(y)}{y-z}dy=
            \frac{1}{2\pi i}\int_0^1 \frac{F_{n,+}(y)-F_{n,-}(y)}{y-z}dy=F_n(z),\qquad
        \mbox{for $z\in\mathbb{C}\setminus[0,1]$.}
    \end{equation}
    Therefore, $F_n(z)=\frac{-1}{\pi i z}\int_0^1
    \hat\psi_n(y)dy+\bigO(z^{-2})$, as $z\to\infty$. Comparing this with (\ref{asymptotics hulpFn})
    we obtain $\int_0^1 \hat\psi_n(y)dy=1$.

    It now remains to prove the Euler-Lagrange variational
    conditions (\ref{variational conditions: eq1}) and
    (\ref{variational conditions: eq2}). Since $\hat\psi_n\in
    L^{3/2}(\mathbb{R})$ it follows from (\ref{proof proposition equilibrium measure: eq3})
    and \cite[Theorems 5.31 and 5.32]{RosRov} that
    \[
        \Im F_{n,+}(x)=\frac{1}{\pi}\,\PVint\frac{\hat\psi_n(y)}{x-y}dy,
            \qquad\mbox{a.e.\ for $x\in\mathbb{R}$,}
    \]
    where the integral is a Cauchy principal value integral.
    Therefore, using also (\ref{definition: hulpFn}), we obtain that
    \begin{align}\label{derivative variational}
        \nonumber
        \frac{d}{dx}\left(2\int\log|x-y|\hat\psi_n(y)dy-V_n(x)\right)
        &= 2 \PVint\frac{\hat\psi_n(y)}{x-y}dy-V_n'(x)
            = 2\pi \Im F_{n,+}(x)-V_n'(x) \\[1ex]
        &= -\Re\left((x-1)_+^{1/2}x_+^{-1/2}h_n(x)\right).
    \end{align}
    Since $h_n$ is positive on $[0,\infty)$, this yields that
    $2\int \log|x-y|\hat\psi_n(y)dy-V_n(x)$
    is constant for $x\in[0,1]$ and decreasing for
    $x\in(1,\infty)$, so that $\hat\psi_n$ satisfies conditions
    (\ref{variational conditions: eq1}) and (\ref{variational conditions: eq2}).
\end{proof}

\subsection{Normalization of the RH problem at infinity: $U\mapsto T$}

In order to normalize the RH problem for $U$ at infinity, we use
the log-transform of the equilibrium measure
$d\mu_n(y)=\hat\psi_n(y)dy$. Define, for all $n\geq n_2$ (with
$n_2$ defined in Proposition \ref{proposition: hn} above),
\begin{equation}\label{definition: gn}
    g_n(z)=\int_0^1\log(z-y)\hat\psi_n(y)dy,
        \qquad \mbox{for $z\in\mathbb{C}\setminus(-\infty,1]$,}
\end{equation}
where we take the principal branch of the logarithm, so that $g_n$
is analytic in $\mathbb{C}\setminus(-\infty,1]$.

We now give properties of $g_n$, cf.~\cite[Proposition
5.4]{DKMVZ1}, which we will need in the following. From the
definition of $g_n$ and from the Euler-Lagrange variational
conditions (\ref{variational conditions: eq1}) and
(\ref{variational conditions: eq2}) it follows that
\begin{align}
    \label{property gn: eq1a}
    & g_{n,+}(x)+g_{n,-}(x)-V_n(x)-\ell_n=0,\qquad\mbox{for
    $x\in[0,1]$,} \\[1ex]
    \label{property gn: eq1b}
    & 2g_n(x)-V_n(x)-\ell_n<0,\qquad\mbox{for $x\in(1,\infty)$.}
\end{align}
Furthermore, using $\int_0^1\hat\psi_n(y)dy=1$ we obtain
\begin{align}
    \label{property gn: eq2a}
    g_{n,+}(x)-g_{n,-}(x) &= 2\pi i,
        \qquad \mbox{for $x\in(-\infty,0)$,} \\
    \label{property gn: eq2b}
    g_{n,+}(x)-g_{n,-}(x) &= 2\pi i \int_x^1 \hat \psi_n(y)dy,
        \qquad \mbox{for $x\in[0,1]$,}
\end{align}
and
\begin{equation}\label{property gn: eq3}
    e^{ng_n(z)}=z^n+ \bigO(z^{n-1}) ,\qquad\mbox{as $z\to\infty$.}
\end{equation}

Now, we are ready to perform the transformation $U\mapsto T$.
Define, for all $n\geq n_2$, the matrix valued function $T$ as,
\begin{equation}\label{definition: T}
    T(z)=e^{-\frac{1}{2}n\ell_n\sigma_3}U(z)e^{-ng_n(z)\sigma_3}e^{\frac{1}{2}n\ell_n\sigma_3},
        \qquad\mbox{for $z\in\mathbb{C}\setminus \mathbb{R}$,}
\end{equation}
where $\ell_n$ is the constant that appears in the Euler-Lagrange
variational conditions (\ref{variational conditions: eq1}) and
(\ref{variational conditions: eq2}). Note that, by (\ref{property
gn: eq2a}), the function $e^{ng_n}$ has no jumps across
$(-\infty,0)$, so that $T$ has an analytic continuation to
$\mathbb{C}\setminus[0,\infty)$. It is then straightforward to
check, using (\ref{property gn: eq1a}), (\ref{property gn: eq3})
and the conditions of the RH problem for $U$, that $T$ is the
unique solution of the following equivalent RH problem.

\subsubsection*{RH problem for $T$:}
\begin{itemize}
    \item[(a)] $T:\mathbb{C}\setminus[0,\infty)\to\mathbb{C}^{2\times 2}$ is analytic.
    \item[(b)] $T_+(x)=T_-(x)v_T(x)$ for $x\in(0,\infty)$, with
        \begin{equation}\label{definition: vT}
            v_T(x)=
            \begin{cases}
                \begin{pmatrix}
                    e^{-n(g_{n,+}(x)-g_{n,-}(x))} & x^\alpha \\
                    0 & e^{n(g_{n,+}(x)-g_{n,-}(x))}
                \end{pmatrix}, & \mbox{for $x\in(0,1)$,} \\[4ex]
                \begin{pmatrix}
                    1 & x^\alpha e^{n(2g_n(x)-V_n(x)-\ell_n)} \\
                    0 & 1
                \end{pmatrix}, & \mbox{for $x\in[1,\infty)$.}
            \end{cases}
        \end{equation}
    \item[(c)] $T(z)=I+\bigO(1/z)$,\qquad as $z\to\infty$.
    \item[(d)] $T$ satisfies the same behavior near $z=0$ as $Y$ and $U$ do,
        given by (\ref{RHP Y: d}).
\end{itemize}

\begin{remark}
    From (\ref{property gn: eq2b}) we see that the diagonal entries
    of $v_T$ on $(0,1)$ are rapidly oscillating for large
    $n$. From (\ref{property gn: eq1b}) the jump matrix $v_T$ on
    $(1,\infty)$ converges exponentially fast to the identity matrix
    as $n\to\infty$.
\end{remark}

\subsection{Opening of the lens: $T\mapsto S$}
\label{subsection: opening lens}

Now, we will transform the oscillatory diagonal entries of the
jump matrix $v_T$ on $(0,1)$ into exponentially decaying
off-diagonal entries. This lies at the heart of the Deift-Zhou
steepest descent method \cite{DeiftZhou}, and this step is
referred to as the opening of the lens.

\medskip

In order to perform the transformation $T\mapsto S$ we will
introduce scalar functions $\psi_n$ and $\xi_n$. Define for all
$n\geq n_2$,
\begin{equation}\label{definition: psin}
    \psi_n(z)=\frac{1}{2\pi
    i}\frac{(z-1)^{1/2}}{z^{1/2}}h_n(z),\qquad\mbox{for
    $z\in\mathbb{C}\setminus[0,1]$,}
\end{equation}
with principal branches of powers. So, the $+$boundary value of
$\psi_n$ on $(0,1)$ is precisely the density $\hat \psi_n$ of the
equilibrium measure $\mu_n$. In particular we have,
\begin{equation}\label{property: psin}
    \psi_{n,+}(x)=-\psi_{n,-}(x)=\hat\psi_n(x),\qquad\mbox{for
    $x\in(0,1)$.}
\end{equation}
Now, define for all $n\geq n_2$,
\begin{equation}\label{definition: xin}
    \xi_n(z)=-\pi i\int_1^z \psi_n(y)dy,
        \qquad\mbox{for $z\in\mathbb{C}\setminus(-\infty,1]$,}
\end{equation}
where the path of integration does not cross the real axis. Note
that in \cite[equation (5.34)]{DKMVZ1} a function $\xi_n$ is
defined which is analytic through the support of the equilibrium
measure. We find it more convenient to define $\xi_n$ with a
branch cut along $(0,1)$. The definition of $\xi_n$ in
(\ref{definition: xin}) and \cite[equation (5.34)]{DKMVZ1} also
differs by a factor 2.

The important feature of the function $\xi_n$ is that, by
(\ref{property: psin}) and (\ref{property gn: eq2b}), $\xi_{n,+}$
and $\xi_{n,-}$ are purely imaginary on $(0,1)$ and satisfy,
\begin{equation}\label{property xin: eq1}
    2\xi_{n,+}(x)=-2\xi_{n,-}(x)=2\pi i\int_x^1\hat\psi_n(y)dy=g_{n,+}(x)-g_{n,-}(x),
                \qquad \mbox{for $x\in(0,1)$.}
\end{equation}
So, $2\xi_n$ and $-2\xi_n$ provide analytic extensions of
$g_{n,+}-g_{n,-}$ into the upper half-plane and lower half-plane,
respectively. On $\mathbb{R}\setminus[0,1)$, $\xi_n$ satisfies
\begin{align}\label{property xin: eq2}
    & \xi_{n,+}(x)-\xi_{n,-}(x)=2\pi i,\qquad\mbox{for
        $x\in(-\infty,0)$,} \\[1ex]
    \label{property xin: eq3}
    & 2\xi_n(x)=2g_n(x)-V_n(x)-\ell_n,
        \qquad \mbox{for $x\in[1,\infty)$.}
\end{align}
Here, equation (\ref{property xin: eq2}) follows from
(\ref{intpsinhat=1}). Equation (\ref{property xin: eq3}) follows
from the fact that the function $2g_n-V_n-2\xi_n-\ell_n$ is
constant on $[1,\infty)$, by (\ref{derivative variational}), and 0
in $x=1$, by (\ref{variational conditions: eq1}). Further, we can
prove the existence of a $\delta_1>0$ such that for all $n\geq
n_2$, cf.~\cite[Proposition 5.4]{DKMVZ1}
\begin{equation}\label{property xin: eq4}
    \Re \xi_n(z)>0,\qquad\mbox{for $0<|\Im z|<\delta_1$ and $0<\Re z<1$.}
\end{equation}
We also need an estimate for $\xi_n$ on $(1,\infty)$. From the
fact that $h_n(x)>h_0$ for all $n\geq n_2$ and $x\in[0,\infty)$,
see Proposition \ref{proposition: hn}, we obtain
\begin{equation}\label{property xin: eq5}
    \xi_n(x)<-\frac{1}{3}h_0 \frac{(x-1)^{3/2}}{x^{1/2}},
        \qquad\mbox{for $x\in(1,\infty)$ and $n\geq n_2$.}
\end{equation}

\begin{remark}\label{remark: xin and ln}
    As in \cite[Proof of Proposition 5.4 (v) and (vi)]{DKMVZ1} we
    can determine $\xi_n$ as well as the constant $\ell_n$
    (which appears in the variational conditions)
    explicitly in terms of the coefficients $v_{n,k}$ of $V_n$.
    With $\arccos$ defined as an analytic function on
    $\mathbb{C}\setminus\left((-\infty,-1]\cup [1,\infty)\right)$,
    as described in Remark \ref{remark: arccos}, one has
    \begin{equation}
        \xi_n(z)=\mp i \left(\frac{1}{2}H_n(z)z^{1/2}(1-z)^{1/2}-2\arccos
            z^{1/2}\right), \qquad\mbox{for $z\in\mathbb{C}_\pm$.}
    \end{equation}
    where
    \begin{equation}\label{definition: Hn}
        H_n(z) = \sum_{k=0}^{m-1}
            \left(\sum_{j=k+1}^m v_{n,j} A_{j-k-1} \right)
            z^k= \frac{1}{m}h(z)+\bigO(n^{-1/m}),
    \end{equation}
    as $n\to\infty$, uniformly for $z$ in compact subsets of
    $\mathbb{C}$. Here, the leading order behavior of $H_n$ has
    been determined by using  Remark \ref{remark: Vn} and (\ref{asymptotics
    hn}).

    Further, one can also verify that $\ell_n=-\sum_{k=0}^m v_{n,k} A_k -4\log 2$.
    Then, since the coefficients $v_{n,k}$ have an explicitly computable power series in
    $n^{-1/m}$, so have the constants $\ell_n$,
    \begin{equation}\label{asymptotics ln}
        \ell_n=\sum_{j=0}^\infty \ell^{(j)}n^{-j/m},
            \qquad \ell^{(0)}= -\frac{2}{m}-4\log 2,
    \end{equation}
    where the leading order behavior $\ell^{(0)}$ of $\ell_n$ has been determined
    by using Remark \ref{remark: Vn}.
\end{remark}

\begin{remark}\label{remark: special case Hn and ln}
    The special case $Q(x)=q_m x^m$: In this case, it follows from
    the previous remark together with Remark \ref{remark: special case Vn} that
    $H_n(z)=\frac{1}{m}h(z)$ and that $\ell_n=-\frac{2}{m}-4\log
    2$.
\end{remark}

\medskip

Inserting (\ref{property xin: eq1}) and (\ref{property xin: eq3})
into (\ref{definition: vT}), the jump matrix $v_T$ for $T$ can be
written in terms of the scalar function $\xi_n$ as,
\begin{equation}\label{definition: vT: bis}
    v_T(x)=\begin{cases}
        \begin{pmatrix}
            e^{-2n\xi_{n,+}(x)} & x^\alpha \\
            0 & e^{-2n\xi_{n,-}(x)}
        \end{pmatrix}, & \mbox{for $x\in(0,1)$,} \\[3ex]
        \begin{pmatrix}
            1 & x^\alpha e^{2n\xi_n(x)} \\
            0 & 1
        \end{pmatrix}, & \mbox{for $x\in[1,\infty)$.}
    \end{cases}
\end{equation}
A simple calculation, using the fact that
$\xi_{n,+}(x)+\xi_{n,-}(x)=0$ for $x\in(0,1)$, see (\ref{property
xin: eq1}), then shows that $v_T$ has on the interval $(0,1)$ the
following factorization,
\begin{equation}\label{factorization vT}
    v_T(x)=
    \begin{pmatrix}
        1 & 0 \\
        x^{-\alpha} e^{-2n\xi_{n,-}(x)} & 1
    \end{pmatrix}
    \begin{pmatrix}
        0 & x^\alpha \\
        -x^{-\alpha} & 0
    \end{pmatrix}
    \begin{pmatrix}
        1 & 0 \\
        x^{-\alpha} e^{-2n\xi_{n,+}(x)} & 1
    \end{pmatrix},\quad\mbox{for $x\in(0,1)$,}
\end{equation}
and the opening of the lens is based on this factorization.

\begin{figure}[t]
\begin{center}
    \setlength{\unitlength}{1mm}
    \begin{picture}(100,31)(0,10)
        \put(10,25){\thicklines\circle*{.8}}
        \put(60,25){\thicklines\circle*{.8}}
        \put(10,25){\line(1,0){80}}
        \put(37,25){\thicklines\vector(1,0){.0001}} \put(77,25){\thicklines\vector(1,0){.0001}}
        \qbezier(10,25)(35,47)(60,25) \put(37,36){\thicklines\vector(1,0){.0001}}
        \qbezier(10,25)(35,3)(60,25) \put(37,14){\thicklines\vector(1,0){.0001}}
        \put(24,37){\small $\Sigma_1$}
        \put(24,26.5){\small $\Sigma_2$}
        \put(24,17.5){\small $\Sigma_3$}
        \put(82,26.5){\small $\Sigma_4$}
        \put(8,26.5){\small 0}
        \put(60.5,26.5){\small 1}
    \end{picture}
    \caption{The lens shaped contour $\Sigma_S=\bigcup_{j=1}^4\Sigma_j$
        oriented from the left to the right.}
    \label{figure: opening of the lens}
\end{center}
\end{figure}
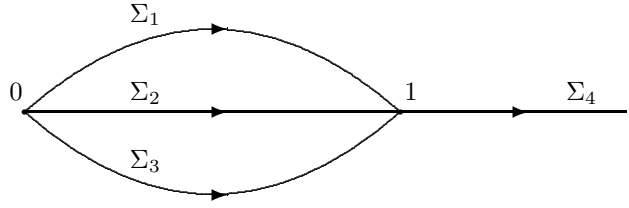

Now, we are ready to do the transformation $T\mapsto S$. Let
$\Sigma_S=\cup_{j=1}^4\Sigma_j$ be the oriented lens shaped
contour as shown in Figure \ref{figure: opening of the lens}. The
precise form of the lens (in fact of the lips $\Sigma_1$ and
$\Sigma_3$) is not yet defined but for now we assume that it will
be contained in the region where (\ref{property xin: eq4}) holds.
We will define the contour $\Sigma_S$ explicitly in the beginning
of Section \ref{subsection: final transformation}, depending on
$n$ and on certain parameters $\delta$ and $\nu$. Define, for all
$n\geq n_2$, an analytic matrix valued function $S$ on
$\mathbb{C}\setminus\Sigma_S$ as,
\begin{equation}\label{definition: S}
    S(z)=\begin{cases}
        T(z), & \mbox{for $z$ outside the lens,} \\[2ex]
        T(z)\begin{pmatrix}
            1 & 0 \\
            -z^{-\alpha}e^{-2n\xi_n(z)} & 1
        \end{pmatrix}, & \mbox{for $z$ in the upper part of the
        lens,} \\[3ex]
        T(z) \begin{pmatrix}
            1 & 0 \\
            z^{-\alpha} e^{-2n\xi_n(z)} & 1
        \end{pmatrix}, & \mbox{for $z$ in the lower part of the
        lens.}
    \end{cases}
\end{equation}
With the upper part of the lens we mean the region between
$\Sigma_1$ and $\Sigma_2$, and with the lower part of the lens the
region between $\Sigma_2$ and $\Sigma_3$.

One can easily check, using (\ref{definition: vT: bis}),
(\ref{factorization vT}) and the conditions of the RH problem for
$T$, that $S$ satisfies the following RH problem.

\subsubsection*{RH problem for $S$:}
\begin{itemize}
    \item[(a)] $S:\mathbb{C}\setminus\Sigma_S\to\mathbb{C}^{2\times
    2}$ is analytic.
    \item[(b)] $S_+(z)=S_-(z)v_S(z)$ for $z\in\Sigma_S$, with
    \begin{equation}\label{definition: vS}
        v_S(z)=\begin{cases}
            \begin{pmatrix}
                1 & 0 \\
                z^{-\alpha}e^{-2n\xi_n(z)} & 1
            \end{pmatrix}, & \mbox{for
            $z\in\Sigma_1\cup\Sigma_3$,} \\[3ex]
            \begin{pmatrix} 0 & z^\alpha \\
            -z^{-\alpha} & 0
            \end{pmatrix}, & \mbox{for $z\in\Sigma_2=(0,1)$,} \\[3ex]
            \begin{pmatrix}
                1 & z^\alpha
                e^{2n\xi_n(z)} \\ 0 &
                1
            \end{pmatrix}, & \mbox{for $z\in\Sigma_4=(1,\infty)$.}
        \end{cases}
    \end{equation}
    \item[(c)] $S(z)=I+\bigO(1/z)$,\qquad as $z\to\infty$.
\end{itemize}

\begin{remark}\label{remark: S}
Note that by (\ref{property xin: eq4}) and (\ref{property xin:
eq5}) the jump matrix $v_S$ on $\Sigma_1,\Sigma_3$ and $\Sigma_4$
converges exponentially fast (as $n\to\infty$) to the identity
matrix.
\end{remark}

\subsection{Parametrix $P^{(\infty)}$ for the outside region}
    \label{section: parametrix outside region}

From Remark \ref{remark: S} we expect that the leading order
asymptotics of $Y$ will be determined by a solution
$P^{(\infty)}$, which will be refered to as the parametrix for the
outside region, of the following RH problem.

\subsubsection*{RH problem for $P^{(\infty)}$:}
    \begin{itemize}
        \item[(a)] $P^{(\infty)}:\mathbb{C}\setminus[0,1]\to
        \mathbb{C}^{2\times 2}$ is analytic.
        \item[(b)] $P^{(\infty)}_+(x)=P^{(\infty)}_-(x)
            \begin{pmatrix}
                0 & x^\alpha \\
                -x^{-\alpha} & 0
            \end{pmatrix}$,\qquad for $x\in(0,1)$.
        \item[(c)] $P^{(\infty)}(z)=I+\bigO(1/z)$,\qquad as
        $z\to\infty$.
    \end{itemize}

As in \cite{Kuijlaars,KMVV,KV2}, we will construct a solution of
this RH problem in terms of the Szeg\H{o} function $D$ associated
with $x^\alpha$ on $(0,1)$. This is a scalar function which is
analytic and non-zero in $\mathbb{C}\setminus[0,1]$, that
satisfies $D_+(x)D_-(x)=x^\alpha$ for $x\in(0,1)$, and which will
not vanish at infinity. One  can easily check that $D$ is given
by,
\begin{equation}\label{definition: D}
    D(z)=\frac{z^{\alpha/2}}{\varphi(z)^{\alpha/2}},
    \qquad \mbox{for $z\in\mathbb{C}\setminus[0,1]$,}
\end{equation}
with principal branches of powers, and where $\varphi$ is the
conformal map from $\mathbb{C}\setminus[0,1]$ onto the exterior of
the unit circle, cf.\ (\ref{definition: varphi: eq1}),
\begin{equation}\label{definition: varphi}
    \varphi(z)=2(z-1/2)+2z^{1/2}(z-1)^{1/2},\qquad\mbox{for
    $z\in\mathbb{C}\setminus[0,1]$.}
\end{equation}
Since $\varphi(z)=4z+\bigO(1)$ as $z\to\infty$, we have
\begin{equation}
    \lim_{z\to\infty}D(z)=2^{-\alpha}.
\end{equation}

The important feature of the Szeg\H{o} function is that the
transformed matrix valued function
$2^{\alpha\sigma_3}P^{(\infty)}D^{\sigma_3}$ will satisfy
conditions (a) and (c) of the RH problem, and
that it will have the jump matrix $\left(\begin{smallmatrix}0 & 1\\
-1 & 0 \end{smallmatrix}\right)$ on $(0,1)$. Then it is well
known, see for example \cite{Deift,DKMVZ1}, that $P^{(\infty)}$ is
given by,
\begin{equation}\label{definition: Pinfinity}
    P^{(\infty)}(z)=2^{-\alpha\sigma_3}
        \begin{pmatrix}
            \frac{a(z)+a(z)^{-1}}{2} & \frac{a(z)-a(z)^{-1}}{2i}
            \\[1ex]
            \frac{a(z)-a(z)^{-1}}{-2i} & \frac{a(z)+a(z)^{-1}}{2}
        \end{pmatrix}
        D(z)^{-\sigma_3}, \qquad\mbox{for
        $z\in\mathbb{C}\setminus[0,1]$,}
\end{equation}
with
\begin{equation}\label{definition: a}
    a(z)=\frac{(z-1)^{1/4}}{z^{1/4}}, \qquad\mbox{for
        $z\in\mathbb{C}\setminus[0,1]$.}
\end{equation}

\begin{remark}\label{remark: detPinfinity=1}
    Since $P^{(\infty)}$ is a product of three matrices all with
    determinant one, we have that $\det P^{(\infty)}\equiv 1$.
    Further, for later reference, note that $P^{(\infty)}(z)
    z^{\frac{1}{2}\alpha\sigma_3}=\bigO(z^{-1/4})$ as $z\to 0$.
\end{remark}

Before we can do the final transformation $S\mapsto R$ we need to do
a local analysis near 0 and 1 since the jump
matrices for $S$ and $P^{(\infty)}$ are not uniformly close
to each other in the neighborhood of these points.

\subsection{Parametrix $P_n$ near the endpoint 1}
    \label{section: parametrix near 1}

In this subsection, we will construct inside the disk
$U_{\delta_2}=\{z\in\mathbb{C}:|z-1|<\delta_2\}$ with center 1 and
radius $\delta_2>0$ (sufficiently small and which will be
determined as part of the problem in Proposition \ref{proposition:
phin} below), a $2\times 2$ matrix valued function $P_n$ that
satisfies the following conditions.

\subsubsection*{RH problem for $P_n$:}

\begin{itemize}
\item[(a)] $P_n:U_{\delta_2}\setminus\Sigma_S\to\mathbb{C}^{2\times 2}$
    is analytic.
\item[(b)] $P_{n,+}(z)=P_{n,-}(z)v_S(z)$ for $z\in\Sigma_S\cap U_{\delta_2}$,
    with $v_S$ the jump matrix (\ref{definition: vS}) for $S$.
\item[(c)] $P_n(z)P^{(\infty)}(z)^{-1}=I+\bigO(1/n)$ as $n\to\infty$,
    uniformly for $z$ on the boundary $\partial U_\delta$ of the
    disk $U_\delta$ and for $\delta$ in compact subsets of $(0,\delta_2)$.
\end{itemize}

The construction of $P_n$ has many similarities to the analogous
construction carried out in \cite{DKMVZ2,DKMVZ1}, see also
\cite{Deift} for an excellent exposition, and will be done using
Airy functions. It typically involves three steps. First, we will
construct a matrix valued function that satisfies conditions (a)
and (b) of the RH problem for $P_n$. In order to do this we will
transform, in the first step, this RH problem into a RH problem
for $P_n^{(1)}$ with constant jump matrices and construct, in the
second step, a solution of the latter RH problem. Afterwards, we
will take in the third step also the matching condition (c) into
account.

\subsubsection*{Step 1: Transformation to constant jump matrices}

In order to transform to constant jump matrices, we seek the
parametrix $P_n$ near 1 in the following form,
\begin{equation}\label{definition: Pn}
    P_n(z) = E_n(z) P_n^{(1)}(z) e^{-n\xi_n(z)\sigma_3}
        z^{-\frac{1}{2}\alpha\sigma_3},
        \qquad \mbox{for $z\in U_{\delta_2}\setminus\Sigma_S$,}
\end{equation}
with $E_n$ an invertible analytic matrix valued function in
$U_{\delta_2}$, to be determined in the last step to ensure that
the matching condition of the RH problem for $P_n$ is satisfied.

The reader can easily verify, using $\xi_{n,+}(x)+\xi_{n,-}(x)=0$
for $x\in(0,1)$, that if $P_n^{(1)}$ is analytic in
$U_{\delta_2}\setminus \Sigma_S$ with jump relations,
\begin{equation}\label{jump relations Pn1}
    P_{n,+}^{(1)}(z) =
    \begin{cases}
        P_{n,-}^{(1)}(z)
            \begin{pmatrix}
                1 & 0 \\
                1  & 1
            \end{pmatrix},
            & \mbox{for $z\in(\Sigma_1\cup\Sigma_3)\cap U_{\delta_2}$,}
        \\[3ex]
        P_{n,-}^{(1)}(z)
            \begin{pmatrix}
                0 & 1 \\
                -1  & 0
            \end{pmatrix},
            & \mbox{for $z\in\Sigma_2\cap U_{\delta_2}=(1-\delta_2,1)$,}
        \\[3ex]
        P_{n,-}^{(1)}(z)
            \begin{pmatrix}
                1 & 1 \\
                0 & 1
            \end{pmatrix},
            & \mbox{for $z\in\Sigma_4\cap U_{\delta_2}=(1,1+\delta_2)$,}
    \end{cases}
\end{equation}
then $P_n$ defined by (\ref{definition: Pn}) satisfies conditions
(a) and (b) of the RH problem for $P_n$. In the next step we will
determine $P_n^{(1)}$ to satisfy these conditions.

\subsubsection*{Step 2: Determine $P_n^{(1)}$ explicitly}

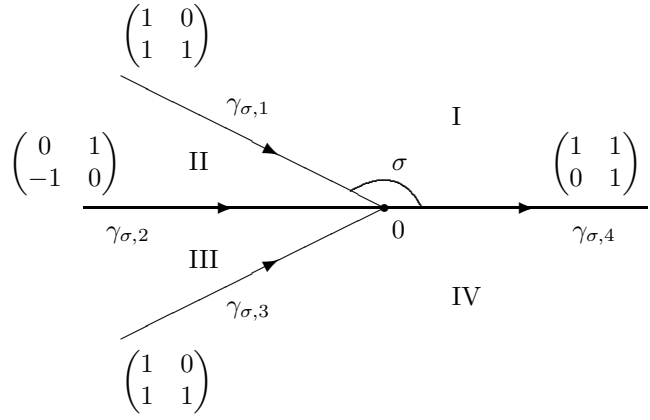
\begin{figure}[t]
    \begin{center}
    \setlength{\unitlength}{1mm}
    \begin{picture}(100,57)(0,-3)
        \put(29,38){\small $\gamma_{\sigma,1}$}
        \put(13,21){\small $\gamma_{\sigma,2}$}
        \put(29,11){\small $\gamma_{\sigma,3}$}
        \put(75,21){\small $\gamma_{\sigma,4}$}

        \put(59,36){\small $\I$}
        \put(24,30){\small $\II$}
        \put(24,17){\small $\III$}
        \put(59,12){\small $\IV$}

        \put(15,47){\small
            $\begin{pmatrix}
                1 & 0 \\
                1 & 1
            \end{pmatrix}$}
        \put(0,30){\small
            $\begin{pmatrix}
                0 & 1 \\
                -1 & 0
            \end{pmatrix}$}
        \put(15,1){\small
            $\begin{pmatrix}
                1 & 0 \\
                1 & 1
            \end{pmatrix}$}
        \put(72,30){\small
            $\begin{pmatrix}
                1 & 1 \\
                0 & 1
            \end{pmatrix}$}

        \put(50,25){\thicklines\circle*{.9}}
        \put(51,21){\small 0}

        \put(50,25){\line(-2,1){35}} \put(36,32){\thicklines\vector(2,-1){.0001}}
        \put(50,25){\line(-2,-1){35}} \put(36,18){\thicklines\vector(2,1){.0001}}
        \put(50,25){\line(-1,0){40}} \put(30,25){\thicklines\vector(1,0){.0001}}
        \put(50,25){\line(1,0){35}} \put(70,25){\thicklines\vector(1,0){.0001}}

        \qbezier(55,25)(52,31)(45.523,27.236)
        \put(51,30){\small $\sigma$}
    \end{picture}
    \caption{The oriented contour $\gamma_\sigma$ and the jump
        matrix $v_1$ for $\Psi$ on $\gamma_\sigma$. The four
        straight rays $\gamma_{\sigma,1},\ldots ,\gamma_{\sigma,4}$ divide
        the complex plane into four regions $\I,\II, \III$ and $\IV$.}
        \label{figure: contour gamma}
    \end{center}
\end{figure}

The construction of $P_n^{(1)}$ is based upon an auxiliary RH
problem for $\Psi$ in the $\zeta$-plane with jumps on the oriented
contour $\gamma_{\sigma}$, shown in Figure \ref{figure: contour
gamma}, consisting of four straight rays
\[
    \gamma_{\sigma,1}:\arg\zeta=\sigma,
        \qquad \gamma_{\sigma,2}:\arg\zeta=\pi,
        \qquad \gamma_{\sigma,3}:\arg\zeta=-\sigma,
        \qquad \gamma_{\sigma,4}:\arg\zeta=0,
\]
with $\sigma\in(\frac{\pi}{3},\pi)$. These four rays divide the
complex plane into four regions $\I,\II,\III$ and $\IV$, also
shown in Figure \ref{figure: contour gamma}. The RH problem for
$\Psi$ is the following, cf.~\cite{Deift,DKMVZ1}.

\subsubsection*{RH problem for $\Psi$:}
\begin{itemize}
    \item[(a)] $\Psi:\mathbb{C}\setminus\gamma_\sigma\to\mathbb{C}^{2\times
    2}$ is analytic.
    \item[(b)] $\Psi_+(\zeta)=\Psi_-(\zeta)v_1(\zeta)$ for
    $\zeta\in\gamma_\sigma$, where $v_1$ is the piecewise
    constant matrix valued function on $\gamma_\sigma$ defined as shown in Figure
    \ref{figure: contour gamma}, i.e.~$v_1(\zeta)=
    \left(\begin{smallmatrix} 1 & 0 \\ 1 & 1\end{smallmatrix}\right)$ for $\zeta\in\gamma_{\sigma,1}$
    and so on. This means that $\Psi$
    has the same jumps on $\gamma_\sigma$ as
    $P_n^{(1)}$ on $\Sigma_S\cap U_{\delta_2}$, see (\ref{jump relations Pn1}).
    \item[(c)] $\Psi$ has the following asymptotic behavior at
    infinity,
    \begin{multline}\label{asymptotics: Psi}
        \Psi(\zeta)\sim\zeta^{-\frac{\sigma_3}{4}}\frac{1}{\sqrt 2}
        \begin{pmatrix}
            1 & 1 \\
            -1 & 1
        \end{pmatrix}\\[1ex]
        \times\,
        \left[
            I+\sum_{k=1}^\infty\frac{1}{2}\left(\frac{2}{3}\zeta^{3/2}\right)^{-k}
            \begin{pmatrix}
                (-1)^k(s_k+t_k) & s_k-t_k \\
                (-1)^k(s_k-t_k) & s_k+t_k
            \end{pmatrix}
       \right]
        e^{-\frac{\pi i}{4}\sigma_3}e^{-\frac{2}{3}\zeta^{3/2}\sigma_3},
    \end{multline}
    as $\zeta\to\infty$, uniformly for
    $\zeta\in\mathbb{C}\setminus\gamma_\sigma$ and $\sigma$ in
    compact subsets of $(\frac{\pi}{3},\pi)$. Here,
    \begin{equation}\label{definition: sktk}
        s_k=\frac{\Gamma(3k+1/2)}{54^k k! \Gamma(k+1/2)},\qquad
        t_k=-\frac{6k+1}{6k-1}s_k,\qquad\mbox{for $k\geq 1$.}
    \end{equation}
\end{itemize}

\medskip

It is well-known, see for example \cite{Deift,DKMVZ1}, that
$\Psi=\Psi^\sigma$ (we suppress $\sigma$ in the notation for
brevity) defined by,
\begin{equation}\label{definition: Psi}
    \Psi(\zeta)=\sqrt{2\pi}e^{-\frac{\pi i}{12}}\times
    \begin{cases}
        \begin{pmatrix}
            \Ai(\zeta) & \Ai(\omega^2\zeta) \\
            \Ai'(\zeta) & \omega^2\Ai'(\omega^2\zeta)
        \end{pmatrix}e^{-\frac{\pi i}{6}\sigma_3}, & \mbox{for $\zeta\in\I$,} \\[3ex]
        \begin{pmatrix}
            \Ai(\zeta) & \Ai(\omega^2\zeta) \\
            \Ai'(\zeta) & \omega^2\Ai'(\omega^2\zeta)
        \end{pmatrix}e^{-\frac{\pi i}{6}\sigma_3}
        \begin{pmatrix}
            1 & 0 \\
            -1 & 1
        \end{pmatrix}, & \mbox{for $\zeta\in\II$,} \\[3ex]
        \begin{pmatrix}
            \Ai(\zeta) & -\omega^2\Ai(\omega\zeta) \\
            \Ai'(\zeta) & -\Ai'(\omega\zeta)
        \end{pmatrix}e^{-\frac{\pi i}{6}\sigma_3}
        \begin{pmatrix}
            1 & 0 \\
            1 & 1
        \end{pmatrix}, & \mbox{for $\zeta\in\III$,} \\[3ex]
        \begin{pmatrix}
            \Ai(\zeta) & -\omega^2\Ai(\omega\zeta) \\
            \Ai'(\zeta) & -\Ai'(\omega\zeta)
        \end{pmatrix}e^{-\frac{\pi i}{6}\sigma_3}, & \mbox{for $\zeta\in\IV$,}
    \end{cases}
\end{equation}
with $\omega=e^{\frac{2\pi i}{3}}$ and $\Ai$ the Airy function,
solves the RH problem for $\Psi$. See for example
\cite{AbramowitzStegun} for definitions and properties of the Airy
functions.

\begin{remark}\label{remark: detPsi=1}
    Using \cite[formulae 10.4.11 and 10.4.12]{AbramowitzStegun} we
    have $\det \Psi\equiv 1$. This can also be seen from a
    different point of view as follows. From condition (b) of the RH problem
    for $\Psi$ and from the fact that the Airy function
    remains bounded near 0 it follows that $\det \Psi$ is entire.
    Further, from condition (c) we have $\det\Psi\to
    1$ as $z\to\infty$. Using Liouville's theorem we then indeed obtain
    $\det\Psi\equiv 1$.
\end{remark}

\medskip

The idea is now to construct $P_n^{(1)}$ out of $\Psi$ as
$P_n^{(1)}(z)=\Psi(f_n(z))$ for appropriate biholomorphic maps
$f_n:U_{\delta_2}\to f_n(U_{\delta_2})$ with $f_n(1)=0$. We will
choose these biholomorphic maps to compensate for the factor
$e^{-n\xi_n(z)\sigma_3}$ in (\ref{definition: Pn}). So, by the
asymptotic behavior (\ref{asymptotics: Psi}) of $\Psi$ at
infinity, we need,
\begin{equation}\label{property fn}
    -\frac{2}{3}f_n(z)^{3/2}=n\xi_n(z),\qquad\mbox{for $z\in
    U_{\delta_2}\setminus(-\infty,1]$.}
\end{equation}
Note that, by (\ref{definition: xin}), this is precisely equation
(\ref{statement of results: fn}). The construction of these
biholomorphic maps is analogous as in \cite[Section 7.1]{DKMVZ1}
and we define, for all $n\geq n_2$,
\begin{equation}\label{definition: fn}
    f_n(z)=n^{2/3}\phi_n(z),\qquad\mbox{for $z\in U_{\delta_2}$,}
\end{equation}
with $\phi_n$ defined in the following proposition,
cf.~\cite[Proposition 7.3]{DKMVZ1}.

\begin{proposition}\label{proposition: phin}
    There exists $\delta_2>0$ such that for all $n\geq n_2$ there
    are biholomorphic maps
    $\phi_n:U_{\delta_2}\to\phi_n(U_{\delta_2})$ satisfying:
    \begin{itemize}
        \item[1.] There exists a constant $c_0>0$ such that for all
            $z\in U_{\delta_2}$ and all $n\geq n_2$ the derivative of
            $\phi_n$ can be estimated by: $c_0<|\phi_n'(z)|<1/c_0$ and
            $|\arg\phi_n'(z)|<\pi/15$.
        \item[2.] $\phi_n(U_{\delta_2} \cap \mathbb{R})
            = \phi_n(U_{\delta_2}) \cap \mathbb{R}$, and
            $\phi_n(U_{\delta_2} \cap \mathbb{C}_\pm)
            = \phi_n(U_{\delta_2})\cap\mathbb{C}_\pm$.
        \item[3.] $-\frac{2}{3}\phi_n(z)^{3/2}=\xi_n(z)$ for $z\in
            U_{\delta_2}\setminus(-\infty,1]$.
    \end{itemize}
\end{proposition}

\begin{proof}
    Define, for all $n\geq n_2$, the auxiliary function,
    \begin{equation}\label{definition: phinhat}
        \hat\phi_n(z) =
        \frac{2}{h_n(1)}\frac{-\frac{3}{2}\xi_n(z)}{(z-1)^{3/2}},
        \qquad\mbox{for $z\in\mathbb{C}\setminus(-\infty,1]$.}
    \end{equation}
    Note that by (\ref{property xin: eq1}) the function
    $\hat\phi_n$ has no jumps across $(0,1)$, so that $\hat\phi_n$
    has an analytic continuation to
    $\mathbb{C}\setminus((-\infty,0]\cup\{1\})$. From (\ref{definition: xin})
    and (\ref{definition: psin}) it follows that,
    \begin{equation}\label{proof: proposition phin: eq2}
        \hat\phi_n(z)
            = 1+\frac{3}{2}(z-1)^{-3/2}\frac{1}{h_n(1)}
            \int_1^z\left(\frac{h_n(s)}{s^{1/2}}-h_n(1)\right)(s-1)^{1/2}ds.
    \end{equation}
    Using Cauchy's theorem and the fact that $h_n$ is uniformly bounded in
    compact subsets of $\mathbb{C}$, which follows from equation
    (\ref{asymptotics hn}), there exists a constant $c>0$ such that for all
    $n\geq n_2$ and $|s-1|\leq 1/4$,
    \begin{align*}
        \left|\frac{h_n(s)}{s^{1/2}}-h_n(1)\right|
            &= \left|(s-1)\frac{1}{2\pi i}\oint_{|w-1|=\frac{1}{2}}
                \frac{h_n(w)w^{-1/2}-h_n(1)}{w-1}\frac{dw}{w-s}\right|
            \\[1ex]
            & \leq 4|s-1|\sup_{|w-1|=\frac{1}{2}}
                \left|\frac{h_n(w)}{w^{1/2}}-h_n(1)\right| \leq c|s-1|.
    \end{align*}
    Inserting this into (\ref{proof: proposition phin: eq2}) we obtain that
    there exists a constant $C_1>0$ such that
    \begin{equation}\label{property: phinhat}
        |\hat\phi_n(z)-1| \leq C_1|z-1|,
            \qquad\mbox{for all $n\geq n_2$ and $|z-1|\leq 1/4$.}
    \end{equation}
    Therefore, the isolated singularity of $\hat\phi_n$ at 1 is
    removable so that $\hat\phi_n$ is analytic in
    $\mathbb{C}\setminus(-\infty,0]$, and there exists $\delta>0$
    such that $\Re\hat\phi_n(z)>0$, for all $n\geq n_2$ and
    $|z-1|<\delta$. This yields,
    \begin{equation}\label{definition: phin}
        \phi_n(z) \equiv \bigl(\frac{1}{2}h_n(1)\bigr)^{2/3}(z-1)
            \hat\phi_n(z)^{2/3}
    \end{equation}
    is analytic for $z\in U_\delta$.

    Observe that, by (\ref{property: phinhat}) and (\ref{definition: phin}),
    $\phi_n(z)$ is uniformly (in $n$ and $z$) bounded in $U_\delta$.
    This implies, by using Cauchy's theorem for derivatives,
    that $\phi_n''(z)$ is also uniformly (in $n$ and $z$) bounded
    in $U_{\delta}$ for a smaller $\delta$. Since
    $\hat\phi_n(1)=1$, see (\ref{property: phinhat}),
    we have $\phi_n'(1)=(\frac{1}{2}h_n(1))^{2/3}$, so that
    \[
        \left|\phi_n'(z)-\bigl(\frac{1}{2}h_n(1)\bigr)^{2/3}\right|
            =\left|\int_1^z\phi_n''(s)ds\right|\leq
            C_2|z-1|,\qquad\mbox{for all $n\geq n_2$ and $z\in U_{\delta}$,}
    \]
    for some constant $C_2>0$. Therefore, since $h_n(1)>h_0>0$, see
    Proposition \ref{proposition: hn}, there exists $0<\delta_2<\delta$
    such that for all $n\geq n_2$ the $\phi_n$ are injective and hence
    biholomorphic in $U_{\delta_2}$ and such
    that they satisfy part 1 of the proposition.

    Part 2 follows from the first part (for a possible smaller $\delta_2$).
    The last part of the proposition follows from the second part
    and from equations (\ref{definition: phin}) and (\ref{definition: phinhat}).
\end{proof}

\begin{remark}\label{remark: fn}
    For later reference, observe that by (\ref{definition: fn})
    and (\ref{definition: phin}) the biholomorphic maps $f_n$ are given by,
    \begin{equation}\label{remark fn: eq1}
        f_n(z)= c_n n^{2/3}
        (z-1)\hat f_n(z),\qquad\mbox{for $z\in
        U_{\delta_2}$,}
    \end{equation}
    where $\hat f_n=\hat\phi_n^{2/3}$ with $\hat\phi_n$ given by
    (\ref{definition: phinhat}), and where
    $c_n=(\frac{1}{2}h_n(1))^{2/3}$. The constant $c_n$ has, by
    Remark \ref{remark: h}, the following asymptotic behavior,
    \begin{equation}
        c_n=(2m)^{2/3}(1+\bigO(n^{-1/m})), \qquad\mbox{as $n\to\infty$.}
    \end{equation}
    Furthermore, from the proof of the proposition it follows that $\hat f_n$ is analytic
    and uniformly (in $n$ and $z$) bounded in $U_\delta$, for some $\delta>\delta_2$, and
    that $\hat f_n(1)=1$. Therefore, there exists a constant $C>0$ such that
    \begin{equation}
        |\hat f_n(z)-1|=\left|\frac{1}{2\pi i}\oint_{|s-1|=\frac{\delta+\delta_2}{2}}
            \frac{\hat
            f_n(s)-1}{s-1}\frac{ds}{s-z}\right||z-1|\leq C|z-1|,
    \end{equation}
    for all $n\geq n_2$ and $z\in U_{\delta_2}$.
\end{remark}

\medskip

We now have introduced the necessary ingredients to define
$P_n^{(1)}$. Let $n\geq n_2$ and $\sigma\in(\frac{\pi}{3},\pi)$,
and recall that the contour $\Sigma_S$ is not yet defined. We
suppose that $\Sigma_S$ is defined in $U_{\delta_2}$ as the
inverse $f_n$-image of $\gamma_\sigma\cap f_n(U_{\delta_2})$.
Define,
\begin{equation}\label{definition: Pnhat}
    P_n^{(1)}(z)=\Psi(f_n(z)),
        \qquad\mbox{for $z\in U_{\delta_2}\setminus f_n^{-1}(\gamma_\sigma)$.}
\end{equation}
Then, we immediately see that $P_n^{(1)}$ is analytic in
$U_{\delta_2}\setminus\Sigma_S$ with jump relations (\ref{jump
relations Pn1}).

\subsubsection*{Step 3: Determine $E_n$ explicitly}

In this final step, we determine the invertible analytic matrix
valued function $E_n$ in equation (\ref{definition: Pn}) such that
the matching condition (c) of the RH problem for $P_n$ is
satisfied. From (\ref{definition: Pn}), (\ref{definition: Pnhat}),
(\ref{asymptotics: Psi}), (\ref{property fn}) and
(\ref{definition: fn}) we see that (to ensure that the matching
condition is satisfied) we have to define $E_n$, for all $n\geq
n_2$ as,
\begin{equation}\label{definition: En}
    E_n(z)= P^{(\infty)}(z) z^{\frac{1}{2}\alpha\sigma_3}e^{\frac{\pi
    i}{4}\sigma_3}\frac{1}{\sqrt 2}
    \begin{pmatrix}
        1 & -1 \\ 1 & 1
    \end{pmatrix}
    f_n(z)^{\frac{\sigma_3}{4}},\qquad\mbox{for $z\in U_{\delta_2}$.}
\end{equation}
This ends the contruction of the parametrix $P_n$.

\begin{remark}\label{remark: Eninvertible}
    Obviously, $E_n$ is analytic in
    $U_{\delta_2}\setminus(-\infty,1]$. Using condition (b) of the
    RH problem for $P^{(\infty)}$ and using the fact that
    $f_{n,+}(x)^{\sigma_3/4}=f_{n,-}(x)^{\sigma_3/4}e^{\frac{\pi
    i}{2}\sigma_3}$ for $x\in(1-\delta_2,1)$, it is easy to
    check that $E_n$ has no jumps on $(1-\delta_2,1)$. So, what remains
    is a possible isolated singularity at 1. However, $E_n$ has at most
    $1/2$-root singularities at 1, which implies that the singularity
    at 1 has to be removable. Therefore, $E_n$ is indeed analytic in
    $U_{\delta_2}$.

    Further, from (\ref{definition: En}) and from the fact that
    $\det P^{(\infty)}\equiv 1$, see Remark \ref{remark: detPinfinity=1}, it
    follows that $\det E_n\equiv 1$, so that $E_n$ is also invertible.
\end{remark}

\subsubsection*{Summary of the obtained result}

We will now briefly summarize the obtained result. Let $n\geq n_2$
and $\sigma\in(\frac{\pi}{3},\pi)$, and suppose that the contour
$\Sigma_S$ satisfies $f_n(\Sigma_S\cap
U_{\delta_2})=\gamma_\sigma\cap f_n(U_{\delta_2})$. Define,
\begin{equation}\label{definition: Pnbis}
    P_n(z)=E_n(z)\Psi(f_n(z))e^{-n\xi_n(z)\sigma_3}z^{-\frac{1}{2}\alpha\sigma_3},
        \qquad\mbox{for $z\in U_{\delta_2}\setminus f_n^{-1}(\gamma_\sigma)$,}
\end{equation}
where the matrix valued function $E_n$ is given by
(\ref{definition: En}), the matrix valued function $\Psi$ by
(\ref{definition: Psi}), and the scalar function $f_n$ by
(\ref{definition: fn}). Then, $P_n$ solves the RH problem for
$P_n$. Furthermore, using (\ref{definition: Pnbis}),
(\ref{definition: En}), (\ref{asymptotics: Psi}), (\ref{property
fn}), (\ref{definition: fn}) and part 3 of Proposition
\ref{proposition: phin}, we have
\begin{equation}\label{asymptotic expansion: Pn}
    P_n(z)P^{(\infty)}(z)^{-1}\sim
    I+\sum_{k=1}^\infty \Delta_k(z)\frac{1}{n^k},\qquad\mbox{as $n\to\infty$,}
\end{equation}
uniformly for $z$ in compact subsets of $\{0<|z-1|<\delta_2\}$ and
$\sigma$ in compact subsets of $(\frac{\pi}{3},\pi)$, where
$\Delta_k$ is a meromorphic $2\times 2$ matrix valued function
given by,
\begin{multline}\label{definition: Deltak}
    \Delta_k(z)=\frac{1}{2(-\xi_n(z))^k}
            P^{(\infty)}(z) z^{\frac{1}{2}\alpha\sigma_3} \\
            \times\,
            \begin{pmatrix}
                (-1)^k(s_k+t_k) & i(s_k-t_k) \\
                -i(-1)^k(s_k-t_k) & s_k+t_k
            \end{pmatrix}z^{-\frac{1}{2}\alpha\sigma_3} P^{(\infty)}(z)^{-1},
\end{multline}
for $z\in \{0<|z-1|<\delta_2\}$. Here, the coefficients $s_k$ and
$t_k$ are defined by (\ref{definition: sktk}).

\begin{remark}\label{remark: Deltak}
    Obviously, $\Delta_k$ is analytic in $U_{\delta_2}
    \setminus(-\infty,1]$. Using the fact that $\xi_{n,+}(x)
    =-\xi_{n,-}(x)$ for $x\in(1-\delta_2,1)$ together with
    condition (b) of the RH problem for $P^{(\infty)}$, the
    reader can verify that $\Delta_k$ has no jumps on
    $(1-\delta_2,1)$, so that $\Delta_k$ is indeed meromorphic
    in $U_{\delta_2}$.

    Furthermore, since $\xi_n(z)=\bigO((z-1)^{3/2})$ and
    $P^{(\infty)}(z)=\bigO((z-1)^{1/4})$ as $z\to 1$ it follows
    that $\Delta_k$ has a pole of order at most
    $\left[\frac{3k+1}{2}\right]$ at 1.
\end{remark}

\subsection{Parametrix $\widetilde P_n$ near the endpoint 0}
    \label{section: parametrix near 0}

Here, we do the local analysis near 0. We will construct inside
the disk $\tilde U_{\delta_3}=\{z\in\mathbb{C}:|z|<\delta_3\}$
with center 0 and radius $\delta_3>0$ (sufficiently small and
which will be determined as part of the problem in Proposition
\ref{proposition: phintilde} below), a $2\times 2$ matrix valued
function $\widetilde P_n$ that satisfies the following conditions.

\subsubsection*{RH problem for $\widetilde P_n$:}
\begin{itemize}
    \item[(a)] $\widetilde P_n:\tilde U_{\delta_3}\setminus\Sigma_S\to
        \mathbb{C}^{2\times 2}$ is analytic.
    \item[(b)] $\widetilde P_{n,+}(z)=\widetilde
        P_{n,-}(z)v_S(z)$ for $z\in\Sigma_S\cap\tilde U_{\delta_3}$,
        with $v_S$ the jump matrix (\ref{definition: vS}) for $S$.
    \item[(c)] $\widetilde P_n(z)P^{(\infty)}(z)^{-1}=I+\bigO(1/n)$,
        as $n\to\infty$, uniformly for $z\in\partial \tilde
        U_\delta$ and for $\delta$ in compact subsets of $(0,\delta_3)$.
\end{itemize}

The construction of $\widetilde P_n$ is similar to the
construction of the parametrix near the endpoints $\pm 1$ of the
modified Jacobi weight, see \cite{Kuijlaars,KMVV}, and will be
done using Bessel function. It involves, like the construction of
$P_n$, three steps. In the first two steps we construct a matrix
valued function that satisfies conditions (a) and (b) of the RH
problem for $\widetilde P_n$. In the last step we take also the
matching condition (c) into account.

\subsubsection*{Step 1: Transformation to constant jump matrices}

Seek $\widetilde P_n$ in the form,
\begin{equation}\label{Pntilde in Pntilde1}
    \widetilde P_n(z)=
    \widetilde E_n(z)
    \widetilde
    P_n^{(1)}(z)e^{-n\xi_n(z)\sigma_3}(-z)^{-\frac{1}{2}\alpha\sigma_3},
    \qquad\mbox{for $z\in \tilde U_{\delta_3}\setminus \Sigma_S$,}
\end{equation}
with $\widetilde E_n$ an invertible analytic matrix valued
function in $\tilde U_{\delta_3}$, which will be determined in the
third step. Note that by (\ref{property xin: eq2}) the function
$e^{-n\xi_n}$ has no jumps across $(-\infty,0)$. Therefore, if
$\widetilde P_n^{(1)}$ is analytic in $\tilde
U_{\delta_3}\setminus\Sigma_S$ so is $\widetilde P_n$.

It is straightforward to check, using
$\xi_{n,+}(x)+\xi_{n,-}(x)=0$ for $x\in(0,1)$, that if $\widetilde
P_n^{(1)}$ is analytic in $\tilde U_{\delta_3}\setminus \Sigma_S$
with jump relations,
    \begin{equation}\label{jump relations Pn1tilde}
    \widetilde P_{n,+}^{(1)}(z)=
    \begin{cases}
        \widetilde P_{n,-}^{(1)}(z)
        \begin{pmatrix}
            1 & 0 \\
            e^{-\pi i\alpha} & 1
        \end{pmatrix}, & \mbox{for $z\in\Sigma_1\cap \tilde
        U_{\delta_3}$,} \\[3ex]
        \widetilde P_{n,-}^{(1)}(z)
        \begin{pmatrix}
            0 & 1 \\
            -1 & 0
        \end{pmatrix}, & \mbox{for $z\in\Sigma_2\cap \tilde U_{\delta_3}=(0,\delta_3)$,} \\[3ex]
        \widetilde P_{n,-}^{(1)}(z)
        \begin{pmatrix}
            1 & 0 \\
            e^{\pi i\alpha} & 1
        \end{pmatrix}, & \mbox{for $z\in\Sigma_3\cap \tilde
        U_{\delta_3}$,}
    \end{cases}
    \end{equation}
then $\widetilde P_n$ defined by (\ref{Pntilde in Pntilde1})
satisfies conditions (a) and (b) of the RH problem for $\widetilde
P_n$.

\subsubsection*{Step 2: Determine $\widetilde P_n^{(1)}$ explicitly}

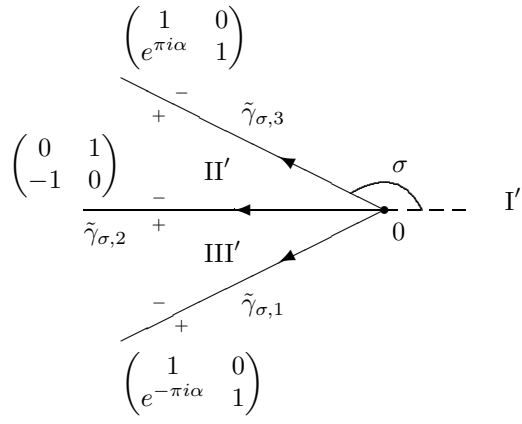
\begin{figure}[t]
\begin{center}
    \setlength{\unitlength}{1mm}
    \begin{picture}(80,58.5)(0,-5)
        \put(31,37){\small $\tilde\gamma_{\sigma,3}$}
        \put(10,21){\small $\tilde\gamma_{\sigma,2}$}
        \put(31,11.5){\small $\tilde\gamma_{\sigma,1}$}
        \put(66,25){\small $\I'$}
        \put(26,29){\small $\II'$}
        \put(26,18){\small $\III'$}
        \put(19,37){\tiny $+$}
        \put(22,40){\tiny $-$}
        \put(22,9){\tiny $+$}
        \put(19,12){\tiny $-$}
        \put(19,22.5){\tiny $+$}
        \put(19,26){\tiny $-$}
        \put(15,47){\small
            $\begin{pmatrix}
                1 & 0 \\
                e^{\pi i\alpha} & 1
            \end{pmatrix}$}
        \put(0,30){\small
            $\begin{pmatrix}
                0 & 1 \\
                -1 & 0
            \end{pmatrix}$}
        \put(15,1){\small
            $\begin{pmatrix}
                1 & 0 \\
                e^{-\pi i\alpha} & 1
            \end{pmatrix}$}
        \put(50,25){\thicklines\circle*{.9}}
        \put(50,25){\line(-2,1){35}} \put(36,32){\thicklines\vector(-2,1){.0001}}
        \put(50,25){\line(-2,-1){35}} \put(36,18){\thicklines\vector(-2,-1){.0001}}
        \put(50,25){\line(-1,0){40}} \put(30,25){\thicklines\vector(-1,0){.0001}}
        \multiput(50,25)(3,0){4}{\line(1,0){1.75}}
        \put(51,21){\small 0}
        \qbezier(55,25)(52,31)(45.523,27.236)
        \put(51,30){\small $\sigma$}
    \end{picture}
    \caption{The oriented contour $\tilde\gamma_\sigma$ and the jump matrix $v_0$ for $\Psi_\alpha$
        on $\tilde \gamma_\sigma$. The
        three straight rays $\tilde\gamma_{\sigma,1},
        \tilde\gamma_{\sigma,2}$ and $\tilde\gamma_{\sigma_3}$ divide the complex plane into
        three regions $\I',\II'$ and $\III'$.}
        \label{figure: contour gammatilde}
\end{center}
\end{figure}

The construction of $\widetilde P_n^{(1)}$ is based upon an
auxiliary RH problem for $\Psi_\alpha$ in the $\zeta$-plane with
jumps on the oriented contour $\tilde\gamma_\sigma$, shown in
Figure \ref{figure: contour gammatilde}, consisting of three
straight rays
\[
    \tilde\gamma_{\sigma,1}:\arg\zeta=-\sigma,
    \qquad\tilde\gamma_{\sigma,2}:\arg\zeta=\pi,
    \qquad\tilde\gamma_{\sigma,3}:\arg\zeta=\sigma,
\]
with $\sigma\in(0,\pi)$. These three rays are oriented to infinity
and divide the complex plane into three regions $\I',\II'$ and
$\III'$, also shown in Figure \ref{figure: contour gammatilde}.
This auxiliary RH problem has been used before
\cite{Kuijlaars,KMVV} in the construction of the parametrix near
the endpoints $\pm 1$ of the modified Jacobi weight, and is the
following.

\medskip

\subsubsection*{RH problem for $\Psi_\alpha$:}
\begin{itemize}
    \item[(a)]
    $\Psi_\alpha:\mathbb{C}\setminus\tilde\gamma_\sigma\to\mathbb{C}^{2\times
    2}$ is analytic.
    \item[(b)]
    $\Psi_{\alpha,+}(\zeta)=\Psi_{\alpha,-}(\zeta)v_0(\zeta)$
    for $\zeta\in\tilde\gamma_\sigma$, where $v_0$ is the piecewise constant matrix valued function
    on $\tilde\gamma_\sigma$ defined as shown in Figure \ref{figure: contour
    gammatilde}, i.e.~$v_0(\zeta)=\left(\begin{smallmatrix} 1 & 0 \\ e^{-\pi i\alpha} & 1\end{smallmatrix}\right)$
    for $\zeta\in\tilde\gamma_{\sigma,1}$, and so on.
    \item[(c)] $\Psi_\alpha$ has the following asymptotic behavior at
    infinity,
    \begin{multline}\label{asymptotics Psialpha}
        \Psi_\alpha(\zeta)\sim (2\pi)^{-\sigma_3/2}
        \zeta^{-\frac{\sigma_3}{4}}\frac{1}{\sqrt 2}
        \begin{pmatrix}
            1 & -i \\
            -i & 1
        \end{pmatrix}\\[1ex]
        \times\left[I
        +\sum_{k=1}^\infty \frac{(\alpha,k-1)}{4^k
            \zeta^{k/2}}
        \begin{pmatrix}
            \frac{(-1)^k}{k}(\alpha^2+\frac{1}{2}k-\frac{1}{4}) & (k-\frac{1}{2})i
            \\[1ex]
            -(-1)^k(k-\frac{1}{2})i & \frac{1}{k}(\alpha^2+\frac{1}{2}k-\frac{1}{4})
        \end{pmatrix}\right]e^{2\zeta^{1/2}\sigma_3},
    \end{multline}
    as $\zeta\to\infty$, uniformly for
    $\zeta\in\mathbb{C}\setminus\tilde\gamma_\sigma$ and $\sigma$
    in compact subsets of $(0,\pi)$. Here, $(\alpha,0)=1$ and
    \begin{equation}\label{definition: alphak}
        (\alpha,k)=\frac{(4\alpha^2-1)(4\alpha^2-9)\cdots(4\alpha^2-(2k-1)^2)}{2^{2k}k!},
            \qquad\mbox{for $k\geq 1$.}
    \end{equation}
\end{itemize}

\medskip

One knows \cite{Kuijlaars,KMVV} that the $2\times 2$ matrix valued
function $\Psi_\alpha=\Psi_\alpha^\sigma$ (we suppress $\sigma$ in
the notation for brevity) defined by,
\begin{equation}
    \label{definition: Psialpha}
    \Psi_\alpha(\zeta)=
    \begin{cases}
        \begin{pmatrix}
            I_\alpha(2\zeta^{\frac{1}{2}}) &
                -\frac{i}{\pi}K_\alpha(2\zeta^{\frac{1}{2}}) \\[1ex]
            -2\pi i\zeta^{\frac{1}{2}}I_\alpha'(2\zeta^{\frac{1}{2}}) &
                -2\zeta^{\frac{1}{2}}K_\alpha'(2\zeta^{\frac{1}{2}})
        \end{pmatrix}, & \mbox{for $\zeta\in\I'$,} \\[5ex]
        \begin{pmatrix}
            \frac{1}{2}H_\alpha^{(1)}(2(-\zeta)^{\frac{1}{2}}) &
                -\frac{1}{2}H_\alpha^{(2)}(2(-\zeta)^{\frac{1}{2}})
                \\[1ex]
            -\pi\zeta^{\frac{1}{2}}(H_\alpha^{(1)})'(2(-\zeta)^{\frac{1}{2}}) &
                \pi\zeta^{\frac{1}{2}}(H_\alpha^{(2)})'(2(-\zeta)^{\frac{1}{2}})
        \end{pmatrix}e^{\frac{1}{2}\alpha\pi
        i\sigma_3}, & \mbox{for $\zeta\in\II'$,} \\[5ex]
        \begin{pmatrix}
            \frac{1}{2}H_\alpha^{(2)}(2(-\zeta)^{\frac{1}{2}}) &
                \frac{1}{2}H_\alpha^{(1)}(2(-\zeta)^{\frac{1}{2}}) \\[1ex]
            \pi\zeta^{\frac{1}{2}}(H_\alpha^{(2)})'(2(-\zeta)^{\frac{1}{2}}) &
                \pi\zeta^{\frac{1}{2}}(H_\alpha^{(1)})'(2(-\zeta)^{\frac{1}{2}})
        \end{pmatrix}e^{-\frac{1}{2}\alpha\pi
        i\sigma_3}, & \mbox{for $\zeta\in\III'$,}
    \end{cases}
\end{equation}
solves the RH problem for $\Psi_\alpha$. Here $I_\alpha$ and
$K_\alpha$ are modified Bessel functions of order $\alpha$, and
$H_\alpha^{(1)}$ and $H_\alpha^{(2)}$ are Hankel functions of
order $\alpha$ of the first and the second kind, respectively. See
for example \cite{AbramowitzStegun} for definitions and properties
of these functions.

\begin{remark}\label{remark: detPsialpha=1}
    From \cite[formulae 9.1.17 and 9.6.15]{AbramowitzStegun} we have
    $\det\Psi_\alpha\equiv 1$.
\end{remark}

\medskip

As in the construction of the parametrix near 1 we construct
$\widetilde P_n^{(1)}$ out of $\Psi_\alpha$ as $\Psi_\alpha(\tilde
f_n(z))$ using appropriate biholomorphic maps $\tilde f_n:\tilde
U_{\delta_3}\to \tilde f_n(\tilde U_{\delta_3})$ with $\tilde
f_n(0)=0$. We will choose them to compensate for the factor
$e^{-n\xi_n(z)\sigma_3}$ in (\ref{Pntilde in Pntilde1}). By
(\ref{asymptotics Psialpha}) we see that a good choice would be to
construct $\tilde f_n$ such that it satisfies,
\begin{equation}\label{property fntilde}
    e^{2\tilde f_n^{1/2}(z)}=(-1)^n e^{n\xi_n(z)},
        \qquad\mbox{for $z\in\tilde U_{\delta_3}\setminus[0,\infty)$.}
\end{equation}
In order to construct $\tilde f_n$ we prove the following lemma
and proposition.

\begin{lemma}\label{lemma: phintildehat}
    Define, for all $n\geq n_2$, the auxiliary scalar function,
    \begin{equation}\label{definition: phintildehat}
        \hat{\tilde\phi}_n(z)= (-z)^{-1/2} \frac{-\pi i}{h_n(0)} \int_0^z\psi_n(s)ds,
    \end{equation}
    with $\psi_n$ given by {\rm (\ref{definition: psin})}. This
    function is, by {\rm (\ref{property: psin})}, well-defined and analytic in
    $\mathbb{C}\setminus[1,\infty)$. Furthermore, there exist
    constants $C_1,\tilde\delta_3>0$ such that
    \begin{align}
        \label{lemma: phintildehat: eq1}
        &\left|\hat{\tilde\phi}_n(z)-1\right|\leq
        C_1|z|,\qquad\mbox{for all $n\geq n_2$ and $|z|\leq 1/4$,}\\[1ex]
        \label{lemma: phintildehat: eq2}
        &\Re\left((-z)^{1/2}\hat{\tilde\phi}_n(z)\right)>0,\qquad
            \mbox{for all $n\geq n_2$ and $z\in\tilde U_{\tilde\delta_3}\setminus[0,\infty)$.}
    \end{align}
\end{lemma}

\begin{proof}
    Equation (\ref{lemma: phintildehat: eq1}) can be proven analogously
    as equation (\ref{property: phinhat}) in the proof of Proposition
    \ref{proposition: phin}.

    By (\ref{lemma: phintildehat: eq1}) there exists $\delta>0$ such that
    $|\arg\hat{\tilde\phi}_n(z)|< \pi/4$, for all $n\geq n_2$ and
    $|z|<\delta$. Since $|\arg(-z)^{1/2}|\leq \pi/4$
    for $|\arg z|\geq \pi/2$, we then have
    \begin{equation}\label{proof lemma phintildehat: eq1}
        \Re \left((-z)^{1/2}\hat{\tilde\phi}_n(z)\right)>0,\qquad
            \mbox{for all $n\geq n_2$ and $|z|<\delta$ with $|\arg z|\geq \frac{\pi}{2}$.}
    \end{equation}
    Now, from (\ref{property: psin}), (\ref{definition: xin}) and the fact that
    $\int_0^1\hat\psi_n(s)ds=1$, we obtain that
    \begin{equation}\label{proof lemma phintildehat: eq2}
        (-z)^{1/2}\hat{\tilde\phi}_n(z)=-\frac{\pi
        i}{h_n(0)}\int_0^z\psi_n(s)ds=
        \begin{cases}
            \frac{1}{h_n(0)}(\xi_n(z)-\pi i), & \mbox{if $\Im z>0$,}
            \\[1ex]
            \frac{1}{h_n(0)}(\xi_n(z)+\pi i), & \mbox{if $\Im z<0$.}
        \end{cases}
    \end{equation}
    So, by (\ref{property xin: eq4}) and since $h_n(0)>0$,
    \begin{equation}\label{proof lemma phintildehat: eq3}
        \Re \left((-z)^{1/2}\hat{\tilde\phi}_n(z)\right)>0,\qquad
            \mbox{for all $n\geq n_2$ and $|z|<\delta_1$ with $0<|\arg z|< \frac{\pi}{2}$.}
    \end{equation}
    with $\delta_1$ defined in (\ref{property xin: eq4}). From
    (\ref{proof lemma phintildehat: eq1}) and (\ref{proof lemma phintildehat: eq3}),
    equation (\ref{lemma: phintildehat: eq2}) is then proven with $\tilde\delta_3=\min\{\delta,\delta_1\}$.
\end{proof}

\begin{proposition}\label{proposition: phintilde}
    There exists $\delta_3>0$ such that for every $n\geq n_2$
    there are biholomorphic maps $\tilde\phi_n:
    \tilde U_{\delta_3}\to\phi_n(\tilde U_{\delta_3})$ satisfying,
    \begin{itemize}
        \item[1.] There exists a constant $c_0>0$ such that for
        all $z\in \tilde U_{\delta_3}$ and all $n\geq n_2$ the derivative
        of $\tilde\phi_n$ can be estimated by:
        $c_0<|\tilde\phi_n'(z)|<1/c_0$ and
        $|\arg\tilde\phi_n'(z)-\pi|<\pi/15$.
        \item[2.] $e^{2n\tilde\phi_n^{1/2}(z)}=(-1)^ne^{n\xi_n(z)}$ for
        $z\in \tilde U_{\delta_3}\setminus[0,\infty)$.
    \end{itemize}
\end{proposition}

\begin{proof}
    Let
    \begin{equation}\label{definition: phintilde}
        \tilde\phi_n(z) \equiv -\frac{1}{4}h_n(0)^2 z \hat{\tilde\phi}_n^2(z)
            = \left(\frac{1}{2}h_n(0)(-z)^{1/2}\hat{\tilde\phi}_n(z)\right)^2,
            \qquad\mbox{for $z\in\mathbb{C}\setminus[1,\infty)$,}
    \end{equation}
    with $\hat{\tilde\phi}_n$ defined by (\ref{definition: phintildehat}).
    From (\ref{lemma: phintildehat: eq1}) and from the fact that $h_n(0)$
    is bounded it follows that $\tilde\phi_n(z)$ is uniformly (in $n$ and
    $z$) bounded in $\{|z|\leq 1/4\}$. So, by using Cauchy's
    theorem for derivatives, $\tilde\phi_n''(z)$ is uniformly
    (in $n$ and $z$) bounded in $\{|z|<\delta\}$ for some
    $0<\delta<1/4$. Since $\hat{\tilde\phi}_n(0)=1$, see (\ref{lemma: phintildehat: eq1}),
    we have $\tilde\phi_n'(0)=-\frac{1}{4}h_n(0)^2$, and thus
    \[
        \left|\tilde\phi_n'(z)+\frac{1}{4}h_n(0)^2\right|=\left|\int_0^z\tilde\phi_n''(s)ds\right|\leq
        C_2|z|,\qquad\mbox{for all $n\geq n_2$ and $|z|<\delta$.}
    \]
    Therefore, since $h_n(0)>h_0>0$, there exists
    $0<\delta_3<\tilde\delta_3$ such that for all $n\geq n_2$ the $\tilde\phi_n$ are injective and hence
    biholomorphic in $\tilde U_{\delta_3}$ and such that they
    satisfy the first part of the proposition.

    The second part of the proposition can be verified by using equations (\ref{lemma: phintildehat:
    eq2}), (\ref{proof lemma phintildehat: eq2}) and (\ref{definition: phintilde}).
\end{proof}

We now define the biholomorphic maps $\tilde f_n$ for all $n\geq
n_2$ as,
\begin{equation}\label{definition: fntilde}
    \tilde f_n(z)=n^2 \tilde \phi_n(z),\qquad \mbox{for $z\in \tilde U_{\delta_3}$.}
\end{equation}
By the second part of the proposition, equation (\ref{property
fntilde}) is then satisfied. Note that, by (\ref{definition:
phintildehat}), (\ref{lemma: phintildehat: eq2}) and
(\ref{definition: phintilde}), equation (\ref{statement of
results: fntilde}) is satisfied.

\begin{remark}\label{remark: fntilde}
    For later reference, we state the analogue of Remark \ref{remark:
    fn}. From (\ref{definition: phintilde}) and (\ref{definition:
    fntilde}), it follows that
    \begin{equation}\label{remark fntilde: eq1}
        \tilde f_n(z)= -\tilde c_n n^2 z
        \hat{\tilde f}_n(z), \qquad\mbox{for $z\in \tilde
        U_{\delta_3}$},
    \end{equation}
    where $\hat{\tilde f}_n=\hat{\tilde\phi}_n^2$ with
    $\hat{\tilde\phi}_n$ given by (\ref{definition: phintildehat}),
    and where $\tilde c_n=(\frac{1}{2}h_n(0))^2$. The constant
    $\tilde c_n$ has, by Remark \ref{remark: h}, the following
    asymptotic behavior,
    \begin{equation}
        \tilde c_n=\left(\frac{2m}{2m-1}\right)^2(1+\bigO(n^{-1/m})),
        \qquad\mbox{as $n\to\infty$.}
    \end{equation}
    Furthermore, from Lemma \ref{lemma: phintildehat}, we have that $\hat{\tilde
    f}_n$ is analytic in $\tilde U_{\delta_3}$ and that $\hat{\tilde f}_n(0)=1$.
    As in Remark \ref{remark: fn}, there exists a constant $C>0$ such that
    \begin{equation}\label{remark fntilde: eq3}
        \left|\hat{\tilde f}_n(z)-1\right|\leq C|z|,
            \qquad\mbox{for all $n\geq n_2$ and $z\in\tilde U_{\delta_3}$.}
    \end{equation}
\end{remark}

\begin{remark}
    Observe that, in contrast to the biholomorphic maps $\phi_n$
    of Proposition \ref{proposition: phin}, the function $\tilde\phi_n$
    maps the upper (lower) part of the disk $\tilde U_{\delta_3}$
    onto the lower (upper) part of $\tilde\phi_n(\tilde U_{\delta_3})$.
\end{remark}

\medskip

We now have all the ingredients to define $\widetilde P_n^{(1)}$.
Let $n\geq n_2$ and $\sigma\in(0,\pi)$, and recall that the
contour $\Sigma_S$ is not yet defined. We suppose that $\Sigma_S$
is defined in $U_{\delta_3}$ as the inverse $\tilde f_n$-image of
$\tilde\gamma_\sigma\cap \tilde f_n(\tilde U_{\delta_3})$. Define,
\begin{equation}\label{definition: Pn1tilde}
    \widetilde P_n^{(1)}(z)=\Psi_\alpha(\tilde f_n(z)),
        \qquad \mbox{for $z\in\tilde U_{\delta_3}\setminus
        \tilde f_n^{-1}(\tilde\gamma_\sigma)$.}
\end{equation}
Then, $\widetilde P_n^{(1)}$ is analytic in $\tilde
U_{\delta_3}\setminus\Sigma_S$ with jump relations (\ref{jump
relations Pn1tilde}).

\subsubsection*{Step 3: Determine $\widetilde E_n$ explicitly}

In this final step, we determine the invertible analytic matrix
valued function $\widetilde E_n$ in equation (\ref{Pntilde in
Pntilde1}) such that the matching condition of the RH problem for
$\widetilde P_n$ is satisfied. From (\ref{Pntilde in Pntilde1}),
(\ref{definition: Pn1tilde}), (\ref{asymptotics Psialpha}),
(\ref{property fntilde}) and (\ref{definition: fntilde}) we have
to define $\widetilde E_n$, for all $n\geq n_2$ as,
\begin{equation}\label{definition: Entilde}
    \widetilde E_n(z)=(-1)^n P^{(\infty)}(z)(-z)^{\frac{1}{2}\alpha\sigma_3}
        \frac{1}{\sqrt 2}
        \begin{pmatrix}
            1 & i \\
            i & 1
        \end{pmatrix}
        \tilde
        f_n(z)^{\frac{\sigma_3}{4}}(2\pi)^{\sigma_3/2},
        \qquad
        \mbox{for $z\in\tilde U_{\delta_3}$.}
\end{equation}
This ends the construction of the parametrix near 0.

\begin{remark}\label{remark: Entildeinvertible}
    Obviously, $\widetilde E_n$ is analytic in $\tilde
    U_{\delta_3}\setminus[0,\infty)$. From condition (b) of the RH
    problem for $P^{(\infty)}$ and from the fact that
    $\tilde f_{n,+}(x)^{\sigma_3/4}=e^{-\frac{\pi i}{2}\sigma_3}\tilde
    f_{n,-}(x)^{\sigma_3/4}$ for $x\in(0,\delta_3)$, one obtains
    that $\widetilde E_n$ has no jumps on $(0,\delta_3)$. So, what
    remains is a possible isolated singularity in 0. However, using
    Remark \ref{remark: detPinfinity=1}, one sees that we have at most
    $1/2$-root singularities in 0, so that the singularity at the
    origin has to be removable. Therefore, $\widetilde E_n$ is
    indeed analytic in $\tilde U_{\delta_3}$.

    From (\ref{definition: Entilde}) and from the fact that
    $\det P^{(\infty)}\equiv 1$, see Remark \ref{remark: detPinfinity=1},
    it follows that $\det \widetilde E_n\equiv 1$, so that $\widetilde E_n$
    is also invertible.
\end{remark}

\subsubsection*{Summary of the obtained result}

We will now briefly summarize the obtained result. Let $n\geq n_2$
and $\sigma\in(0,\pi)$, and suppose that the contour $\Sigma_S$
satisfies $\tilde f_n(\Sigma_S\cap\tilde
U_{\delta_3})=\tilde\gamma_\sigma\cap\tilde f_n(\tilde
U_{\delta_3})$. Define,
\begin{equation}\label{definition: Pntilde}
    \widetilde P_n(z)=\widetilde E_n(z)\Psi_\alpha(\tilde
    f_n(z))e^{-n\xi_n(z)\sigma_3}(-z)^{-\frac{1}{2}\alpha\sigma_3},
    \qquad\mbox{for $z\in\tilde U_{\delta_3}\setminus\tilde
    f_n^{-1}(\tilde\gamma_\sigma)$,}
\end{equation}
where the matrix valued function $\widetilde E_n$ is given by
(\ref{definition: Entilde}), the matrix valued function
$\Psi_\alpha$ by (\ref{definition: Psialpha}), and the scalar
function $\tilde f_n$ by (\ref{definition: fntilde}). Then,
$\widetilde P_n$ solves the RH problem for $\widetilde P_n$.
Furthermore, using (\ref{definition: Pntilde}), (\ref{definition:
Entilde}), (\ref{asymptotics Psialpha}), (\ref{property fntilde})
and (\ref{definition: fntilde}), we have
\begin{equation}\label{asymptotic expansion: Pntilde}
    \widetilde P_n(z) P^{(\infty)}(z)^{-1}\sim
    I+\sum_{k=1}^\infty \widetilde
    \Delta_k(z)\frac{1}{n^k},\qquad\mbox{as $n\to\infty$,}
\end{equation}
uniformly for $z$ in compact subsets of $\{0<|z|<\delta_3\}$ and
for $\sigma$ in compact subsets of $(0,\pi)$, where
$\widetilde\Delta_k$ is a meromorpic $2\times 2$ matrix valued
function given by,
\begin{multline}\label{definition: Deltaktilde}
    \widetilde \Delta_k(z)=\frac{(\alpha,k-1)}{4^k\tilde\phi_n(z)^{k/2}}
        P^{(\infty)}(z)(-z)^{\frac{1}{2}\alpha\sigma_3} \\
            \times
            \begin{pmatrix}
                \frac{(-1)^k}{k}(\alpha^2+\frac{1}{2}k-\frac{1}{4}) & (k-\frac{1}{2})i
                \\[1ex]
                -(-1)^k(k-\frac{1}{2})i & \frac{1}{k}(\alpha^2+\frac{1}{2}k-\frac{1}{4})
            \end{pmatrix}
            (-z)^{-\frac{1}{2}\alpha\sigma_3}P^{(\infty)}(z)^{-1},
\end{multline}
for $z\in\{0<|z|<\delta_3\}$. Here, $(\alpha,0)=1$ and
$(\alpha,k)$ is defined by (\ref{definition: alphak}) for $k\geq
1$. The function $\tilde\phi_n$ is defined by (\ref{definition:
phintilde}) and (\ref{definition: phintildehat}).

\begin{remark}\label{remark: Deltaktilde}
    We have the analogue of Remark \ref{remark: Deltak}. Namely,
    one can check that $\widetilde \Delta_k$ is indeed meromorphic in
    $\tilde U_{\delta_3}$ and has a pole of order at most
    $\left[\frac{k+1}{2}\right]$ at 0.
\end{remark}

\subsection{Final transformation: $S\mapsto R$}
\label{subsection: final transformation}

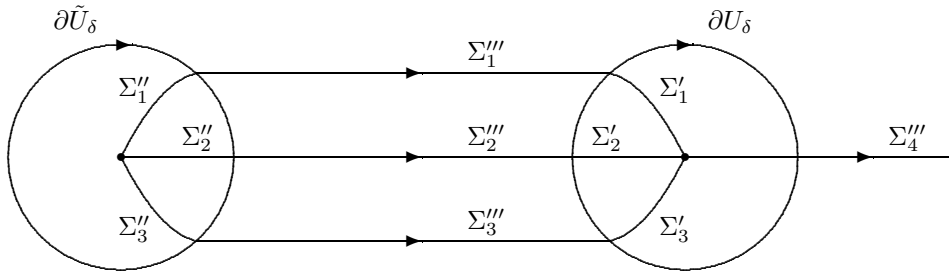
\begin{figure}[t]
\begin{center}
    \setlength{\unitlength}{1mm}
    \begin{picture}(130,41.5)(2.5,15)
        \put(20,35){\line(1,0){110}}
            \put(60,35){\thicklines\vector(1,0){0.0001}}
            \put(120,35){\thicklines\vector(1,0){0.0001}}
            \put(66,36.7){\small $\Sigma_2'''$}
            \put(122,36.7){\small $\Sigma_4'''$}
        \put(20,35){\bigcircle{30}}
            \put(21.5,50){\thicklines\vector(1,0){0.0001}}
            \put(20,35){\circle*{.9}}
            \qbezier(20,35)(25,45)(30,46.18)
            \qbezier(20,35)(25,25)(30,23.82)
            \put(19.6,43.1){\small $\Sigma_1''$}
            \put(28,36.7){\small $\Sigma_2''$}
            \put(19.6,24.6){\small $\Sigma_3''$}
            \put(11,52){\small $\partial \tilde U_\delta$}
        \put(95,35){\bigcircle{30}}
            \put(96,50){\thicklines\vector(1,0){0.0001}}
            \put(95,35){\circle*{.9}}
            \qbezier(95,35)(90,45)(85,46.18)
            \qbezier(95,35)(90,25)(85,23.82)
            \put(91.5,43.1){\small $\Sigma_1'$}
            \put(82.5,36.7){\small $\Sigma_2'$}
            \put(91.5,24.6){\small $\Sigma_3'$}
            \put(98,52){\small $\partial U_\delta$}
        \put(30,46.18){\line(1,0){55}}
            \put(60,46.18){\thicklines\vector(1,0){0.0001}}
            \put(66,47.88){\small $\Sigma_1'''$}
        \put(30,23.82){\line(1,0){55}}
            \put(60,23.82){\thicklines\vector(1,0){0.0001}}
            \put(66,25.52){\small $\Sigma_3'''$}
    \end{picture}
    \caption{The contour $\Sigma_R=\Sigma_S\cup\partial U_\delta\cup\partial \tilde U_\delta$
        depending on the parameters $n,\delta$ and $\nu$.
        Here $\Sigma_S=\bigcup_{j=1}^4\Sigma_j$ with $\Sigma_1=\Sigma_1'\cup\Sigma_1''\cup\Sigma_1'''$
        and so on.}
    \label{figure: contour SigmaR}
\end{center}
\end{figure}

In this subsection, we will perform the final transformation of
our RH problem. Recall that the contour $\Sigma_S$ is still not
yet explicitly defined. We will now define it in terms of the
parameters $n,\delta$ and $\nu$ (a new parameter replacing
$\sigma$). Here, we follow \cite[Section 7.2]{DKMVZ1}.

Let $\delta_0=\min\{\delta_1,\delta_2,\delta_3\}$ (cf.\
(\ref{property xin: eq4}), Proposition \ref{proposition: phin} and
\ref{proposition: phintilde}). Fix $\delta\in(0,\delta_0)$, $n\geq
n_2,$ (cf. Proposition \ref{proposition: hn}) and
$\nu\in(\frac{2\pi}{3},\frac{5\pi}{6})$. From Proposition
\ref{proposition: phin} we know that there exists a
$\sigma=\sigma(n,\nu,\delta)\in(\frac{2\pi}{3}
-\frac{\pi}{15},\frac{5\pi}{6}+\frac{\pi}{15})$ such that
$f_n^{-1}(\gamma_{\sigma,1})\cap\partial U_\delta=\{1+\delta
e^{i\nu}\}$. By the symmetry $\overline{f_n(z)}=f_n(\bar z)$ we
then also have $f_n^{-1}(\gamma_{\sigma,3})\cap\partial
U_\delta=\{1+\delta e^{-i\nu}\}$. We then define $\Sigma_S$ in
$U_\delta$ as the inverse $f_n$-image of $\gamma_\sigma$. We can
do an analogous construction near 0, and define $\Sigma_S$ in the
disk $U_\delta$ near 0 as the inverse $\tilde f_n$-image of
$\tilde\gamma_{\tilde\sigma(n,\nu,\delta)}$ such that $\tilde
f_n^{-1}(\tilde\gamma_{\tilde\sigma,1})\cap \partial \tilde
U_{\delta}=\{-\delta e^{-i\nu}\}$ and $\tilde
f_n^{-1}(\tilde\gamma_{\tilde\sigma,3})\cap \partial \tilde
U_{\delta}=\{-\delta e^{i\nu}\}$.

Further, define a contour $\Sigma_R$ in terms of the contour
$\Sigma_S$. Let $\Sigma_R=\Sigma_S\cup
\partial U_\delta\cup\partial \tilde U_\delta$. This leads to
Figure \ref{figure: contour SigmaR}. Note that the contour
$\Sigma_R$ depends on $n$ (and also on $\delta$ and $\nu$).
However, we immediately see that $\Sigma_1''',\ldots ,\Sigma_4'''$
are independent of $n$.

\medskip

Now, we are ready to do the transformation $S\mapsto R$. Define a
matrix valued function
$R:\mathbb{C}\setminus\Sigma_R\to\mathbb{C}^{2\times 2}$
(depending on the parameters $n,\delta$ and $\nu$) as,
\begin{equation}\label{definition: R}
    R(z)=
    \begin{cases}
        S(z) P_n(z)^{-1},
            & \mbox{for $z\in U_\delta\setminus\Sigma_S$,} \\[1ex]
        S(z) \widetilde P_n(z)^{-1},
            & \mbox{for $z\in \tilde U_\delta\setminus\Sigma_S$,} \\[1ex]
        S(z) P^{(\infty)}(z)^{-1},
            & \mbox{for $z$ elsewhere,}
    \end{cases}
\end{equation}
where $P_n$ is the parametrix near 1 given by (\ref{definition:
Pnbis}), $\widetilde P_n$ is the parametrix near 0, see
(\ref{definition: Pntilde}), $P^{(\infty)}$ is the parametrix for
the outside region given by (\ref{definition: Pinfinity}), and $S$
solves the RH problem for $S$.

\begin{remark}
    The inverses of the matrices $P_n$, $\widetilde P_n$ and
    $P^{(\infty)}$ used in (\ref{definition: R}) exist, since the
    determinants of these matrices are 1. For $P^{(\infty)}$ see Remark \ref{remark:
    detPinfinity=1}, for $P_n$ see (\ref{definition: Pnbis}) and
    Remarks \ref{remark: detPsi=1} and \ref{remark: Eninvertible}, and finally,
    for $\widetilde P_n$ see (\ref{definition: Pntilde}) and
    Remarks \ref{remark: detPsialpha=1} and \ref{remark: Entildeinvertible}.
\end{remark}

By definition, $R$ has jumps on the contour $\Sigma_R$. However,
in the next proposition we will show that $R$ has only jumps on
the reduced contour $\widehat\Sigma_R$, see Figure \ref{figure:
contour SigmaRhat},
\begin{equation}\label{definition: SigmaRhat}
    \widehat\Sigma_R=\Sigma_1'''\cup\Sigma_3'''\cup\Sigma_4'''\cup\partial
    U_\delta\cup\partial \tilde U_\delta.
\end{equation}

\begin{figure}[t]
\begin{center}
    \setlength{\unitlength}{1mm}
    \begin{picture}(130,41.5)(2.5,15)
        \put(110,35){\line(1,0){20}}
            \put(120,35){\thicklines\vector(1,0){0.0001}}
            \put(122,36.7){\small $\Sigma_4'''$}
        \put(20,35){\bigcircle{30}}
            \put(21.5,50){\thicklines\vector(1,0){0.0001}}
            \put(20,35){\circle*{.9}}
            \put(17.5,36.7){\small 0}
            \put(11,52){\small $\partial \tilde U_\delta$}
        \put(95,35){\bigcircle{30}}
            \put(96,50){\thicklines\vector(1,0){0.0001}}
            \put(95,35){\circle*{.9}}
            \put(96,36.7){\small 1}
            \put(98,52){\small $\partial U_\delta$}
        \put(30,46.18){\line(1,0){55}}
            \put(60,46.18){\thicklines\vector(1,0){0.0001}}
            \put(66,47.88){\small $\Sigma_1'''$}
        \put(30,23.82){\line(1,0){55}}
            \put(60,23.82){\thicklines\vector(1,0){0.0001}}
            \put(66,25.52){\small $\Sigma_3'''$}
    \end{picture}
    \caption{The reduced contour $\widehat\Sigma_R$ depending only on
        the parameters $\delta$ and $\nu$ (so independent of $n$).}
    \label{figure: contour SigmaRhat}
\end{center}
\end{figure}
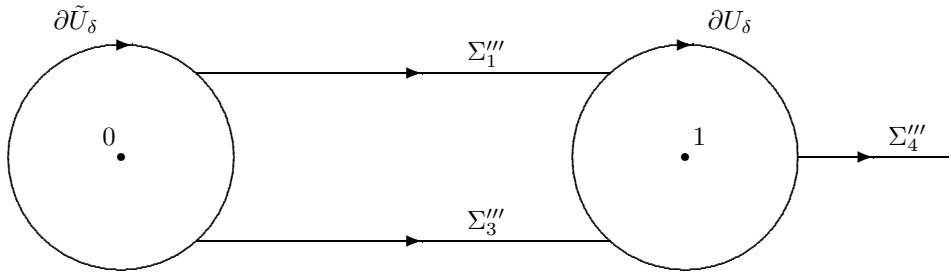

\begin{proposition}
    The matrix valued function $R$ defined by {\rm (\ref{definition: R})} is analytic in
    $\mathbb{C}\setminus\widehat\Sigma_R$, where $\widehat\Sigma_R$ is defined by
    {\rm (\ref{definition: SigmaRhat})}, see also
    Figure {\rm \ref{figure: contour SigmaRhat}}.
\end{proposition}

\begin{proof}
    By construction, the jumps of $S$ and $P^{(\infty)}$ agree on
    $\Sigma_2'''$, the jumps of $S$ and $P_n$ agree on
    $\sum_{j=1}^3\Sigma_j'$, and the jumps of $S$ and $\widetilde
    P_n$ agree on $\sum_{j=1}^3\Sigma_j''$. Therefore, $R$ has an
    analytic continuation to
    $\mathbb{C}\setminus(\widehat\Sigma_R\cup\{0,1\})$. It now
    suffices to show that the isolated singularities of $R$ at 0 and 1
    are removable.

    Since the Airy function is bounded near 0, it follows from
    (\ref{definition: Pn}) that $P_n$ is bounded near 1, and since $\det P_n\equiv 1$
    thus also $P_n^{-1}$. Furthermore, $S$ remains
    also bounded near 1. So, by (\ref{definition: R}), $R$ has a
    removable singularity at 1.

    It remains to prove that the possible isolated singularity of
    $R$ at 0 is removable. We work as in \cite{Kuijlaars,KMVV} using
    the behavior of $R$ near 0. This will be determined by
    multiplying the behavior of $S$ near 0 (which follows from
    (\ref{definition: S}) and condition (d) of the RH problem for $T$) with the behavior of
    $\widetilde P_n^{-1}$ near 0 (which follows from (\ref{definition: Pntilde})
    together with the behavior of $\Psi_\alpha$ near 0 given by
    \cite[equations (6.19)--(6.21)]{KMVV} and the fact that $\det \widetilde P_n\equiv 1$).
    After a straightforward calculation we obtain,
    \begin{align}
        \label{behavior R near 0: eq1}
        R(z) &=
        \bigO\begin{pmatrix}
            z^\alpha & z^\alpha \\
            z^\alpha & z^\alpha
        \end{pmatrix},\qquad\mbox{as $z\to 0$, if $\alpha<0$,}
        \\[1ex]
        \label{behavior R near 0: eq2}
        R(z) &=
        \bigO\begin{pmatrix}
            \log^2 z & \log^2 z \\
            \log^2 z & \log^2 z
        \end{pmatrix},\qquad\mbox{as $z\to 0$, if $\alpha=0$,}
    \end{align}
    and, if $\alpha>0$, that
    \begin{equation}\label{behavior R near 0: eq3}
        R(z)=
        \begin{cases}
            \bigO\begin{pmatrix}
                1 & 1 \\
                1 & 1
            \end{pmatrix}, & \mbox{as $z\to 0$, from outside the
            lens,} \\[3ex]
            \bigO\begin{pmatrix}
                z^{-\alpha} & z^{-\alpha} \\
                z^{-\alpha} & z^{-\alpha}
            \end{pmatrix}, & \mbox{as $z\to 0$, from inside the
            lens.}
        \end{cases}
    \end{equation}
    From this it follows that, in all cases, $R$ has a removable singularity at 0.
    If $\alpha\leq 0$, this is clear from (\ref{behavior R near 0: eq1})
    and (\ref{behavior R near 0: eq2}), since $\alpha>-1$. If $\alpha>0$,
    it follows from (\ref{behavior R near 0: eq3}) that $R$ remains
    bounded if we approach 0 from outside the lens. Therefore,
    $R$ cannot have a pole at 0. Furthermore, we also get from
    (\ref{behavior R near 0: eq3}) that $z^k R(z)$ is bounded near 0 for
    any integer $k>\alpha$. Then, $R$ cannot have an essential singularity at 0 either, so that
    0 is a removable singularity of $R$. This ends the proof of
    the proposition.
\end{proof}

From (\ref{definition: R}) and condition (c) of the RH problems
for $S$ and $P^{(\infty)}$, one then obtains that $R$ is a
solution of the following RH problem on the contour $\widehat
\Sigma_R$.

\subsubsection*{RH problem for $R$}
\begin{itemize}
    \item[(a)] $R:\mathbb{C}\setminus \widehat \Sigma_R\to\mathbb{C}^{2\times
    2}$ is analytic.
    \item[(b)] $R_+(z)=R_-(z)v_R(z)$, for $z\in\widehat \Sigma_R$, with
        \begin{equation}
            v_R(z)=
            \begin{cases}
                P_n(z)P^{(\infty)}(z)^{-1}, & \mbox{for $z\in \partial
                U_\delta$,} \\[1ex]
                \widetilde P_n(z)P^{(\infty)}(z)^{-1}, &
                \mbox{for $z\in\partial \tilde U_\delta$,} \\[1ex]
                P^{(\infty)}(z)v_S(z) P^{(\infty)}(z)^{-1}, & \mbox{for
                $z\in\Sigma_1'''\cup\Sigma_3'''\cup\Sigma_4'''$,}
            \end{cases}
        \end{equation}
        where $v_S$ is the jump matrix of $S$ given by (\ref{definition: vS}).
    \item[(c)] $R(z)=I+\bigO(1/z)$,\qquad as $z\to\infty$.
\end{itemize}

\medskip

Now, we will briefly explain that $R$ defined by (\ref{definition:
R}) is the \textit{unique} solution of the RH problem for $R$ and
that $R$ is uniformly close to the identity matrix as
$n\to\infty$. Again we follow \cite[Section 7.2]{DKMVZ1}.

Introduce the matrix valued function $\Delta_R=v_R-I$ on
$\widehat\Sigma_R$. This matrix satisfies the following estimates
as $n\to\infty$ (with $c_1>0$ some constant and
$\Sigma'''=\Sigma_1'''\cup\Sigma_3'''\cup\Sigma_4'''$), cf.\
\cite[Proposition 7.7]{DKMVZ1},
\begin{align}
    \label{estimates DeltaR: eq1}
        & \|\Delta_R\|_{L^1(\partial U_\delta\cup\partial \tilde U_\delta)}
            + \|\Delta_R\|_{L^2(\partial U_\delta\cup\partial \tilde U_\delta)}
            + \|\Delta_R\|_{L^\infty(\partial U_\delta\cup\partial \tilde U_\delta)}
        =\bigO(1/n), \\[1ex]
    \label{estimates DeltaR: eq2}
        & \|\Delta_R\|_{L^1(\Sigma''')} + \|\Delta_R\|_{L^2(\Sigma''')} +
            \|\Delta_R\|_{L^\infty(\Sigma''')}=\bigO(e^{-nc_1}),
    \end{align}
uniformly for $\delta$ in compact subsets of $(0,\delta_0)$ and for
$\nu\in(\frac{2\pi}{3},\frac{5\pi}{6})$. Here, (\ref{estimates
DeltaR: eq1}) follows from equations (\ref{asymptotic expansion:
Pn}) and (\ref{asymptotic expansion: Pntilde}). Estimate
(\ref{estimates DeltaR: eq2}) follows from equations
(\ref{property xin: eq4}), (\ref{property xin: eq5}) and
(\ref{definition: vS}).

Let $C_-$ be the Cauchy operator on $\widehat\Sigma_R$ given by
$C_-f=(Cf)_-$. The estimates for $\Delta_R$ above imply that the
integral operator $C_{\Delta_R}$ defined as
\[
    C_{\Delta_R}(f)\equiv C_-(f\Delta_R),
        \qquad \mbox{for $f\in L^2(\widehat\Sigma_R,\mathbb{C}^{2\times 2})$,}
\]
is a bounded linear operator from $L^2(\widehat\Sigma_R,\mathbb{C}^{2\times
2})$ into itself with operator norm $\|C_{\Delta_R}\|=\bigO(1/n)$, as
$n\to\infty$. Therefore, $Id-C_{\Delta_R}$ can be inverted by a Neumann series
for $n$ sufficiently large, and we define,
\begin{equation}
    \mu_R \equiv (Id-C_{\Delta_R})^{-1}(C_-\Delta_R)\in
    L^2(\widehat\Sigma_R,\mathbb{C}^{2\times 2}).
\end{equation}
As in \cite[Theorem 7.8]{DKMVZ1} one can then show that for $n$
sufficiently large the RH problem for $R$ has a unique solution,
and that
\begin{equation}
    R=I+C(\Delta_R+\mu_R\Delta_R).
\end{equation}

Note that $\Delta_R$ is exponentially small on $\Sigma'''$.
Therefore, the contribution to $\mu_R$ from $\Sigma'''$ is
exponentially small. Furthermore, from (\ref{asymptotic expansion:
Pn}), (\ref{asymptotic expansion: Pntilde}) and the facts that
$\xi_n$ and $\tilde \phi_n$ can be expanded in powers of
$n^{-1/m}$, it follows that $\Delta_R$ possesses on the disks an
aymptotic expansion in powers of $n^{-1/m}$. This will imply that
$\mu_R$ possesses an asymptotic expansion in powers of $n^{-1/m}$.
Similar as in \cite[Theorem 7.10]{DKMVZ1}, this discussion then
leads to the following theorem.

\begin{theorem}\label{theorem: asymptotics R}
    The matrix valued function $R$ has the following asymptotic
    expansion in powers of $n^{-\frac{1}{m}}$,
    \begin{equation}
        R(z)\sim I+\frac{1}{n}\sum_{k=0}^\infty
            r_k(z)n^{-\frac{k}{m}},\qquad \mbox{as $n\to\infty$,}
    \end{equation}
    uniformly for $\delta$ in compact subsets of $(0,\delta_0)$,
    for $\nu\in(\frac{2\pi}{3},\frac{5\pi}{6})$ and for
    $z\in\mathbb{C}\setminus\widehat\Sigma_R$. Furthermore, the scalar
    functions $r_k$ are bounded functions which are analytic in
    $\mathbb{C}\setminus(\partial U_\delta\cup\partial\tilde U_\delta)$
    and which can be computed explicitly.
\end{theorem}

This theorem states that $R$ is uniformly close to the identity matrix as
$n\to\infty$. By going back in the series of transformations $Y\mapsto U\mapsto
T\mapsto S\mapsto R$ we then find the asymptotics of $Y$. This ends the
asymptotic analysis of the RH problem for $Y$.

\section{Asymptotics of the recurrence coefficients $a_n,b_{n-1}$ and
the leading coefficient $\gamma_n$}

In order to determine the asymptotics (as $n\to\infty$) of
$a_n,b_{n-1}$ and $\gamma_n$ we will make use of the following
result, see for example \cite{Deift,DKMVZ1}. Let $Y$ be the unique
solution of the RH problem for $Y$. There exist $2\times 2$
constant matrices $Y_1$ and $Y_2$ such that
\begin{equation}\label{expansion Y infinity}
    Y(z)
    \begin{pmatrix}
        z^{-n} & 0 \\
        0 & z^n
    \end{pmatrix}=
    I+\frac{Y_1}{z}+\frac{Y_2}{z^2}+\bigO(1/z^3),\qquad\mbox{as
    $z\to\infty$,}
\end{equation}
and
\begin{equation}\label{coefficients in Y}
    a_n=(Y_1)_{11}+\frac{(Y_2)_{12}}{(Y_1)_{12}},\qquad
    b_{n-1}=\sqrt{(Y_1)_{12}(Y_1)_{21}},\qquad \gamma_n=\sqrt{\frac{1}{-2\pi
    i(Y_1)_{12}}}.
\end{equation}

\medskip

We will now rewrite the above expressions for $a_n,b_{n-1}$ and
$\gamma_n$ in terms of the solution of the rescaled RH problem for
$U$. From (\ref{definition: U}) and (\ref{expansion Y infinity})
we have,
\begin{align*}
    U(z)
    \begin{pmatrix}
        z^{-n} & 0 \\
        0 & z^n
    \end{pmatrix}
    &=
        \beta_n^{-(n+\frac{\alpha}{2})\sigma_3}Y(\beta_n z)
        \begin{pmatrix}
            (\beta_n z)^{-n} & 0 \\
            0 & (\beta_n z)^n
        \end{pmatrix}\beta_n^{(n+\frac{\alpha}{2})\sigma_3} \\[1ex]
    &= \beta_n^{-(n+\frac{\alpha}{2})\sigma_3}\left(I+\frac{Y_1}{\beta_n z}+\frac{Y_2}{(\beta_n z)^2}
        +\bigO(1/z^3)\right)\beta_n^{(n+\frac{\alpha}{2})\sigma_3},\qquad
        \mbox{as $z\to\infty$.}
\end{align*}
Therefore, $U$ has an expansion of the form (\ref{expansion Y
infinity}) at infinity. Further, with $U_1$ and $U_2$ the
analogues of $Y_1$ and $Y_2$, respectively, we have,
\[
    Y_1=\beta_n \beta_n^{(n+\frac{\alpha}{2})\sigma_3} U_1
    \beta_n^{-(n+\frac{\alpha}{2})\sigma_3}, \qquad
    Y_2=\beta_n^2 \beta_n^{(n+\frac{\alpha}{2})\sigma_3} U_2
    \beta_n^{-(n+\frac{\alpha}{2})\sigma_3}.
\]
Inserting this into (\ref{coefficients in Y}) we arrive at,
\begin{equation}\label{coefficients in U}
    a_n=\beta_n\left((U_1)_{11}+\frac{(U_2)_{12}}{(U_1)_{12}}\right),\qquad
    b_{n-1}=\beta_n\sqrt{(U_1)_{12}(U_1)_{21}},
\end{equation}
and
\begin{equation}\label{coefficient gamman in U}
    \gamma_n=\beta_n^{-(n+\frac{\alpha}{2}+\frac{1}{2})}\sqrt{\frac{1}{-2\pi i(U_1)_{12}}}.
\end{equation}

We thus need to determine the constant matrices $U_1$ and $U_2$.
For large $|z|$ we have by (\ref{definition: T}),
(\ref{definition: S}) and (\ref{definition: R}),
\begin{equation}\label{U at infinity}
    U(z)=e^{\frac{1}{2}n\ell_n\sigma_3} R(z) P^{(\infty)}(z)
    e^{ng_n(z)\sigma_3}e^{-\frac{1}{2}n\ell_n\sigma_3}.
\end{equation}
So, in order to get $U_1$ and $U_2$ we need the asymptotics of
$P^{(\infty)}(z),e^{ng_n(z)\sigma_3}$ and $R(z)$ as $z\to\infty$.

\subsubsection*{Asymptotics of $P^{(\infty)}(z)$ as $z\to\infty$:}

From (\ref{definition: D}) and (\ref{definition: varphi}) it is easy to check
that that the Szeg\H{o} function $D$ has an asymptotic expansion in powers of
$z^{-1}$ given by $D(z)=2^{-\alpha}\left(1+\frac{1}{4}\alpha z^{-1}+
\frac{1}{32}\alpha(\alpha+3) z^{-2}+\bigO(z^{-3})\right)$, as $z\to\infty$.
Clearly, the scalar function $a$, defined by (\ref{definition: a}), has also an
asymptotic expansion in powers of $z^{-1}$ given by
$a(z)=1-\frac{1}{4}z^{-1}-\frac{3}{32}z^{-2}+\bigO(z^{-3})$, as $z\to\infty$.
Inserting these asymptotics into (\ref{definition: Pinfinity}) we find after a
straightforward calculation,
\begin{equation}\label{asymptotics Pinfinity}
    P^{(\infty)}(z)=I+\frac{P_1^{(\infty)}}{z}+\frac{P_2^{(\infty)}}{z^2}+\bigO(1/z^3),
        \qquad\mbox{as $z\to\infty$,}
\end{equation}
with
\begin{equation}\label{definition: P1P2}
    P^{(\infty)}_1=
        2^{-\alpha\sigma_3}\begin{pmatrix}
            -\frac{\alpha}{4} & \frac{i}{4} \\[1ex]
            -\frac{i}{4} & \frac{\alpha}{4}
        \end{pmatrix}2^{\alpha\sigma_3},\qquad
    P^{(\infty)}_2=
        2^{-\alpha\sigma_3}\begin{pmatrix}
            \frac{\alpha(\alpha-3)+1}{32} & i\frac{(\alpha+2)}{16} \\[1ex]
            i\frac{(\alpha-2)}{16} & \frac{\alpha(\alpha+3)+1}{32}
        \end{pmatrix}2^{\alpha\sigma_3}.
\end{equation}

\subsubsection*{Asymptotics of $e^{ng_n(z)\sigma_3}$ as $z\to\infty$:}

By (\ref{definition: gn}) we have, cf.~\cite[equations (8.8) and
(8.9)]{DKMVZ1},
\begin{equation}\label{asymptotics engn}
    e^{ng_n(z)\sigma_3}
    \begin{pmatrix}
        z^{-n} & 0 \\
        0 & z^n
    \end{pmatrix}=I+\frac{G_1}{z}+\frac{G_2}{z^2}+\bigO(1/z^3),\qquad\mbox{as
    $z\to\infty$,}
\end{equation}
with
\begin{equation}\label{definition: G1G2}
    G_1=
    \begin{pmatrix}
        -n\int_0^1 t\hat\psi_n(t)dt & 0 \\[1ex]
        0 & n\int_0^1 t\hat\psi_n(t)dt
    \end{pmatrix},\qquad G_2=\begin{pmatrix} * & 0 \\ 0 & *
    \end{pmatrix}.
\end{equation}

\subsubsection*{Asymptotics of $R(z)$ as $z\to\infty$:}

Analogous as in \cite{DKMVZ1} the matrix valued function $R$ has
the following asymptotic expansion at infinity,
\begin{equation}\label{asymptotics R}
    R(z)=I+\frac{R_1}{z}+\frac{R_2}{z^2}+\bigO(1/z^3),\qquad\mbox{as
    $z\to\infty$,}
\end{equation}
where the $2\times 2$ constant matrices $R_1$ and $R_2$ satisfy,
\begin{align*}
    R_1&=\left(-\frac{1}{2\pi i}\int_{\partial U_\delta}\Delta_1(y)dy
        -\frac{1}{2\pi i}\int_{\partial \tilde U_\delta}\widetilde\Delta_1(y)dy\right)\frac{1}{n}+\bigO(1/n^2),
    \qquad\mbox{as $n\to\infty$,}
    \\[2ex]
    R_2&=\left(-\frac{1}{2\pi i}\int_{\partial U_\delta}y\Delta_1(y)dy
        -\frac{1}{2\pi i}\int_{\partial \tilde U_\delta}y\widetilde\Delta_1(y)dy\right)\frac{1}{n}+\bigO(1/n^2),
    \qquad\mbox{as $n\to\infty$,}
\end{align*}
with $\Delta_1$ and $\widetilde \Delta_1$ given by
(\ref{definition: Deltak}) and (\ref{definition: Deltaktilde}),
respectively. Recall that $\Delta_1$ is meromorphic in
$\{0<|z-1|<\delta_0\}$ with pole of order 2 at 1, and that
$\widetilde\Delta_1$ is meromorphic in $\{0<|z|<\delta_0\}$ with
pole of order 1 at 0. So, one can evaluate the above integrals
using the residue theorem. After a rather long (but
straightforward) calculation we obtain,
\begin{multline}
    \label{definition: R1}
    R_1=2^{-\alpha\sigma_3}\left[\frac{1}{16 h_n(0)}(4\alpha^2-1)
        \begin{pmatrix}
            1 & i \\
            i & -1
        \end{pmatrix}
        +\frac{h_n'(1)}{16h_n(1)^2}
        \begin{pmatrix}
            1 & -i \\
            -i & -1
        \end{pmatrix}\right. \\[2ex]
        +\left. \frac{1}{48h_n(1)}
        \begin{pmatrix}
            -12\alpha^2+3 & i(12\alpha^2+24\alpha+11) \\
            i(12\alpha^2-24\alpha+11) & 12\alpha^2-3
        \end{pmatrix}\right]
        2^{\alpha\sigma_3}\frac{1}{n}+\bigO(1/n^2),
\end{multline}
 and
\begin{equation}\label{definition: R2}
    (R_2)_{12}=i4^{-\alpha}\left(\frac{12\alpha^2+24\alpha+16}{48h_n(1)}
        -\frac{h_n'(1)}{16h_n(1)^2}\right)\frac{1}{n}+\bigO(1/n^2).
\end{equation}

\begin{remark}\label{remark: asymptotic expansion R1R2}
    The error terms in (\ref{definition: R1}) and
    (\ref{definition: R2}) have an asymptotic expansion
    in powers of $n^{-1/m}$ with coefficients that can be explicitly determined,
    cf.\ \cite[Section 8.1]{DKMVZ1}.
\end{remark}

\medskip

We now have all the ingredients to determine the asymptotics of
$a_n,b_{n-1}$ and $\gamma_n$. Note that by (\ref{U at infinity}),
(\ref{asymptotics Pinfinity}), (\ref{asymptotics engn}) and
(\ref{asymptotics R}) the constant matrices $U_1$ and $U_2$ are
given by
\begin{equation}\label{definition: U1}
    U_1=
    e^{\frac{1}{2}n\ell_n\sigma_3}\left[P^{(\infty)}_1+G_1+R_1\right]e^{-\frac{1}{2}n\ell_n\sigma_3},
\end{equation}
and
\begin{equation}\label{definition: U2}
    U_2=
    e^{\frac{1}{2}n\ell_n\sigma_3}\left[P^{(\infty)}_2+G_2+R_2+R_1P^{(\infty)}_1+
        \left(P_1^{(\infty)}+R_1\right)G_1\right]
        e^{-\frac{1}{2}n\ell_n\sigma_3}.
\end{equation}

\begin{varproof}\textbf{of Theorem \ref{theorem: asymptotics coefficients}.}
    First, we determine the asymptotics of $b_{n-1}$.
    Inserting (\ref{definition: U1}) into (\ref{coefficients in U}),
    and using the fact that $(G_1)_{12}=(G_1)_{21}=0$, see
    (\ref{definition: G1G2}), the recurrence coefficient $b_{n-1}$ is
    given by,
    \[
        b_{n-1} = \beta_n\left[(P_1^{(\infty)})_{12}(P_1^{(\infty)})_{21}+
                (P_1^{(\infty)})_{12}(R_1)_{21}+(R_1)_{12}(P_1^{(\infty)})_{21}+(R_1)_{12}(R_1)_{21}\right]^{1/2}.
    \]
    From equations (\ref{definition: P1P2}) and (\ref{definition: R1}), and from the fact that
    $(R_1)_{12}(R_2)_{21}=\bigO(1/n^2)$ as $n\to\infty$, we then arrive at
    \begin{align}
        \nonumber
        b_{n-1}
                &=\frac{\beta_n}{4}
                    \left[1+4i\left(4^{-\alpha}(R_1)_{21}-4^\alpha(R_1)_{12}\right)+\bigO(1/n^2)\right]^{1/2} \\[1ex]
                &=\frac{\beta_n}{4}\left[1+\frac{4\alpha}{h_n(1)}\frac{1}{n}+\bigO(1/n^2)\right]^{1/2}
                ,\qquad\mbox{as $n\to\infty$.}
    \end{align}
    This proves equation (\ref{theorem: asymptotics coefficients: bn-1}).

    Next, we consider $a_n$. Inserting (\ref{definition: U1}) and
    (\ref{definition: U2}) into (\ref{coefficients in U}),
    and using the facts that $(G_1)_{12}=(G_2)_{12}=0$ and
    $(G_1)_{11}+(G_1)_{22}=0$, see (\ref{definition: G1G2}), the
    recurrence coefficient $a_n$ is given by,
    \begin{equation}\label{proof an: eq1}
        a_n=\beta_n\left((P_1^{(\infty)})_{11}+(R_1)_{11}
            +\frac{(P_2^{(\infty)})_{12}+(R_2)_{12}+(R_1 P_1^{(\infty)})_{12}}
            {(P_1^{(\infty)})_{12}+(R_1)_{12}}\right).
    \end{equation}
    Now, from (\ref{definition: P1P2}) and from the fact that
    $(R_1)_{12}=\bigO(1/n)$ as $n\to\infty$, it follows that
    \[
        \frac{1}{(P_1^{(\infty)})_{12}+(R_1)_{12}}=-4i 4^\alpha
            \left(1+4i 4^\alpha
            (R_1)_{12}+\bigO(1/n^2)\right),\qquad\mbox{as $n\to\infty$.}
    \]
    Inserting this into (\ref{proof an: eq1}), we then find, from
    (\ref{definition: P1P2}), (\ref{definition: R1}) and
    (\ref{definition: R2}), after a straightforward calculation,
    \begin{align*}
        a_n&=\frac{\beta_n}{2}\left(1+4(R_1)_{11}+4i
        4^\alpha(R_1)_{12}-8i4^\alpha(R_2)_{12}+\bigO(1/n^2)\right)\\[1ex]
        &= \frac{\beta_n}{2}\left(1+\frac{2(\alpha+1)}{h_n(1)}\frac{1}{n}+\bigO(1/n^2)\right),
        \qquad\mbox{as $n\to\infty$,}
    \end{align*}
    which proves (\ref{theorem: asymptotics coefficients: an}).

    Finally, we consider the leading coefficients $\gamma_n$.
    Inserting (\ref{definition: U1}) into (\ref{coefficient gamman in U}),
    and using (\ref{definition: P1P2}), (\ref{definition: G1G2})
    and the fact that $(R_1)_{12}=\bigO(1/n)$ as $n\to\infty$, we find,
    \begin{align}\label{asymptotics gamman}
        \nonumber
        \gamma_n
        &= \beta_n^{-(n+\frac{\alpha}{2}+1)}
        e^{-\frac{1}{2}n\ell_n}\sqrt{\frac{2}{\pi}}\,2^\alpha(1-4i 4^\alpha
        (R_1)_{12})^{-1/2} \\[1ex]
        &= \beta_n^{-(n+\frac{\alpha}{2}+\frac{1}{2})}
        e^{-\frac{1}{2}n\ell_n}\sqrt{\frac{2}{\pi
        }}\,2^\alpha\left(1+2i4^\alpha
        (R_1)_{12}+\bigO(1/n^2)\right),\qquad\mbox{as $n\to\infty$.}
    \end{align}
    From (\ref{definition: R1}) we then obtain (\ref{theorem:
    asymptotics coefficients: gamman}).

    Note that the statement that all the error terms have an explicit asymptotic
    expansion in powers of $n^{-1/m}$ follows from Remark \ref{remark: asymptotic expansion R1R2}.
    This ends the proof of Theorem \ref{theorem: asymptotics
    coefficients}.
\end{varproof}

\section{Plancherel-Rotach type asymptotics for the orthonormal polynomials $p_n$}

In this section we will determine the asymptotic behavior (as
$n\to\infty$) of $p_n(\beta_n z)$ in the four asymptotic regions
$A_\delta$, $B_\delta$, $C_\delta$ and $D_\delta$, see Figure
\ref{figure: asymptotic regions}. We will do this by rewriting
$p_n(\beta_n z)$ in terms of the solution $U$ of the rescaled RH
problem for $U$. From (\ref{definition: Y}), (\ref{definition:
U}), (\ref{definition: Vn}) and (\ref{theorem: asymptotics
coefficients: gamman}), we obtain
\begin{equation}\label{asymptotics pn in U}
    p_n(\beta_nz) = (\beta_n z)^{-\frac{\alpha}{2}} e^{\frac{1}{2}Q(\beta_n z)}
    e^{-\frac{1}{2}n\ell_n}2^\alpha\sqrt{\frac{2}{\pi\beta_n}}
    \left(z^{\frac{\alpha}{2}}e^{-\frac{1}{2}nV_n(z)} U_{11}(z)\right)(1+\bigO(1/n)),
\end{equation}
as $n\to\infty$, uniformly for $z\in\mathbb{C}$.

So, we need to determine the asymptotics of $U_{11}$. In order to
do this we will make use of the following properties. Since
$g_n-\xi_n-\frac{1}{2}\ell_n-\frac{1}{2}V_n$ is analytic in
$\mathbb{C}\setminus(-\infty,1]$ and zero on $(1,\infty)$, see
(\ref{property xin: eq3}), it follows from the uniqueness theorem
that
\begin{equation}\label{asymptotics U: prop1}
    g_n(z)-\xi_n(z)-\frac{1}{2}\ell_n=\frac{1}{2}V_n(z),\qquad\mbox{for
    $z\in\mathbb{C}\setminus(-\infty,1]$.}
\end{equation}
The reader can easily verify that
$z^{1/2}+(z-1)^{1/2}=\varphi(z)^{1/2}$ and
$z^{1/2}-(z-1)^{1/2}=\varphi(z)^{-1/2}$ for
$z\in\mathbb{C}\setminus(-\infty,1]$, where $\varphi$ is defined
by (\ref{definition: varphi: eq1}). From (\ref{definition: a}) we
then obtain,
\[
    \frac{a(z)+a(z)^{-1}}{2}=\frac{\varphi(z)^{1/2}}{2z^{1/4}(z-1)^{1/4}},
    \qquad
    \frac{a(z)-a(z)^{-1}}{-2i}=-i\frac{\varphi(z)^{-1/2}}{2z^{1/4}(z-1)^{1/4}},
\]
for $z\in\mathbb{C}\setminus[0,1]$ by analytic continuation.
Inserting this into (\ref{definition: Pinfinity}) and using the
fact that $\varphi(z)=e^{i\arccos(2z-1)}$ for $z\in\mathbb C_+$,
we have
\begin{align}\label{asymptotics U: prop2}
    \nonumber
    P^{(\infty)}(z)&=\frac{2^{-\alpha\sigma_3}}{2z^{1/4}(z-1)^{1/4}}
    \begin{pmatrix}
        \varphi(z)^{\frac{1}{2}(\alpha+1)} &
            i\varphi(z)^{-\frac{1}{2}(\alpha+1)} \\[1ex]
        -i\varphi(z)^{\frac{1}{2}(\alpha-1)} &
            \varphi(z)^{-\frac{1}{2}(\alpha-1)}
    \end{pmatrix}z^{-\frac{1}{2}\alpha\sigma_3} \\[2ex]
    &= \frac{2^{-\alpha\sigma_3}}{2z^{1/4}(z-1)^{1/4}}
    \begin{pmatrix}
        e^{i\eta_1(z)} & ie^{-i\eta_1(z)} \\[1ex]
        -ie^{i\eta_2(z)} & e^{-i\eta_2(z)}
    \end{pmatrix}z^{-\frac{1}{2}\alpha\sigma_3}, \qquad
    \mbox{for $z\in\mathbb C_+$,}
\end{align}
where we have introduced (for the sake of brevity) the notation,
\[
    \eta_1(z)=\frac{1}{2}(\alpha+1)\arccos(2z-1),\qquad
    \eta_2(z)=\frac{1}{2}(\alpha-1)\arccos(2z-1).
\]

\subsection{The outside region $A_\delta$}

In this subsection we determine the asymptotics of the orthonormal
polynomials $p_n$ in the outside region $A_\delta$.

\begin{varproof}\textbf{of Theorem \ref{theorem: asymptotics pn} (a).}
    For $z\in A_\delta$ we have from (\ref{definition: T}), (\ref{definition:
    S}), (\ref{definition: R}), from the first equality in (\ref{asymptotics U:
    prop2}), and from (\ref{asymptotics U: prop1}),
    \begin{equation}
        U(z)\begin{pmatrix}
            1 \\
            0
        \end{pmatrix}= \frac{1}{2z^{1/4}(z-1)^{1/4}}z^{-\frac{\alpha}{2}} e^{\frac{1}{2}nV_n(z)}
            e^{\frac{1}{2}n\ell_n\sigma_3}R(z)
            2^{-\alpha\sigma_3}
            \begin{pmatrix}
                \varphi(z)^{\frac{1}{2}(\alpha+1)} \\
                -i\varphi(z)^{\frac{1}{2}(\alpha-1)}
            \end{pmatrix} e^{n\xi_n(z)}.
    \end{equation}
    So, from Theorem \ref{theorem: asymptotics R}
    and from the fact that $1/\varphi$ is
    uniformly bounded in $A_\delta$, we then obtain
    \begin{align*}
        U_{11}(z)&=z^{-\frac{\alpha}{2}}e^{\frac{1}{2}n V_n(z)}e^{\frac{1}{2}n\ell_n}2^{-\alpha}
                \frac{\varphi(z)^{\frac{1}{2}(\alpha+1)}}{2z^{1/4}(z-1)^{1/4}}
                e^{n\xi_n(z)}
            \left(R_{11}(z)-i4^\alpha R_{12}(z)\varphi(z)^{-1}\right) \\[1ex]
            &= z^{-\frac{\alpha}{2}}e^{\frac{1}{2}n V_n(z)}e^{\frac{1}{2}n\ell_n}2^{-\alpha}
                \frac{\varphi(z)^{\frac{1}{2}(\alpha+1)}}{2z^{1/4}(z-1)^{1/4}}
                e^{n\xi_n(z)}\left(1+\bigO(1/n)\right),
    \end{align*}
    as $n\to\infty$, uniformly for $\delta$ in compact subsets of $(0,\delta_0)$ and for $z\in A_\delta$.
    Inserting this into (\ref{asymptotics pn in U}) and using the definition of $\xi_n$, see
    (\ref{definition: xin}), the first part of the theorem is proven.
\end{varproof}

\subsection{The bulk region $B_\delta$}

Here, we will determine the asymptotics of $p_n$ in the bulk
region $B_\delta$. We will make use of the following proposition,
which will also be used (in the next section) to prove the
universality result in the bulk of the spectrum, cf.~Theorem
\ref{theorem: universality} (a).

\begin{proposition}\label{proposition: Bdelta}
    For $z\in B_\delta$, the first column of $U$ satisfies,
    \begin{equation}\label{proposition: Bdelta: eq1}
        U(z)\begin{pmatrix}1 \\ 0 \end{pmatrix}=
        z^{-\frac{\alpha}{2}}e^{\frac{1}{2}nV_n(z)}e^{\frac{1}{2}n\ell_n\sigma_3}R(z)
        P^{(\infty)}(z)z^{\frac{\alpha}{2}\sigma_3}
            \begin{pmatrix}
                e^{n\xi_n(z)} \\
                e^{-n\xi_n(z)}
            \end{pmatrix},
    \end{equation}
    where $\xi_n$ is given by {\rm(\ref{definition: xin})},
    $P^{(\infty)}$ is the parametrix for the outside region given
    by {\rm(\ref{definition: Pinfinity})}, and $R$
    is the solution of the RH problem for $R$ given by {\rm (\ref{definition: R})}.
\end{proposition}

\begin{proof}
    The proposition can easily be verified using (\ref{definition:
    T}), (\ref{definition: S}), (\ref{definition: R}) and (\ref{asymptotics U: prop1}).
\end{proof}

\begin{varproof}\textbf{of Theorem \ref{theorem: asymptotics pn} (b).}
    Inserting the expression (\ref{asymptotics U: prop2}) for $P^{(\infty)}$
    into equation (\ref{proposition: Bdelta: eq1}), and using the fact
    that $(z-1)^{1/4}=(1-z)^{1/4}e^{\frac{\pi i}{4}}$ for $z\in B_\delta$,
    and the fact that $\xi_n(z)=-\pi i\int_1^z\psi_n(s)ds$, we arrive at,
    \begin{multline}
        U(z)\begin{pmatrix}
            1 \\
            0
        \end{pmatrix}
        = \frac{1}{z^{1/4}(1-z)^{1/4}} z^{-\frac{\alpha}{2}}
        e^{\frac{1}{2}nV_n(z)}e^{\frac{1}{2}n\ell_n\sigma_3} \\
        \times\, R(z)
        2^{-\alpha\sigma_3}
        \begin{pmatrix}
            \cos\left(\eta_1(z)-\pi n\int_1^z\psi_n(s)ds-\frac{\pi}{4}\right)
            \\[1ex]
            -i\cos\left(\eta_2(z)-\pi n\int_1^z\psi_n(s)ds-\frac{\pi}{4}\right)
        \end{pmatrix}.
    \end{multline}
    The second part of the theorem now follows
    easily from this equation together with (\ref{asymptotics pn in U})
    and Theorem \ref{theorem: asymptotics R}.
\end{varproof}

\subsection{The Airy region $C_\delta$}

In order to determine the asymptotics of the orthonormal
polynomials in the Airy region $C_\delta$ we start, as in the
previous subsection, with a proposition which will also be used to
prove the universality result at the soft edge of the spectrum,
cf.~Theorem \ref{theorem: universality} (b).

\begin{proposition}\label{proposition: Cdelta}
    For $z\in U_\delta$, the first column of $U$ satisfies,
    \begin{equation}\label{proposition: Cdelta: eq1}
        U(z)\begin{pmatrix}
            1 \\
            0
        \end{pmatrix}=\sqrt{2\pi}e^{-\frac{\pi i}{4}}z^{-\frac{\alpha}{2}}
        e^{\frac{1}{2}nV_n(z)}e^{\frac{1}{2}n\ell_n\sigma_3}R(z)E_n(z)
        \begin{pmatrix}
            \Ai(f_n(z)) \\
            \Ai'(f_n(z))
        \end{pmatrix},
    \end{equation}
    where $f_n$ is the scalar function given by {\rm (\ref{definition: fn})},
    $E_n$ is the $2\times 2$ matrix valued function given by
    {\rm (\ref{definition: En})}, and $R$ is the solution of the RH problem for $R$
    given by {\rm (\ref{definition: R})}.
\end{proposition}

\begin{proof}
    Let $z\in U_\delta$ be such that $f_n(z)\in\II$, cf.\ Figure \ref{figure: contour gamma}.
    This means that $z$ lies inside
    the disk $U_\delta$ and belongs to the upper part of
    the lens. Then, from (\ref{definition: T}), (\ref{definition: S}), (\ref{definition: R}),
    (\ref{definition: Pnbis}) and (\ref{asymptotics U: prop1}) we
    find after a straightforward calculation,
    \begin{equation}\label{proof proposition Cdelta: eq1}
        U(z)\begin{pmatrix}1 \\ 0\end{pmatrix}=z^{-\frac{\alpha}{2}}e^{\frac{1}{2}nV_n(z)}
        e^{\frac{1}{2}n\ell_n\sigma_3}R(z)
        E_n(z)\Psi(f_n(z))
        \begin{pmatrix}
            1 \\
            1
        \end{pmatrix},
    \end{equation}
    where the $2\times 2$ matrix valued function $\Psi$ is given by (\ref{definition: Psi}).
    Now, since $f_n(z)\in\II$, we have by (\ref{definition: Psi}),
    \[
        \Psi(f_n(z))\begin{pmatrix} 1 \\ 1 \end{pmatrix}
        =\sqrt{2\pi} e^{-\frac{\pi i}{4}}
        \begin{pmatrix}
            \Ai(f_n(z)) \\
            \Ai'(f_n(z))
        \end{pmatrix}.
    \]
    Inserting this into (\ref{proof proposition Cdelta: eq1}), equation (\ref{proposition: Cdelta: eq1}) is
    proven in this sector of $U_\delta$. The proof in the other sectors (i.e.\ $z\in U_\delta$
    such that $f_n(z)\in \I\cup\III\cup\IV$) is analogous.
\end{proof}

\begin{varproof}\textbf{of Theorem \ref{theorem: asymptotics pn} (c).}
    Inserting expression (\ref{asymptotics U: prop2})
    for $P^{(\infty)}$ into the definition (\ref{definition: En}) of $E_n$, we
    find
    \[
        \sqrt{2\pi}e^{-\frac{\pi i}{4}}E_n(z)=
        \frac{\sqrt\pi}{z^{1/4}(z-1)^{1/4}}2^{-\alpha\sigma_3}
        \begin{pmatrix}
            \cos\eta_1(z) & -i\sin\eta_1(z) \\
            -i\cos\eta_2(z) & -\sin\eta_2(z)
        \end{pmatrix}f_n(z)^{\frac{\sigma_3}{4}},
    \]
    for $z\in C_\delta$.
    Inserting this in turn into (\ref{proposition: Cdelta: eq1}) we arrive at,
    \begin{multline}
        U(z)\begin{pmatrix}
            1 \\
            0
        \end{pmatrix}
        = \frac{\sqrt\pi}{z^{1/4}(z-1)^{1/4}} z^{-\frac{\alpha}{2}}e^{\frac{1}{2}nV_n(z)}
            e^{\frac{1}{2}n\ell_n\sigma_3} \\
            \,\times R(z) 2^{-\alpha\sigma_3}
            \begin{pmatrix}
            \cos\eta_1(z) & -i\sin\eta_1(z) \\
            -i\cos\eta_2(z) & -\sin\eta_2(z)
        \end{pmatrix}f_n(z)^{\frac{\sigma_3}{4}}
        \begin{pmatrix}
            \Ai(f_n(z)) \\
            \Ai'(f_n(z))
        \end{pmatrix},
    \end{multline}
    for $z\in C_\delta$.
    Now, one can easily check that
    the scalar functions
    $\cos\eta_2/\cos\eta_1$ and
    $\sin\eta_2/\sin\eta_1$ are bounded in
    $C_\delta$, so that by Theorem \ref{theorem: asymptotics R},
    \begin{multline}
        U_{11}(z)
        =
        z^{-\frac{\alpha}{2}}e^{\frac{1}{2}nV_n(z)}e^{\frac{1}{2}n\ell_n}
        2^{-\alpha}\frac{\sqrt\pi}{z^{1/4}(z-1)^{1/4}} \\
        \times\, \begin{pmatrix}
            \cos\eta_1(z)(1+\bigO(1/n)) & -i\sin\eta_1(z)(1+\bigO(1/n))
        \end{pmatrix}f_n(z)^{\frac{\sigma_3}{4}}
        \begin{pmatrix}
            \Ai(f_n(z)) \\
            \Ai'(f_n(z))
        \end{pmatrix},
    \end{multline}
    as $n\to\infty$, uniformly for $\delta$ in compact subsets of
    $(0,\delta_0)$ and for $z\in C_\delta$.
    Inserting this into (\ref{asymptotics pn in U}) part (c) of the theorem is proven.
\end{varproof}

\subsection{The Bessel region $D_\delta$}

We will start with a result similar to Propositions
\ref{proposition: Bdelta} and \ref{proposition: Cdelta}. This
result will also be used to prove the universality result at the
hard edge of the spectrum, cf.\ Theorem \ref{theorem:
universality} (c) and Theorem \ref{theorem: averages}.

\begin{proposition}\label{proposition: Ddelta}
    For $z\in \tilde U_\delta$, the first column of $U$ satisfies,
    \begin{equation}\label{proposition: Ddelta: eq1}
        U(z)\begin{pmatrix} 1 \\ 0 \end{pmatrix}=
        z^{-\frac{\alpha}{2}}e^{\frac{1}{2}nV_n(z)}
        e^{\frac{1}{2}n\ell_n\sigma_3}R(z)\widetilde
        E_n(z)\begin{pmatrix}
            J_\alpha(2(-\tilde f_n(z))^{1/2}) \\[1ex]
            -2\pi i(-\tilde f_n(z))^{1/2} J_\alpha'(2(-\tilde f_n(z))^{1/2})
        \end{pmatrix},
    \end{equation}
    where $\tilde f_n$ is the scalar function given by {\rm(\ref{definition:
    fntilde})}, $\widetilde E_n$ is the $2\times 2$ matrix valued
    function given by {\rm (\ref{definition: Entilde})}, and $R$ is
    the solution of the RH problem for $R$ given by {\rm (\ref{definition: R})}.
\end{proposition}

\begin{proof}
    Let $z\in\tilde U_\delta$ be such that $\tilde f_n(z)\in \III'$,
    cf.\ Figure \ref{figure: contour gammatilde}. This means that $z$ lies inside the disk
    $\tilde U_\delta$ and belongs to the upper part of the lens. Then,
    from (\ref{definition: T}), (\ref{definition: S}), (\ref{definition:
    R}), (\ref{definition: Pntilde}), (\ref{asymptotics U: prop1}),
    and from the fact that
    $(-z)^{-\frac{1}{2}\alpha\sigma_3}=z^{-\frac{1}{2}\alpha\sigma_3}e^{\frac{1}{2}\alpha\pi i\sigma_3}$
    for $z\in \mathbb{C}_+$,
    we find after a straightforward calculation,
    \begin{equation}\label{proof: proposition Ddelta: eq1}
        U(z)\begin{pmatrix} 1 \\ 0 \end{pmatrix}=
        z^{-\frac{\alpha}{2}}e^{\frac{1}{2}nV_n(z)}
        e^{\frac{1}{2}n\ell_n\sigma_3}R(z)\widetilde
        E_n(z)\Psi_\alpha\left(\tilde f_n(z)\right)
        e^{\frac{1}{2}\pi i\alpha\sigma_3}
        \begin{pmatrix}
            1 \\ 1
        \end{pmatrix},
    \end{equation}
    where the $2\times 2$ matrix valued function $\Psi_\alpha$
    is given by (\ref{definition: Psialpha}). Now, since $\tilde
    f_n(z)\in \III'$, we have by (\ref{definition: Psialpha}) and
    \cite[formulae 9.1.3 and 9.1.4]{AbramowitzStegun}
    \[
        \Psi_\alpha\left(\tilde f_n(z)\right)
        e^{\frac{1}{2}\pi i\alpha\sigma_3}
        \begin{pmatrix}
            1 \\ 1
        \end{pmatrix}=
        \begin{pmatrix}
            J_\alpha(2(-\tilde f_n(z))^{1/2}) \\[1ex]
            2\pi \tilde f_n(z)^{1/2} J_\alpha'(2(-\tilde f_n(z))^{1/2})
        \end{pmatrix}.
    \]
    Inserting this into (\ref{proof: proposition Ddelta: eq1}) and using the fact that
    $\tilde f_n(z)^{1/2}=-i(\tilde f_n(z))^{1/2}$ for $z\in\mathbb{C}_+$,
    the proposition is proven in this sector
    of $U_\delta$. The proof in the other sectors (i.e. $z\in\tilde U_\delta$ such
    that $\tilde f_n(z)\in\I'\cup\II'$) is analogous. For the case $\tilde f_n(z)\in \I'$
    one has to use \cite[formulae 9.1.35 and 9.6.3]{AbramowitzStegun}.
\end{proof}

\begin{varproof}\textbf{of Theorem \ref{theorem: asymptotics pn}
(d).} Inserting (\ref{asymptotics U: prop2}) into the definition
(\ref{definition: Entilde}) of $\widetilde E_n$, and using the
facts that for $z\in D_\delta$,
\[
    (-z)^{\frac{1}{2}\alpha\sigma_3}=z^{\frac{1}{2}\alpha\sigma_3}e^{-\frac{1}{2}\pi\alpha\sigma_3},
    \quad (z-1)^{1/4}=(1-z)^{1/4} e^{\frac{\pi i}{4}},
    \quad \tilde f_n(z)^{1/4}=(-\tilde f_n(z))^{1/4}e^{-\frac{\pi i}{4}},
\]
we obtain
\begin{equation}\label{proof asymptotics Ddelta: Entilde}
    \widetilde E_n(z)= \frac{(-1)^n 2^{-\alpha\sigma_3}}{\sqrt 2
    z^{1/4}(1-z)^{1/4}}
    \begin{pmatrix}
        \sin\zeta_1(z) &  i\cos\zeta_1(z) \\
        -i\sin\zeta_2(z) & \cos\zeta_2(z)
    \end{pmatrix}(-\tilde f_n(z))^{\frac{\sigma_3}{4}}(2\pi)^{\frac{\sigma_3}{2}},
\end{equation}
for $z\in D_\delta$. Here the scalar functions $\zeta_1$ and
$\zeta_2$ are given by $\zeta_{1,2}(z)=\eta_{1,2}(z)-\pi
\alpha/2$. Plugging this in into (\ref{proposition: Ddelta: eq1}),
we arrive at,
\begin{multline}\label{proof asymptotics Ddelta: eq1}
    U(z)\begin{pmatrix}
        1 \\
        0
    \end{pmatrix}=
    (-1)^n \frac{\sqrt\pi(-\tilde
    f_n(z))^{1/4}}{z^{1/4}(1-z)^{1/4}}
    z^{-\frac{\alpha}{2}}e^{\frac{1}{2}nV_n(z)}e^{\frac{1}{2}n\ell_n\sigma_3} \\
    \times\,R(z)
    2^{-\alpha\sigma_3}
    \begin{pmatrix}
        \sin\zeta_1(z) &  \cos\zeta_1(z) \\
        -i\sin\zeta_2(z) & -i\cos\zeta_2(z)
    \end{pmatrix}
    \begin{pmatrix}
            J_\alpha(2(-\tilde f_n(z))^{1/2}) \\[1ex]
            J_\alpha'(2(-\tilde f_n(z))^{1/2})
        \end{pmatrix},
\end{multline}
for $z\in D_\delta$. Now, since the functions
$\sin\zeta_2/\sin\zeta_1$ and $\cos\zeta_2/\cos\zeta_1$ are
bounded in $D_\delta$, we then have by Theorem \ref{theorem:
asymptotics R},
\begin{multline}
    U_{11}(z)=z^{-\frac{\alpha}{2}}e^{\frac{1}{2}nV_n(z)}e^{\frac{1}{2}n\ell_n}2^{-\alpha}
    (-1)^n\frac{\sqrt{\pi}(-\tilde f_n(z))^{1/4}}{z^{1/4}(1-z)^{1/4}}
    \\
    \times\,\begin{pmatrix}
        \sin\zeta_1(z)(1+\bigO(1/n)) &  \cos\zeta_1(z)(1+\bigO(1/n))
    \end{pmatrix}
    \begin{pmatrix}
            J_\alpha(2(-\tilde f_n(z))^{1/2}) \\[1ex]
            J_\alpha'(2(-\tilde f_n(z))^{1/2})
        \end{pmatrix},
\end{multline}
as $n\to\infty$, uniformly for $\delta$ in compact subsets of
$(0,\delta_0)$ and for $z\in D_\delta$. The last part of the
theorem then follows from (\ref{asymptotics pn in U}).
\end{varproof}

\medskip

Proposition \ref{proposition: Ddelta} gives the behavior of the
first column of $U$ near the origin. To get the universal behavior
of the three kernels $W_{\I,n}$, $W_{\II,n}$ and  $W_{\III,n}$ at
the the hard edge, see Theorem \ref{theorem: averages}, we also
need the behavior of the second column of $U$ near the origin,
cf.\ Remark \ref{remark: kernels in Y}. This will be given by the
next proposition.

\begin{proposition}\label{proposition: second column U}
    For $z\in\tilde U_\delta$, the second column of $U$ satisfies,
    \begin{multline}\label{proposition: second column U: eq1}
        U(z)
        \begin{pmatrix}
            0 \\ 1
        \end{pmatrix}
        = z^{\frac{\alpha}{2}}e^{-\frac{1}{2}nV_n(z)} e^{\frac{1}{2}n\ell_n\sigma_3}
        R(z)\widetilde E_n(z) \\
        \times\,
        \begin{cases}
            \begin{pmatrix}
                \frac{1}{2} H_\alpha^{(1)}(2(-\tilde f_n(z))^{1/2}) \\[1ex]
                -\pi i (-\tilde f_n(z))^{1/2} (H_\alpha^{(1)})'(2(-\tilde f_n(z))^{1/2})
            \end{pmatrix}, &\mbox{if $z\in\tilde
            U_\delta\cap\mathbb{C}_+$,} \\[4ex]
            \begin{pmatrix}
                -\frac{1}{2} H_\alpha^{(2)}(2(-\tilde f_n(z))^{1/2}) \\[1ex]
                \pi i (-\tilde f_n(z))^{1/2} (H_\alpha^{(2)})'(2(-\tilde f_n(z))^{1/2})
        \end{pmatrix}, &\mbox{if $z\in\tilde
            U_\delta\cap\mathbb{C}_-$.}
        \end{cases}
    \end{multline}
\end{proposition}

\begin{proof}
    The proof is analogous to the proof of Proposition
    \ref{proposition: Ddelta}.
\end{proof}

\begin{remark}
    In work in progress \cite{DGKV} we need to evaluate $U(0)$ for the case $\alpha>0$. From
    the fact that $-\tilde f_n(z)=\tilde c_n n^2 z (1+\bigO(z))$ as $z\to 0$, and from
    (\ref{proposition: Ddelta: eq1}), (\ref{proposition: second column U: eq1}),
    (\ref{definition: Vn}) and
    \cite[formulae 9.1.7, 9.1.9 and 9.1.27]{AbramowitzStegun}, we obtain,
    \[
        U(0)= e^{\frac{1}{2}n\ell_n\sigma_3}R(0)\widetilde E_n(0)
        \begin{pmatrix}
            \frac{1}{\alpha} & -\frac{i}{2\pi} \\[1ex]
            -\pi i & \frac{\alpha}{2}
        \end{pmatrix}
        \left(\frac{\tilde c_n^{\frac{\alpha}{2}}n^\alpha}{\Gamma(\alpha)}\right)^{\sigma_3}
        e^{\frac{1}{2}Q(0)\sigma_3}.
    \]
    So, it remains to evaluate $\widetilde E_n(0)$. This can be
    done using (\ref{proof asymptotics Ddelta: Entilde}) and the
    facts $\sin\zeta_1(z)=1+\bigO(z)$,
    $\sin\zeta_2(z)=-1+\bigO(z)$, $\cos\zeta_1(z)=(\alpha+1)\sqrt
    z(1+\bigO(z))$ and $\cos\zeta_2(z)=(1-\alpha)\sqrt z(1+\bigO(z))$ as
    $z\to 0$. After a straightforward calculation we then arrive
    at,
    \begin{multline}
        U(0)=(-1)^n e^{\frac{1}{2}n\ell_n\sigma_3}R(0)2^{-\alpha\sigma_3}
        \begin{pmatrix}
            1 & i(\alpha+1) \\
            i & 1-\alpha
        \end{pmatrix}(\tilde c_n n^2)^{\frac{1}{4}\sigma_3}
        \\
        \times\,
        \begin{pmatrix}
            \frac{1}{\alpha} & -\frac{i}{2} \\[1ex]
            -\frac{i}{2} & \frac{\alpha}{4}
        \end{pmatrix}
        \left(
        \frac{\tilde c_n^{\frac{\alpha}{2}} n^\alpha e^{\frac{1}{2}Q(0)} \sqrt\pi}
            {\Gamma(\alpha)}\right)^{\sigma_3}.
    \end{multline}
    Note that $\det U(0)=1$ as it should be.
\end{remark}

\section{Applications in random matrix theory}

In this final section we will prove Theorems \ref{theorem:
universality} and \ref{theorem: averages}. We will only consider
the proof of Theorem \ref{theorem: universality} (c), i.e.\ the
universality of $K_n$ at the hard edge of the spectrum, in
detail. The other results can be proven analogously, and are left
as an easy exercise for the reader.

\medskip

We start by writing the kernel $K_n$ in terms of the first column
of $U$. Recall that $K_n$ is already written in terms of the first
column of $Y$ by (\ref{Kn in Y}). Using (\ref{definition: U}) we
then obtain,
\begin{equation}\label{Kn in U}
    \beta_n K_n(\beta_n x,\beta_n y)
        = x^{\frac{\alpha}{2}}e^{-\frac{1}{2}nV_n(x)}
            y^{\frac{\alpha}{2}}e^{-\frac{1}{2}nV_n(y)}
            \frac{1}{2\pi i(x-y)}
            \begin{pmatrix}
                0 & 1
            \end{pmatrix} U^{-1}(y) U(x)
            \begin{pmatrix}
                1 \\ 0
            \end{pmatrix}.
\end{equation}
Then, we prove Theorem \ref{theorem: universality} (c) by
inserting the expression for the first column of $U$ derived in
Proposition \ref{proposition: Ddelta} into this expression for
$K_n$. In order to prove part (a) and (b) of the theorem we have
to use Propositions \ref{proposition: Bdelta} and
\ref{proposition: Cdelta}, respectively.

\medskip

Introduce, for the sake of brevity, some notation. With $\tilde
c_n=\left(\frac{1}{2}h_n(0)\right)^2$, cf.\ Remark \ref{remark:
fntilde}, define
\[
    u_n=\frac{u}{4\tilde c_n n^2},\quad
    v_n=\frac{v}{4\tilde c_n n^2},\quad
    \tilde u_n=2(-\tilde f_n(u_n))^{1/2},\quad\mbox{and}\quad
    \tilde v_n=2(-\tilde f_n(v_n))^{1/2}.
\]
Note that by (\ref{remark fntilde: eq1}) and (\ref{remark fntilde:
eq3}),
\begin{equation}\label{untildevntilde}
    \tilde u_n=\sqrt u(1+\bigO(u/n))
        \qquad\mbox{and}\qquad \tilde v_n=\sqrt v(1+\bigO(v/n)),
        \qquad\mbox{as $n\to\infty$,}
\end{equation}
uniformly for $u,v$ in bounded subsets of
$(0,\infty)$. Further, we also need the following lemma.
\begin{lemma}
    Let $L_n = \widetilde E_n R$, where $\widetilde E_n$ is
    given by {\rm (\ref{definition: Entilde})}, and where $R$
    is the solution of the RH problem for $R$. Then, with the
    notation of above,
    \begin{equation}\label{lemma: Ln: eq1}
        L_n^{-1}(v_n)L_n(u_n)=I+\left(\frac{u-v}{n}\right),
            \qquad\mbox{as $n\to\infty$,}
    \end{equation}
    uniformly for $u,v$ in bounded subsets of $(0,\infty)$.
\end{lemma}

\begin{proof}
    Let $u,v\in (0,\infty)$ and let $\gamma$ be a closed contour
    in $\tilde U_\delta$ encircling the origin once in the
    positive direction. Then, since $L_n$ is analytic in
    $\tilde U_\delta$ we have by Cauchy's formula,
    \[
        L_n^{-1}(v_n)L_n(u_n)=I+\frac{1}{2\pi i}
        \oint_\gamma
        \frac{L_n^{-1}(v_n)L_n(\zeta)}{(\zeta-u_n)(\zeta-v_n)}d\zeta
        (u_n-v_n),
    \]
    for $u,v$ in bounded subsets of $(0,\infty)$ and $n$
    sufficiently large. Now, from (\ref{definition: Entilde}) and
    from the fact that $R$ is uniformly (in $n$ and $z$) bounded in $\tilde U_\delta$, we
    have that $L_n(\zeta)=\bigO(n^{1/2})$, as $n\to\infty$,
    uniformly for $\zeta\in\gamma$. Furthermore, since $\det L_n=1$ we obtain that
    $L_n^{-1}(v_n)=\bigO(n^{1/2})$, as $n\to\infty$, uniformly for
    $v$ in bounded subsets of $(0,\infty)$. Together with
    $u_n-v_n=\bigO(\frac{u-v}{n^2})$, this proves the lemma.
\end{proof}

Now, we have introduced the necessary ingredients to prove the
universal behavior of the kernel $K_n$ at the hard edge of the
spectrum.

\begin{varproof}\textbf{of Theorem \ref{theorem: universality} (c).}
    With the notation of above, we have by
    (\ref{proposition: Ddelta: eq1}),
    \begin{equation}\label{proof theorem universality hard edge: eq1}
        U(u_n)
        \begin{pmatrix}
            1 \\
            0
        \end{pmatrix}
        = u_n ^{-\frac{\alpha}{2}} e^{\frac{1}{2}nV_n(u_n)}
        e^{\frac{1}{2}n\ell_n\sigma_3} L_n(u_n)
        \begin{pmatrix}
            J_\alpha(\tilde u_n) \\
            -\pi i \tilde u_n J_\alpha'(\tilde u_n)
        \end{pmatrix}.
    \end{equation}
    Furthermore, since $A^t=\left(\begin{smallmatrix}0 & 1 \\ -1 & 0\end{smallmatrix}\right)
    A^{-1}\left(\begin{smallmatrix}0 & -1 \\ 1 &
    0\end{smallmatrix}\right)$ for any $2\times 2$ matrix $A$ with
    $\det A=1$, we obtain (note that $\det U\equiv\det L_n\equiv
    1$),
    \begin{align}
        \nonumber
        \begin{pmatrix}
            0 & 1
        \end{pmatrix}
        U^{-1}(v_n)
        &=
            \left[
            \begin{pmatrix}
                0 & -1 \\
                1 & 0
            \end{pmatrix}
            U(v_n)
            \begin{pmatrix}
                1 \\
                0
            \end{pmatrix}\right]^t
        \\[2ex]
        \label{proof theorem universality hard edge: eq2}
        &= v_n ^{-\frac{\alpha}{2}} e^{\frac{1}{2}nV_n(v_n)}
        \begin{pmatrix}
            \pi i \tilde v_n J_\alpha'(\tilde v_n) & J_\alpha(\tilde v_n)
        \end{pmatrix}L_n^{-1}(v_n)e^{-\frac{1}{2}n\ell_n\sigma_3}.
    \end{align}
    Inserting (\ref{proof theorem universality hard edge: eq1}) and
    (\ref{proof theorem universality hard edge: eq2}) into
    (\ref{Kn in U}), we find from (\ref{lemma: Ln: eq1}) and from the facts
    that $J_\alpha(\tilde u_n)=\bigO(u^{\alpha/2})=\tilde u_n J_\alpha'(\tilde
    u_n)$ and $J_\alpha(\tilde v_n)=\bigO(v^{\alpha/2})=\tilde v_n J_\alpha'(\tilde
    v_n)$ as $n\to\infty$,
    \begin{equation}\label{proof theorem universality hard edge: eq3}
        \frac{\beta_n}{4\tilde c_n n^2} K_n(\beta_n u_n,\beta_n v_n)
        =
            \frac{1}{2(u-v)}
            \begin{pmatrix}
                \tilde v_n J_\alpha'(\tilde v_n)
                    & J_\alpha(\tilde v_n)
            \end{pmatrix}
            \begin{pmatrix}
                J_\alpha(\tilde u_n) \\
                -\tilde u_n J_\alpha'(\tilde u_n)
            \end{pmatrix}+\bigO\left(\frac{u^{\frac{\alpha}{2}}v^{\frac{\alpha}{2}}}{n}\right),
    \end{equation}
    as $n\to\infty$, uniformly for $u,v$ in bounded subsets of $(0,\infty)$.
    We can now replace $\tilde u_n$ and $\tilde v_n$ by $\sqrt u$
    and $\sqrt v$, respectively, see (\ref{untildevntilde}), and we then arrive at the Bessel kernel.
    However, then we make an error
    which could not be estimated uniformly for $u-v$ close to zero.
    So, we need to be more careful, and we refer the interested reader
    to the presentation in \cite[Proof of Theorem 1.1 (c)]{KV1} where this problem was solved
    (note that (\ref{proof theorem universality hard edge: eq3})
    corresponds to \cite[equation (3.31)]{KV1}). Then, part (c) of the theorem is
    proven.
\end{varproof}

\begin{varproof}\textbf{of Theorem \ref{theorem: universality} (a) and
(b).} The proof of the first two parts of the theorem is analogous
to the proof of part (c). For part (a) we have to use
(\ref{proposition: Bdelta: eq1}), and for part (b) we use
(\ref{proposition: Cdelta: eq1}). Further, the scaling with $c_n
n^{2/3}$ at the soft edge of the spectrum is clear from
(\ref{remark fn: eq1}).
\end{varproof}

\begin{varproof}\textbf{of Theorem \ref{theorem: averages}.}
    The proof is analogous to the proof of Theorem \ref{theorem: universality}.
    We first have to write (like we have done in (\ref{Kn in U}) with $K_n$) the kernels
    $W_{\I,n}$, $W_{\II,n}$ and $W_{\III,n}$ in terms of $U$. From
    (\ref{Wn in Y}) and (\ref{definition: U}) we have
    \begin{multline}
        \begin{pmatrix}
            \beta_n W_{\II,n}(\beta_n v,\beta_n u) & \beta_n W_{\III,n}(\beta_n u,\beta_n v) \\
            -\beta_n W_{\I,n}(\beta_n u,\beta_n v) & -\beta_n W_{\II,n}(\beta_n u,\beta_n v)
        \end{pmatrix}\\[1ex]
        =\frac{1}{-2\pi i\gamma_{n-1}^2 (u-v)}\beta_n^{\frac{1}{2}\alpha\sigma_3}U^{-1}(v)U(u)
            \beta_n^{-\frac{1}{2}\alpha\sigma_3}.
    \end{multline}
    Then, we prove the theorem by inserting the expressions for the first and second
    column of $U$ (derived in Propositions \ref{proposition:
    Ddelta} and \ref{proposition: second column U}, respectively)
    into this expression. This is left as an easy exercise for the
    reader.
\end{varproof}

\begin{remark}
    Since $W_{\I,n}$ depends only on the first column of
    $U$, it is clear (by Proposition \ref{proposition: Ddelta}) that
    $W_{\I,n}$ depends only on the $J$-Bessel functions. Further, since
    $W_{\II,n}$ depends on both the first and second column of $U$,
    this kernel depends (by Propositions \ref{proposition: Ddelta} and
    \ref{proposition: second column U}) on both the $J$-Bessel and
    Hankel functions. And finally, since $W_{\III,n}$ depends only
    on the second column of $U$, it is clear (by Proposition
    \ref{proposition: second column U}) that this kernel depends
    only on the Hankel functions.
\end{remark}

\section*{Acknowledgements}

The author is grateful to Thomas Kriecherbauer for careful reading
and for many useful discussions and remarks. The author would also
like to thank the Department of Mathematics of the Ruhr
Universit\"at Bochum for hospitality.

\end{document}